\documentclass{amsart}
\usepackage{amsmath}
\usepackage{amssymb}
\usepackage{amsthm}
\usepackage{xcolor}
\usepackage{mathtools}
\usepackage{graphicx}
\usepackage{enumerate}
\usepackage[ruled]{algorithm2e}
\graphicspath{{figures/}}

\DeclareMathOperator*{\sgn}{\textup{sgn}}

\DeclareMathOperator*{\dif}{d\!}

\newtheorem{thm}{Theorem}

\newtheorem{lemma}[thm]{Lemma}
\newtheorem{prop}[thm]{Proposition}
\newtheorem{coro}[thm]{Corollary}

\author[Chung]{Francis Chung}
\address{Department of Mathematics, University of Kentucky, Lexington, KY 40506, USA}
\curraddr{}
\email{fj.chung@uky.edu}

\author[Yang]{Tianyu Yang}
\address{Department of Computational Mathematics Science and Engineering, Michigan State University, East Lansing, MI 48824, USA}
\curraddr{}
\email{yangti27@msu.edu}

\author[Yang]{Yang Yang}
\address{Department of Computational Mathematics Science and Engineering, Michigan State University, East Lansing, MI 48824, USA}
\curraddr{}
\email{yangy5@msu.edu}

\thanks{The research of F. Chung is partially supported by Simons Collaboration Grant 582020.
The research of T. Yang and Y. Yang is partially supported by the NSF grant DMS-1715178, DMS-2006881 and the start-up fund from Michigan State University.
}
 
\title[UMBLT with a Single Optical Measurement]{Ultrasound Modulated Bioluminescence Tomography with a Single Optical Measurement}

\begin{document}

\maketitle

\begin{abstract}
Ultrasound modulated bioluminescence tomography (UMBLT) is an imaging method which can be formulated as a hybrid inverse source problem. In the regime where light propagation is modeled by a radiative transfer equation, previous approaches to this problem require large numbers of optical measurements ~\cite{bal2016ultrasound}. Here we propose an alternative solution for this inverse problem which requires only a single optical measurement in order to reconstruct the isotropic source. Specifically, we derive two inversion formulae based on Neumann series and Fredholm theory respectively, and prove their convergence under sufficient conditions. The resulting numerical algorithms are implemented and experimented to reconstruct both continuous and discontinuous sources in the presence of noise.
\end{abstract}

\section{Introduction}

\textit{BioLuminescence Tomography (BLT)} is a technology that uses light emitted by optical probes to report activity at the molecular level. It has experienced rapid development in the past few decades due to its non-invasiveness and high optical contrast ~\cite{ConBac, NtzRipWanWei}.
However, BLT often suffers from low spatial resolution. This is because of the inherent ill-posedness of the inverse problem in BLT, where reconstruction of the internal distribution of bioluminescent molecules has to be implemented from data measured on the surface -- see ~\cite{BalTam, SteUhl} for more on this problem. 

An effective approach to enhance the spatial resolution of BLT is by ultrasound modulation. This leads to the hybrid imaging modality known as \textit{Ultrasound Modulated BioLumnescence Tomography (UMBLT)} \cite{HuyHayZhaMor, bal2010inverse, bal2014ultrasound,bal2016ultrasound}. In UMBLT, typical BLT is performed while the optical properties of the object-of-interest undergoes a series of perturbation caused by acoustic vibrations. The inverse problem in UMBLT is to recover the spatial distribution of the optical probes from the perturbed boundary measurement of the emitted light. It turns out, as elucidated below, that the perturbed measurement allows retrieval of an internal functional, which helps mitigate the ill-posedness of the inverse problem and enhance the spatial resolution.

This basic idea, in which ultrasound modulation helps improve an otherwise ill-posed problem, has received quite a bit of recent theoretical attention in a number of different contexts. An early example is the problem of ultrasound modulated electrical impedance tomography discussed in ~\cite{BalUMEIT}. In the context of optical tomography, in which one seeks to reconstruct coefficients instead of sources, see for example ~\cite{AmmBosGarNguSep, AmmGarNguSep, AmmNguSep, BalMos, ChuSch, ChuLaiLi}. Other related optical problems include fluorescent ultrasound modulated optical tomography (~\cite{li2019hybrid, li2020inverse}) and multifrequency acousto-optic tomography (~\cite{ChuHosSchDiff, ChuHosSchRTE}). Ultrasound modulated hybrid problems are also part of a broader group of hybrid inverse problems in which the interactions of multiple imaging modalities create well-posed problems; for a survey of such ideas, see ~\cite{BalSurvey}.



We turn to the mathematical formulation of the inverse problem in UMBLT. Let $X$ be a bounded open subset in $\mathbb{R}^n$ with smooth boundary $\partial X$, $n\geq 2$. We model the propagation of light in the medium using the standard \textit{Radiative Transfer Equation (RTE)}:
\begin{equation}\label{eqn:RTE}
\theta\cdot\nabla u+\sigma(x)u-\int_{S^{n-1}}k(x,\theta,\theta')u(x,\theta')\dif\theta' = S(x).
\end{equation}
Here $u=u(x,\theta)$ represents the intensity of light at the point $x\in X$ in the direction $\theta\in \mathbb{S}^{n-1}$,
$S(x)$ is an isotropic source that is independent of $\theta$, $\sigma$ is the \textit{attenuation coefficient} and $k$ is the \textit{scattering kernel}.
Let $\Gamma_+$ and $\Gamma_-$ be the outgoing boundary and the incoming boundary respectively, that is,
\begin{equation}
\Gamma_\pm \coloneqq  \left\{(x,\theta)\in X\times \mathbb{S}^{n-1}\mid\pm\theta\cdot n>0\right\}
\end{equation}
where $n(x)$ is the unit outer normal vector at $x$. 
Assume no light flows through the boundary so that the intensity $u$ obeys the boundary condition
\begin{equation} \label{eq:bc}
u|_{\Gamma_-}=0.
\end{equation}

\bigskip
Next, we take the effect of acoustic modulation into account. Suppose the incident acoustic wave is of the form $\cos(q\cdot x + \varphi)$ where $q$ is the wave vector and $\varphi$ is the phase. The time scale of the acoustic field propagation is generally much greater than that of the optical field, hence
the acoustic field can effectively modulate the time independent RTE. In the presence of the acoustic modulation, the optical coefficients $\sigma$, $k$, and the source $S$ become $\sigma_\varepsilon$, $k_\varepsilon$, and $S_\varepsilon$, respectively. Following~\cite{bal2010inverse, bal2014ultrasound}, the effect of the acoustic modulation on the optical properties can be modeled as
\begin{align}
\sigma_\varepsilon(x) & \coloneqq (1+\varepsilon\cos(q\cdot x+\varphi))\sigma(x) \label{eq:msigma}\\
k_\varepsilon(x,\theta,\theta') & \coloneqq (1+\varepsilon\cos(q\cdot x+\varphi))k(x,\theta,\theta') \label{eq:mk}\\
S_\varepsilon(x,\theta) & \coloneqq (1+\varepsilon\cos(q\cdot x+\varphi))S(x,\theta) \label{eq:mS}
\end{align}
where $0<\varepsilon\leq 1$ is the dimensionless amplitude of the pressure wave. The modulated RTE and boundary condition take the form
\begin{align}
\theta\cdot\nabla u_\varepsilon+\sigma_\varepsilon(x)u_\varepsilon-\int_{\mathbb{S}^{n-1}}k_\varepsilon(x,\theta,\theta')u_\varepsilon(x,\theta')\dif\theta' & = S_\varepsilon(x) \label{eqn:mRTE} \\
u_{\varepsilon}|_{\Gamma_-} & = 0 \label{eq:mbc}
\end{align}
where $u_\varepsilon$ is the modulated RTE solution. We henceforth write $u_0$ for the RTE solution without modulation $(\varepsilon=0)$, that is, $u_0 = u$, the solution of the system ~\eqref{eqn:RTE}, \eqref{eq:bc}.

Under suitable assumptions on $\sigma$, $k$ and $S$ (see (A1)(A2) below), the boundary value problem~\eqref{eqn:mRTE} \eqref{eq:mbc} admits a unique solution $u_\varepsilon$. The measurement is the operator $\Lambda^\varepsilon_S:\mathbb{R}^n \times \{0,\frac{\pi}{2}\} \rightarrow C(\Gamma_+)$ 
defined as
\begin{equation} \label{eq:Lambda}
\Lambda^\varepsilon_S (q,\varphi) \coloneqq  u_\varepsilon|_{\Gamma_+}, \quad\quad\quad \varepsilon \geq 0.
\end{equation}
This is the light that flows out through the boundary during various acoustic modulation

The inverse problem in UMBLT is to reconstruct the non-modulated source term $S(x)$ from the operator $\Lambda^\varepsilon_S$, provided the non-modulated optical coefficients $\sigma$ and $k$ are a-priori known.

\textbf{Our contribution:} The inverse problem of UMBLT was first studied in the special case of the diffusion approximation to the RTE~\cite{bal2014ultrasound}. It was shown that the source can be reconstructed with Lipschitz-type stability. The problem with full RTE model was later considered in~\cite{bal2016ultrasound}, where uniqueness and stability results were established. These results are constructive and valid for general anisotropic sources. Nevertheless the reconstruction algorithm in ~\cite{bal2016ultrasound} has serious drawbacks, stemming from the fact that it requires a point-by-point reconstruction in the interior of the domain.  The reconstruction at each point requires a separate boundary integral to be calculated from boundary observations, which means in practice, $u(x,\theta)$ and $u_{\varepsilon}(x,\theta)$ need to be known at each $(x,\theta) \in \Gamma_-$. This means precision relies on a huge volume of observations, and all of these observations need to be angularly resolved, which can be difficult to guarantee in practice. Moreover the precise boundary integral required for the reconstruction at a given point needs to be calculated by a process described in the proof of Theorem 1.3 of ~\cite{bal2016ultrasound}. This requires repeated calculations of separate solutions to the RTE for \emph{every individual point} of the domain, which in practice is extremely computationally demanding. 

In contrast, the main result of the present paper is that under reasonable conditions (see Theorems \ref{thm:Neumann} and \ref{thm:Fredholm} for precise statements), we can reconstruct an isotropic source from the knowledge of any single boundary integral of the form 
\[
\int_{\Gamma_+} (u - u_{\varepsilon}) v \, dS,
\]
where $u$ is the solution to the RTE \eqref{eqn:RTE} and $v$ is any uniformly positive continuous function on $\Gamma_+$.  This eliminates the requirement for multiple angularly resolved measurements on the boundary. Moreover it drastically reduces the computational demand, as demonstrated below. 


The paper is structured as follows. In Section 2 we describe the derivation of an internal functional from the boundary data collected in the ultrasound modulated experiment, following the ideas of ~\cite{bal2016ultrasound}. In Section 3, we present and discuss the main results, in which two inversion formulae are proved, one based on the Neumann series and the other on Fredholm inversion theory. Each of these formulae allows recovery of the isotropic source from measurement of a single boundary integral.  Finally in Section 4, we describe numerical algorithms for implementing the ideas of Section 3, and present numerical results. 

\section{Derivation of the Internal Functional}

Throughout the paper, we make the following assumptions to ensure well-posedness of some forward boundary value problems. 
\begin{enumerate}
\item[(A1):] $\sigma,k$ and $S$ are continuous on $\overline{X}$ 
\item[(A2):] Set $\rho \coloneqq  \left\|\int_{\mathbb{S}^{n-1}}k(x,\theta,\theta')\dif\theta'\right\|_{L^\infty(X\times \mathbb{S}^{n-1})}$, one of the following inequalities holds:\begin{equation}\label{eqn:X1}
\left( \inf\limits_{x\in \overline{X}}\sigma \right)-\rho\geq\alpha
\end{equation}
where $\alpha >0$ is a positive constant, or
\begin{equation}\label{eqn:X2}
\textup{diam}(X) \rho<1
\end{equation}
where $\textup{diam}(X)\coloneqq \sup\{|x-y|: \; x, y \in X\}$ is the diameter of $X$. 
\end{enumerate}

In order to derive the internal functional, we make the additional assumption that $k$ is invariant under rotations, so 
\begin{equation}\label{eq:isotropy}
k(x,\theta, \theta') = k(x, \theta \cdot \theta').
\end{equation}
This ensures that the integral operator appearing in the RTE is self adjoint over $X \times S^{n-1}$.

Under the assumptions (A1) (A2), well-posedness of the RTE with a prescribed continuous incoming boundary condition and an anisotropic source is proved in~\cite[Theorem 2.1]{bal2016ultrasound}. We will apply it to the special case~\eqref{eqn:RTE} where the source $S=S(x)$ is isotropic. In order to state the result, we define the norm $\|u\|_{L^p(\mathbb{S}^{n-1},C(X))}$ ($1\leq p\leq\infty$)
\[
\left\|u\right\|_{L^p(\mathbb{S}^{n-1},C(X))}\coloneqq\left(\int_{\mathbb{S}^{n-1}}\|u(x,\theta)\|_{C(X)}^p\dif\theta\right)^{\frac1p}
\]
and the function space $L^p(\mathbb{S}^{n-1},C(X))$ by
\[
L^p(\mathbb{S}^{n-1},C(X))\coloneqq \left\{
u: \left\|u\right\|_{L^p(\mathbb{S}^{n-1},C(X))}
<\infty\right\}.
\]

\begin{prop}[{\cite[Theorem 2.1]{bal2016ultrasound}}]\label{thm:WP}
Suppose the assumptions (A1)(A2) hold. Then for any $f_-\in C(\Gamma_-)$, the RTE~\eqref{eqn:RTE} has a unique solution $u\in L^p(\mathbb{S}^{n-1},C(X))$ ($1\leq p\leq\infty$) with the boundary condition $u|_{\Gamma_-}=f_-$.
%
Moreover, if \eqref{eqn:X1} holds, we have the estimate
\[\|u\|_{L^p(\mathbb{S}^{n-1},C(X))}\leq\frac1\alpha\left((\rho+\alpha)\|f_-\|_{L^p(\mathbb{S}^{n-1},C(\partial X))}+\textup{Vol}(\mathbb{S}^{n-1})^{\frac1p}\|S\|_{C(X)}\right).\]
If instead \eqref{eqn:X2} holds, we have the estimate
\[\|u\|_{L^p(\mathbb{S}^{n-1},C(X))}\leq\frac1{1-\tau\rho}\left(\|f_-\|_{L^p(\mathbb{S}^{n-1},C(\partial X))}+\textup{diam}(X)\textup{Vol}(\mathbb{S}^{n-1})^{\frac1p}\|S\|_{C(X)}\right).\]
\end{prop}
Here $\textup{Vol}(\mathbb{S}^{n-1})$ denotes the volume of $\mathbb{S}^{n-1}$. Note that for an isotropic source $S=S(x)$, $\textup{Vol}(\mathbb{S}^{n-1})^{\frac1p}\|S\|_{C(X)} = \|S\|_{L^p(\mathbb{S}^{n-1},C(X))}$.

\bigskip

Let $v=v(x,\theta)$ be the solution to the following adjoint RTE with prescribed outgoing boundary condition $g$:
\begin{align}
-\theta\cdot\nabla v+\sigma v & - \int_{\mathbb{S}^{n-1}}k(x,\theta,\theta')v(x,\theta')\dif\theta' = 0 \label{eqn:aRTE} \\
v|_{\Gamma_+} & = g \label{eq:aRTEbc}.
\end{align}
Since we assume $\sigma$ and $k$ are known, we can solve this boundary value problem to find $v$ for any given $g$.

Next, we derive an internal functional of $u$ from the boundary measurement $\Lambda^\varepsilon_S$ in~\eqref{eq:Lambda}. To this end, we multiply \eqref{eqn:mRTE} by $v(x,\theta)$ to get
\[v\theta\cdot\nabla u_\varepsilon+\sigma_\varepsilon(x)u_\varepsilon v-\int_{\mathbb{S}^{n-1}}k_\varepsilon(x,\theta,\theta')u_\varepsilon(x,\theta')v(x,\theta)\dif\theta'  = S_\varepsilon(x)v(x,\theta)\]
and multiply \eqref{eqn:aRTE} by $u_\varepsilon(x,\theta)$ to get
\[-u_\varepsilon\theta\cdot\nabla v+\sigma u_\varepsilon v - \int_{\mathbb{S}^{n-1}}k(x,\theta,\theta')v(x,\theta')u_\varepsilon(x,\theta)\dif\theta' = 0.\]
Thanks to the condition \eqref{eq:isotropy}, the roles of $\theta$ and $\theta'$ can be interchanged in the integrals. Therefore subtracting these two equalities and then integrating over $X\times \mathbb{S}^{n-1}$ gives
\[\begin{aligned}\int_X\int_{\mathbb{S}^{n-1}}v\theta\cdot\nabla u_\varepsilon+u_\varepsilon\theta\cdot\nabla v\dif\theta\dif x=\int_X\int_{\mathbb{S}^{n-1}}\int_{\mathbb{S}^{n-1}}(k_\varepsilon-k)v(x,\theta)u_\varepsilon(x,\theta')\dif\theta\dif\theta'\dif x\\-\int_X\int_{\mathbb{S}^{n-1}}(\sigma_\varepsilon-\sigma)u_\varepsilon v\dif\theta\dif x+\int_X\int_{\mathbb{S}^{n-1}}S_\varepsilon v\dif\theta\dif x\end{aligned}\]
On the left-hand side, we apply the following integration-by-parts formula
\begin{equation} \label{eq:ibp}
\int_Xv\theta\cdot\nabla u_\varepsilon\dif x=-\int_Xu_\varepsilon\theta\cdot\nabla v\dif x+\int_{\partial X}u_\varepsilon vn\cdot\theta\dif x
\end{equation}
to obtain
\begin{equation}\label{eqn:11}
\begin{aligned}
\int_{\mathbb{S}^{n-1}}\int_{\partial X}u_\varepsilon vn\cdot\theta\dif x\dif\theta=&\int_X\int_{\mathbb{S}^{n-1}}\int_{\mathbb{S}^{n-1}}(k_\varepsilon-k)v(x,\theta)u_\varepsilon(x,\theta')\dif\theta\dif\theta'\dif x\\&+\int_X\int_{\mathbb{S}^{n-1}}vS_\varepsilon\dif\theta\dif x-\int_X\int_{\mathbb{S}^{n-1}}(\sigma_\varepsilon-\sigma)u_\varepsilon v\dif\theta\dif x
\end{aligned}
\end{equation}
When $\varepsilon=0$, that is, in the absence of acoustic modulation, Equation~\eqref{eqn:11} gives
\begin{equation}\label{eqn:12}
\int_{\mathbb{S}^{n-1}}\int_{\partial X}uvn\cdot\theta\dif x\dif\theta=\int_{\mathbb{S}^{n-1}}\int_{X}vS\dif x\dif\theta
\end{equation}
Subtract~\eqref{eqn:12} from~\eqref{eqn:11} to get
\begin{equation}
\begin{aligned}
\int_{\mathbb{S}^{n-1}}\int_{\partial X}(u_\varepsilon-u) vn\cdot\theta\dif x\dif\theta=&\int_X\int_{\mathbb{S}^{n-1}}\int_{\mathbb{S}^{n-1}}(k_\varepsilon-k)v(x,\theta)u_\varepsilon(x,\theta')\dif\theta\dif\theta'\dif x\\&+\int_X\int_{\mathbb{S}^{n-1}}v(S_\varepsilon-S)\dif\theta\dif x-\int_X\int_{\mathbb{S}^{n-1}}(\sigma_\varepsilon-\sigma)u_\varepsilon v\dif\theta\dif x
\end{aligned}
\end{equation}
To separate the $O(\varepsilon)$-term in $u_\varepsilon$, we write $u_\varepsilon = u_0 + \varepsilon \delta u$.
Substituting the expressions~\eqref{eq:msigma} \eqref{eq:mk} \eqref{eq:mS} and comparing the $O(\varepsilon)$-terms yield
\begin{equation}
\begin{aligned}
\int_{\mathbb{S}^{n-1}}\int_{\partial X}& \delta u vn\cdot\theta\dif x\dif\theta\\=&-\int_X\int_{\mathbb{S}^{n-1}}\cos(q\cdot x+\varphi)\sigma u v\dif\theta\dif x+\int_X\int_{\mathbb{S}^{n-1}}\cos(q\cdot x+\varphi)vS\dif\theta\dif x\\&+\int_X\int_{\mathbb{S}^{n-1}}\int_{\mathbb{S}^{n-1}}\cos(q\cdot x+\varphi)k(x,\theta,\theta')v(x,\theta)u(x,\theta')\dif\theta\dif\theta'\dif x+O(\varepsilon)
\end{aligned}
\end{equation}
Since $\delta u|_{\Gamma_-} = \frac{1}{\varepsilon} (u_\varepsilon - u_0)|_{\Gamma_-}=0$ and $\delta u|_{\Gamma_+} = \frac{1}{\varepsilon} (u_\varepsilon - u_0)|_{\Gamma_+} = \frac{1}{\varepsilon} (\Lambda^\varepsilon_S(q,\varphi) - \Lambda^0_S(q,\varphi))$, 
the left-hand side is known from the measurement for any $q$ and $\varphi$. Varying $q$ and $\varphi$, we obtain from the right-hand side the Fourier transform of the quantity $H_v$ defined by
\begin{equation}\label{eqn:17}
\begin{aligned}
H_v(x)\coloneqq&-\int_{\mathbb{S}^{n-1}}\sigma uv\dif\theta+\int_{\mathbb{S}^{n-1}}vS\dif\theta\\&+\int_{\mathbb{S}^{n-1}}\int_{\mathbb{S}^{n-1}}k(x,\theta,\theta')v(x,\theta)u(x,\theta')\dif\theta'\dif\theta
\end{aligned}
\end{equation}
Substitute the RTE \eqref{eqn:RTE} to get
\begin{equation} \label{eq:intfunc}
\begin{aligned}
H_v(x)=&-\int_{\mathbb{S}^{n-1}}\sigma uv\dif\theta+\int_{\mathbb{S}^{n-1}}vS\dif\theta+\int_{\mathbb{S}^{n-1}}v(x,\theta)(\theta\cdot\nabla u+\sigma u-S)\dif\theta\\
= & \int_{\mathbb{S}^{n-1}}v(x,\theta)\theta\cdot\nabla u(x,\theta)\dif\theta \\
= & \int_{\mathbb{S}^{n-1}}v(x,\theta) 
[ \mathcal{A} u(x,\theta) + S(x)]
\dif\theta \\
\end{aligned}
\end{equation}
where the operator $\mathcal{A}$ is defined as 
\begin{equation} \label{eq:opA}
\mathcal{A} u(x,\theta)\coloneqq-\sigma(x)u(x,\theta)+\int_{\mathbb{S}^{n-1}}k(x,\theta,\theta')u(x,\theta')\dif\theta'.
\end{equation}
We therefore have extracted the internal functional $H_v$ from the measurement $\Lambda^\varepsilon_S$.

\section{Inversion Theory and Formulae}

In this section, we assume knowledge of the quantity $H_v$ and derive two algorithms to reconstruct the isotropic source $S$. The first is based on computation of a Neumann series, and the second amounts to solving a Fredholm equation. Our starting point is the following relation, see~\eqref{eq:intfunc}.
$$
H_{v}(x)=\int_{\mathbb{S}^{n-1}} \mathcal{A}u(x,\theta)v(x,\theta)\dif\theta + S(x)\int_{\mathbb{S}^{n-1}}v(x,\theta)\dif\theta.
$$
We need the following simple fact -- see the appendix for a short proof, or ~\cite{DautrayLions, DavisonSykes} for similar results. 

\begin{lemma}
For any uniformly positive function $f_0 \in C(\Gamma_+)$, there exists a continuous adjoint RTE solution $v_0$ to~\eqref{eqn:aRTE} with boundary condition $v_0|_{\Gamma_+} = f_0$, and a constant $c>0$ such that $v_0(x,\theta) \geq c > 0$ and $\int_{\mathbb{S}^{n-1}}v_0(x,\theta)\dif\theta\geq c>0$ for any $(x,\theta)\in X\times \mathbb{S}^{n-1}$.
\end{lemma}

Note that in particular we can choose $f_0 \equiv 1$ on the boundary, in which case the internal functional $H_{v_0}$ corresponds to the measurements obtained from the integral
\[
\int_{\Gamma_+} (u - u_{\varepsilon}) \, dS;
\]
in other words it can be obtained from measurements of the angular average of $u$ and $u_{\varepsilon}$ on the boundary. 

Let $v_0$ be an adjoint RTE solution as in the above lemma. Dividing the internal functional $H_{v_0}$ by $\int_{\mathbb{S}^{n-1}}v_0(x,\theta)\dif\theta$, we obtain
\begin{equation} \label{eq:F}
\frac{H_{v_0}(x)}{\int_{\mathbb{S}^{n-1}}v_0(x,\theta)\dif\theta} \coloneqq  S(x) + \frac{\int_{\mathbb{S}^{n-1}}\mathcal{A}u(x,\theta)v_0(x,\theta)\dif\theta}{\int_{\mathbb{S}^{n-1}}v_0(x,\theta)\dif\theta}.
\end{equation}
We will regard the second term on the right-hand side as a linear operator of $S$.

\subsection{Neumann Series Inversion} \label{sec:Neumann}
We derive a Neumann series inversion formula based on~\eqref{eq:F}. To this end, let us introduce three linear operators. The first operator is
\begin{equation}
\mathcal{S}: C(X)\rightarrow L^p(\mathbb{S}^{n-1}, C(X)), \quad\quad\quad S \mapsto u
\end{equation}
where $u$ is the solution to the boundary value problem~\eqref{eqn:RTE} \eqref{eq:bc}. Here $\mathcal{S}$ is just the source-to-solution operator. It is bounded under the assumptions (A1)(A2), and by Proposition~\ref{thm:WP},
\begin{equation}\label{eq:Sbound}
\|\mathcal{S}\|_{C(X)\rightarrow L^p(\mathbb{S}^{n-1}, C(X))}\leq
\begin{cases}
\frac{\textup{Vol}(\mathbb{S}^{n-1})^{\frac 1 p}}{\alpha} & \left(\inf\limits_{x\in \overline{X}}\sigma\right)-\rho\geq\alpha\\
\frac{\textup{diam}(X)\textup{Vol}(\mathbb{S}^{n-1})^{\frac 1 p}}{1-\textup{diam}(X)\rho} &\textup{diam}(X)\rho<1
\end{cases}
\end{equation}

The second operator is
\begin{equation}
\mathcal{K}_{v_0}: L^p(\mathbb{S}^{n-1}, C(X)) \rightarrow C(X), \quad\quad\quad u(x,\theta) \mapsto \int_{\mathbb{S}^{n-1}}\mathcal{A}u(x,\theta)v_0(x,\theta)\dif\theta
\end{equation}
where the operator $\mathcal{A}$ is introduced in~\eqref{eq:opA}. Based on the estimate


\begin{equation}
\begin{aligned}
& \|\mathcal{K}_{v_0}u\|_{C(X)} = 
\left\|\int_{\mathbb{S}^{n-1}}\mathcal{A}u(x,\theta)v_0(x,\theta)\dif\theta\right\|_{C(X)} \vspace{1ex}\\
\leq&\left\|\int_{\mathbb{S}^{n-1}} (\sigma uv_0)(x,\theta) \dif\theta\right\|_{C(X)}+\left\|\int_{\mathbb{S}^{n-1}}\int_{\mathbb{S}^{n-1}}k(x,\theta,\theta')u(x,\theta')v_0(x,\theta)\dif\theta'\dif\theta\right\|_{C(X)} \vspace{1ex} \\
\leq&\|v_0\|_{C(X)}\left(\left\|\int_{\mathbb{S}^{n-1}}(\sigma u)(x,\theta)\dif\theta\right\|_{C(X)}+\left\|\int_{\mathbb{S}^{n-1}}\int_{\mathbb{S}^{n-1}}k(x,\theta,\theta')u(x,\theta')\dif\theta'\dif\theta\right\|_{C(X)}\right) \vspace{1ex} \\
\leq&\|v_0\|_{C(X)}\left(\|\sigma\|_{C(X)}\left\|\int_{\mathbb{S}^{n-1}}u(x,\theta)\dif\theta\right\|_{C(X)}+\rho\left\|\int_{\mathbb{S}^{n-1}}u(x,\theta')\dif\theta'\right\|_{C(X)}\right) \vspace{1ex} \\
=&\|v_0\|_{C(X)}(\|\sigma\|_{C(X)}+\rho)\left\|\int_{\mathbb{S}^{n-1}}u(x,\theta)\dif\theta\right\|_{C(X)} \vspace{1ex} \\
\leq&\|v_0\|_{C(X)}(\|\sigma\|_{C(X)}+\rho)\textup{Vol}(\mathbb{S}^{n-1})^{1-\frac1p}\left\|u\right\|_{L^p(\mathbb{S}^{n-1},C(X))},
\end{aligned}
\end{equation}
where the last line follows from H\"{o}lder's inequality.
we see that $\mathcal{K}_{v_0}$ is a bounded operator and 
\begin{equation} \label{eq:Kbound}
\|\mathcal{K}_{v_0}\|_{L^p(\mathbb{S}^{n-1}, C(X)) \rightarrow C(X)} \leq \|v_0\|_{C(X)}(\|\sigma\|_{C(X)}+\rho)\textup{Vol}(\mathbb{S}^{n-1})^{1-\frac1p}
\end{equation}

The third operator is the multiplication operator
\begin{equation}
\mathcal{M}_{v_0}: C(X) \rightarrow C(X), \quad\quad\quad f(x) \mapsto \frac{1}{\int_{\mathbb{S}^{n-1}} v_0(x,\theta)\dif\theta} f(x).
\end{equation}
It is bounded since $v_0$ is chosen in such a way that $\int_{\mathbb{S}^{n-1}} v_0(x,\theta)\dif\theta$ is bounded away from zero. We have
\begin{equation} \label{eq:Mbound}
\|\mathcal{M}_{v_0}\|_{C(X) \rightarrow C(X)} \leq 
\frac{1}{ \inf_{x\in\overline{X}} \left( \int_{\mathbb{S}^{n-1}} v_0(x,\theta)\dif\theta \right)}.
\end{equation}

Using these operators, the equation~\eqref{eq:F} can be written as
$$
\mathcal{M}_{v_0}[H_{v_0}] = (Id + \mathcal{M}_{v_0}\circ\mathcal{K}_{v_0}\circ\mathcal{S})[S].
$$
where $Id$ is the identity operator. Here the left-hand side is the known from the internal functional and the choice of $v_0$. It remains to invert the operator $Id + \mathcal{M}_{v_0}\circ\mathcal{K}_{v_0}\circ\mathcal{S}$ to find the source $S$. This leads naturally to a Neumann series reconstruction if the operator $\mathcal{M}_{v_0}\circ\mathcal{K}_{v_0}\circ\mathcal{S}$ is a contraction. Note from~\eqref{eq:Sbound} \eqref{eq:Kbound} \eqref{eq:Mbound} that
\begin{equation} \label{eq:twobounds}
\|\mathcal{M}_{v_0}\circ\mathcal{K}_{v_0}\circ\mathcal{S}\|_{C(X)\rightarrow C(X)} \leq 
\begin{cases}
\frac{\|v_0\|_{C(X)}(\|\sigma\|_{C(X)}+\rho)\textup{Vol}(\mathbb{S}^{n-1})}{\alpha \inf_{x\in\overline{X}} \left( \int_{\mathbb{S}^{n-1}} v_0(x,\theta)\dif\theta \right)} & \left(\inf\limits_{x\in \overline{X}}\sigma\right)-\rho\geq\alpha\\
\frac{\|v_0\|_{C(X)}(\|\sigma\|_{C(X)}+\rho)\textup{diam}(X)\textup{Vol}(\mathbb{S}^{n-1})}{(1-\textup{diam}(X)\rho) \inf_{x\in\overline{X}} \left( \int_{\mathbb{S}^{n-1}} v_0(x,\theta)\dif\theta \right)} &\textup{diam}(X)\rho<1
\end{cases}
\end{equation}

If either bound on the right-hand side is strictly less than $1$, then the operator $\mathcal{M}_{v_0}\circ\mathcal{K}_{v_0}\circ\mathcal{S}$ is a contraction. 
Here the first bound in ~\eqref{eq:twobounds} is not helpful, since 
\[\frac{1}{\textup{Vol}(\mathbb{S}^{n-1})}\int_{\mathbb{S}^{n-1}} v_0(x,\theta)\dif\theta\leq\|v_0\|_{C(X)},
\quad\quad \text{ and }
\quad\alpha\leq\left(\inf\limits_{x\in \overline{X}}\sigma\right)-\rho\leq\|\sigma\|_{C(X)}+\rho,\]
which imply
\[\frac{\|v_0\|_{C(X)}(\|\sigma\|_{C(X)}+\rho)\textup{Vol}(\mathbb{S}^{n-1})}{\alpha \inf_{x\in\overline{X}} \left( \int_{\mathbb{S}^{n-1}} v_0(x,\theta)\dif\theta \right)}
\geq 
\frac{\int_{\mathbb{S}^{n-1}} v_0(x,\theta)\dif\theta}{ \inf_{x\in\overline{X}} \left( \int_{\mathbb{S}^{n-1}} v_0(x,\theta)\dif\theta \right)}
\geq1.\]
On the other hand, the second bound in ~\eqref{eq:twobounds} shows that $\mathcal{M}_{v_0}\circ\mathcal{K}_{v_0}\circ\mathcal{S}$ is a contraction if the domain $X$ is small enough, meaning that \eqref{eq:F} can be inverted through a Neumann series.  This is numerically demonstrated in Section~\ref{sec:numerics}.

It is also not necessarily clear that the first bound is always sharp (see Experiment 2 in Section 4.2).

Summarizing the discussion above, we have
\begin{thm} \label{thm:Neumann}
Suppose the assumptions (A1)(A2) hold. If the following inequality holds
\begin{equation} \label{eq:combound}
\begin{aligned}
\frac{\|v_0\|_{C(X)}(\|\sigma\|_{C(X)}+\rho)\textup{diam}(X)\textup{Vol}(\mathbb{S}^{n-1})}{(1-\textup{diam}(X)\rho) \inf_{x\in\overline{X}} \left( \int_{\mathbb{S}^{n-1}} v_0(x,\theta)\dif\theta \right)}<1 &\text{ when }\textup{diam}(X)\rho<1,
\end{aligned}
\end{equation}
then the operator $\mathcal{M}_{v_0}\circ\mathcal{K}_{v_0}\circ\mathcal{S}$ is a contraction, and the source $S$ can be computed from the following Neumann series:
$$
S = \sum^\infty_{j=0} (-\mathcal{M}_{v_0}\circ\mathcal{K}_{v_0}\circ\mathcal{S})^j (\mathcal{M}_{v_0}[H_{v_0}]).
$$
\end{thm}

\subsection{Fredholm Inversion} \label{sec:Fredholm}
The assumption~\eqref{eq:combound} is a bit too strong and may be invalid in certain circumstances. In this section, we derive another inversion formula which removes such restriction. Let us begin by introducing some function spaces. For $1\leq p \leq \infty$, define
$$
\mathcal{H}^1_p\coloneqq\left\{u\in L^p(X\times\mathbb{S}^{n-1}) \mid \theta\cdot\nabla u\in L^p(X\times\mathbb{S}^{n-1})\right\}.
$$

For $1\leq p<\infty,\theta\in(0,1),f\in L^p(X),$ the Slobodeckij seminorm is defined by
\[[f]_{\theta,p,X}\coloneqq\left(\int_X\int_X\frac{|f(x)-f(y)|^p}{|x-y|^{\theta p+n}}\dif x\dif y\right)^{\frac1p}\]
Let $s>0$ be a non-integer and set $\theta=s-[s]$, the Sobolev space $W^{s,p}$ is defined as
\[W^{s,p}(X)\coloneqq\left\{u\in W^{[s],p}(X)\middle|\sup_{|\alpha|=[s]}[D^\alpha u]_{\theta,p,X}<\infty\right\}\] 
with norm $\|f\|_{W^{s,p}(X)}\coloneqq\|f\|_{W^{[s],p}(X)}+\sup\limits_{|\alpha|=[s]}[D^\alpha u]_{\theta,p,X}.$

Henceforth, we restrict to the case $p=2$ and work on $L^2$-based spaces.
We make the following further assumptions on the optical coefficients. 

\smallskip
\begin{enumerate}
\item[(A3):] $\sigma(x) \geq \sigma_0>0$ everywhere in $X$ for some constant $\sigma_0$.
\item[(A4):] $\|\frac1{\sigma(x)}\int_{\mathbb{S}^{n-1}}k(x,\theta,\theta')\dif\theta'\|_{L^\infty(X \times S^{n-1})}\leq k_0<1$ for some constant $k_0$.
\item[(A5):] $\sigma(x)\in W^{1,2}(X),$ $k(x,\theta,\theta')\in W^{1,2}(X)$ for any $\theta,\theta'\in\mathbb{S}^{n-1}.$
\end{enumerate}
Here (A3) and (A4) are imposed to ensure solvability of the forward boundary value problem~\eqref{eqn:RTE}~\eqref{eq:bc}
in the space $\mathcal{H}_2^1$, see Proposition~\ref{thm:4} below. (A5) is needed when applying the averaging lemma.

\begin{prop}[{\cite[Theorem 3.2]{10.5555/286912}}]\label{thm:4}
For any $S(x)\in L^2(X)$, the boundary value problem~\eqref{eqn:RTE}~\eqref{eq:bc} admits a unique solution $u\in\mathcal{H}_2^1$. Moreover, 
the following estimate holds for some constants $C,\;\tilde C>0$ independent of $S$ and $u:$
\[C\|S\|_{L^2(X)}\leq\|u\|_{\mathcal{H}_2^1}\leq\tilde C\|S\|_{L^2(X)}.\]
\end{prop}

Since $X$ is bounded and $S(x)\in C(X)$, we have $S(x)\in L^2(X)$, hence $u\in\mathcal{H}_2^1$ by Proposition~\ref{thm:4}. Similarly, we have $v_0,\sigma v_0\in L^2(X\times\mathbb{S}^{n-1}).$ Moreover,
\begin{align*}
&\left(\int_X\int_{\mathbb{S}^{n-1}}\left|\int_{\mathbb{S}^{n-1}}k(x,\theta,\theta')v_0(x,\theta')\dif\theta'\right|^2\dif\theta\dif x\right)^{\frac12}\\
\leq&\left(\int_X\int_{\mathbb{S}^{n-1}}\int_{\mathbb{S}^{n-1}}(\sup|k|)^2\left|v_0(x,\theta')\right|^2\dif\theta'\dif\theta\dif x\right)^{\frac12}\\
=&\sup|k|\textup{Vol}(\mathbb{S}^{n-1})^{\frac12}\|v_0\|_{L^2(X\times\mathbb{S}^{n-1})}<\infty,
\end{align*}
then from \eqref{eqn:aRTE}, we have $\theta\cdot\nabla v_0(x,\theta)\in L^2(X\times\mathbb{S}^{n-1})$. Thus $v_0\in\mathcal{H}_2^1.$
%
By the assumption (A5), we conclude $\sigma uv_0\in\mathcal{H}_2^1$ and $\int_{\mathbb{S}^{n-1}}k(x,\theta,\theta')u(x,\theta)v_0(x,\theta)\dif\theta'\in\mathcal{H}_2^1.$ By the Averaging Lemma (see~\cite[Theorem 1.1]{devore2001averaging}), $\mathcal{K}_{v_0}\circ\mathcal{S}[S]\in W^{\frac12,2}(X).$
As the embedding $W^{\frac12,2}(X) \xhookrightarrow{} L^2(X)$ is compact, the operator $\mathcal{K}_{v_0}\circ\mathcal{S}$ is a compact operator from $(C(X),\|\cdot\|_2)$ to $L^2(X)$, which can be extend to be a compact operator defined on the entire space $L^2(X)$. We slightly abuse the notation and denote such extension again by $\mathcal{K}_{v_0}\circ\mathcal{S}$. On the other hand, the multiplication operator $\mathcal{M}_{v_0}$ can be extended to be a bounded operator on $L^2(X)$. Thus, the operator $\mathcal{M}_{v_0}\circ\mathcal{K}_{v_0}\circ\mathcal{S}: L^2(X)\rightarrow L^2(X)$, as the composition of a bounded operator with a compact operator, is compact as well.
We therefore have the following result due to the Fredholm alternative.

%
%

\begin{thm} \label{thm:Fredholm}
Suppose the assumptions (A1)\textasciitilde(A5) hold. If $0$ is not an eigenvalue of the Fredholm operator $Id + \mathcal{M}_{v_0}\circ\mathcal{K}_{v_0}\circ\mathcal{S}$, then $(Id + \mathcal{M}_{v_0}\circ\mathcal{K}_{v_0}\circ\mathcal{S})^{-1}$ is a bounded linear operator on $L^2(X)$, and the source $S$ can be computed as
$$
S = (Id + \mathcal{M}_{v_0}\circ\mathcal{K}_{v_0}\circ\mathcal{S})^{-1} (\mathcal{M}_{v_0}[H_{v_0}]).
$$
\end{thm}

The following stability estimate is an immediate consequence of this inversion formula.
\begin{coro}
Suppose the assumptions (A1)\textasciitilde(A5) hold. Let $S$ and $\tilde{S}$ be two different sources with corresponding internal functional $H_{v_0}$ and $\tilde{H}_{v_0}$, respectively. If $0$ is not an eigenvalue of the operator $Id + \mathcal{M}_{v_0}\circ\mathcal{K}_{v_0}\circ\mathcal{S}$, then the following stability estimate holds
$$
\|S-\tilde{S}\|_{L^2(X)} \leq C \|H_{v_0}-\tilde{H}_{v_0}\|_{L^2(X)}
$$
for some constant $C>0$ depending on $\sigma$, $k$, $v_0$, $X$ yet independent of $S$ and $\tilde{S}$.
\end{coro}

\section{Algorithms and Numerical Experiments}. \label{sec:numerics}

In this section, we implement the proposed source reconstruction procedures in 2D for Theorem~\ref{thm:Neumann} and Theorem~\ref{thm:Fredholm} . We write $(x_1,x_2)$ for the coordinates of a point. The computational domain $X$ is a square whose size will be individually specified in each experiment. The scattering kernel is chosen as the Henyey-Greenstein function
\begin{equation} \label{eq:HGfunc}
k(x,\theta,\theta') = \frac{1}{2\pi} \frac{1-g^2}{1+g^2-2g\cos\phi},
\end{equation}
where $\phi$ is the angle between $\theta$ and $\theta'$, and $-1\leq g \leq 1$ is the anisotropy parameter of the medium.

\subsection{Description of the Algorithms.}

We briefly explain the forward and inverse solvers involved in the numerical experiments below. The forward solver is used to solve the RTE and adjoint RTE, while the inverse solvers implement the Neumann series reconstruction in  Theorem~\ref{thm:Neumann} and the Fredholm inversion in Theorem~\ref{thm:Fredholm}.

\subsubsection{Radiative Transfer Equation.}.
The RTE~\eqref{eqn:RTE} with the zero boundary condition~\eqref{eq:bc} is solved using the discrete ordinate method~\cite{wang_sheng_han_2016}.
Firstly, we uniformly discretize the angular space $[0,2\pi)$ into $M$ angles. To this end, set $\Delta\omega=\frac{2\pi}{M}$ and choose the discrete angles
$\omega_i = (i-1)\Delta\omega$, $i=1,2,\dots, M$ and denote $\theta_i=(\cos\omega_i,\sin\omega_i)$.
Using the trapezoidal rule, we have the approximation
\[\int_{\mathbb{S}^1}k(x,\theta',\theta)u(x,\theta)\dif\theta\approx\sum_{i=1}^Mk(x,\theta',\theta_i)u(x,\theta_i)\Delta\omega,\]
After the angular discretization, the resulting equations form a hyperbolic system:
\[\begin{aligned}\theta_i\cdot\nabla u(x,\theta_i)+\sigma(x) u(x,\theta_i)-\sum_{j=1}^Mk(x,\theta_i,\theta_j)u(x,\theta_j)\Delta\omega&=S(x)&&1\leq i\leq M\\u(x,\theta_i)&=0&&(x,\theta_i)\in\Gamma_-\end{aligned}\].

Secondly, we use the upwind scheme for spatial discretization, that is, 
\begin{align}
\frac{\partial u}{\partial x_1}(x_1,x_2,\theta_i)  &\approx \sgn(\cos\omega_i)\frac{u(x_1+\sgn(\cos\omega_i)\Delta x_1,x_2,\theta_i)-u(x_1,x_2,\theta_i)}{\Delta x_1}, \label{eq:partialx1}\\
\frac{\partial u}{\partial x_2}(x_1,x_2,\theta_i) &\approx \sgn(\sin\omega_i)\frac{u(x_1,x_2+\sgn(\sin\omega_i)\Delta x_2,\theta_i)-u(x_1,x_2,\theta_i)}{\Delta x_2}. \label{eq:partialx2}
\end{align}
where $\Delta x_1$ and $\Delta x_2$ are the spacings along the $x_1$-direction and $x_2$-direction, respectively.
We remark that for an angle $\omega_i$ such that $\cos\omega_i\neq 0$, the right-hand side of~\eqref{eq:partialx1} is a valid approximation of the derivative $\frac{\partial u}{\partial x_1}(x_1,x_2,\theta_i)$; for an angle $\omega_i$ such that $\cos\omega_i = 0$, the right-hand side of~\eqref{eq:partialx1} becomes zero and the approximation fails. However, this does not affect numerical calculation of the directional derivative $\theta_i\cdot\nabla u(x,\theta_i)$ since $\frac{\partial u}{\partial x_1}(x_1,x_2,\theta_i)$ is multiplied by $\cos\omega_i = 0$ there. Similar remark applies to~\eqref{eq:partialx2} for $\omega_i$ such that $\sin\omega_i=0$.

The spatial discretization ends up with a linear system with prescribed zero boundary values on $\Gamma_-$, which is then solved using the Jacobi iteration. The adjoint boundary value problem~\eqref{eqn:aRTE} \eqref{eq:aRTEbc} is solved in a similar manner, yet with upwind directions specified by $-\theta$ and boundary values specified by the function $g$.

Given a known source $S$, we generate the measurement $H_v(x)$ in the following steps. First, we solver the forward problem~\eqref{eqn:RTE} \eqref{eq:bc} using the RTE solver to find the solution $u(x,\theta)$. This, together with the known attenuation coefficient and scattering kernel, is employed to compute $\mathcal{A}u(x,\theta)$ in~\eqref{eq:opA}. Finally, we solve the adjoint RTE~\eqref{eqn:aRTE} \eqref{eq:aRTEbc} to get $v$, and compute $H_v(x)$ in~\eqref{eq:intfunc} with the trapezoidal rule.

%
%


\subsubsection{Neumann Series Inversion.}

In order to implement the Neumann series inversion in Theorem~\ref{thm:Neumann}, we discretize the operator $\mathcal{S}$ by solving the forward RTE, and the operators 
$\mathcal{K}_{v_0}$ and $\mathcal{M}_{v_0}$ using the trapezoidal rule.

The algorithm for Theorem~\ref{thm:Neumann} is simple. The operator $\mathcal{S}$ can be implemented using the forward RTE solver, the operator $\mathcal{K}_{v_0}$ and $\mathcal{M}_{v_0}$ can be discretized using the trapezoidal rule, then the reconstruction can be done by an iteration.
\begin{algorithm}[H]
\KwData{adjoint RTE solution $v_0$, measurement $H_{v_0}$, scattering kernel $k(x,\theta,\theta')$, attenuation coefficient $\sigma(x)$, domain $X$.}
$S\leftarrow0$;\\
$\Delta S\leftarrow\mathcal{M}_{v_0}[H_{v_0}]$;\\
$\varepsilon\leftarrow10^{-6}$;

 \While{$\|\Delta S\|_{L^2} > \varepsilon$}{
  $S\leftarrow S+\Delta S$;\\
  $\Delta S\leftarrow\mathcal{M}_{v_0}\circ\mathcal{K}_{v_0}\circ\mathcal{S}[\Delta S]$;
 }
 \Return $S$;
 \caption{Neumann Series Reconstruction}
\end{algorithm}

\subsubsection{Fredholm Inversion.}

The Fredholm inversion in Theorem~\ref{thm:Fredholm} boils down to solving the linear system \eqref{eqn:Fredholm}. For this purpose, we descretize the source $S$ with respect to some basis functions. Two types of basis functions are used, one is polynomial functions of the form $\{x_1^i x_2^j\}_{i,j\geq 0,\; i+j\leq 10}$; the other is the pyramid-shaped functions
\[f_{ij}=\max\left\{1-\max\left\{20\left|x_1-\frac{i}{20}\right|,20\left|x_2-\frac{j}{20}\right|\right\},0\right\},\quad i,j\in\{0,1,\dots,20\}.\]
Polynomials capture the smooth feature of the source, while the pyramid-shaped functions capture some information of singularities. We write the expansion of a source $S$ with respect to these basis functions as
\begin{equation}
\label{eqn:approx}S(x_1,x_2)\approx\sum\limits_{i,j\geq0,i+j\leq10}c_{ij}x_1^ix_2^j+\sum_{0\leq i,j\leq20}c_{ij}'f_{ij}\eqqcolon\sum_i\tilde{c}_ib_i,
\end{equation}
where $c_{ij}$, $c'_{ij}$ are the coefficients of the expansion. We use $\{b_i(x_1,x_2)\}$ to denote these basis functions and $\{\tilde{c}_i\}$ the correponding coefficients.


Denote $\mathcal{T} := Id + \mathcal{M}_{v_0}\circ\mathcal{K}_{v_0}\circ\mathcal{S},$ then the internal measurement can be represented as \[\mathcal{M}_{v_0}[H_{v_0}]=\mathcal{T}[S]\approx\sum_i\tilde{c}_i\mathcal{T}[b_i].\]
We can compute the inner product with $\mathcal{T}[b_j]$ as follows:
\begin{equation}
\label{eqn:Fredholm}\langle\mathcal{M}_{v_0}[H_{v_0}],\mathcal{T}[b_j]\rangle\approx\sum_i\tilde{c}_i\langle\mathcal{T}[b_i],\mathcal{T}[b_j]\rangle.
\end{equation}
Solving the linear equation \eqref{eqn:Fredholm} gives the coefficient $\tilde{c}_i$, and then we can  numerically reconstruct the source $S$.

\subsection{Numerical Experiments.}

We demonstrate several numerical experiments in this section. For the forward problem, we discretize the angular space into $M=8$ directions, and the spatial domain into a $121\times121$ uniform grid. For the reconstruction, we interpolate the measurement with a spatial $61\times61$ uniform grid to avoid the inverse crime.

\textbf{Experiment 1: Inversion within the Assumption of Theorem~\ref{thm:Neumann}.}
In this experiment, we choose the quantities to satisfy the assumption~\eqref{eq:combound} in Theorem~\ref{thm:Neumann}.
The computational domain is $X=[0,0.2]\times [0,0.2]$; the attenuation coefficient is $\sigma_1(x_1,x_2) = 0.1+0.1x_1$; the anisotropy parameter is $g=0.5$ in the scattering kernel~\eqref{eq:HGfunc}; the function $v_0$ is the solution of \eqref{eqn:aRTE} with the boundary condition $v_0|_{\Gamma_+}=1$.
Such choice gives the following numerical values:
$$
\|v_0\|_{C(X)}\approx1.2603, \quad\quad\quad 
\inf_{x\in\overline{X}} \left( \int_{\mathbb{S}^{1}} v_0(x,\theta)\dif\theta \right)\approx6.4870.
$$
On the other hand, we have $\rho=1$ for any anisotropy parameter between $-1$ and $1$, thus
\[\frac{\|v_0\|_{C(X)}(\|\sigma_1\|_{C(X)}+\rho)\textup{diam}(X)\textup{Vol}(\mathbb{S}^{1})}{(1-\textup{diam}(X)\rho) \inf_{x\in\overline{X}} \left( \int_{\mathbb{S}^{1}} v_0(x,\theta)\dif\theta \right)}\approx0.5392<1,\]
so the assumption~\eqref{eq:combound} in Theorem~\ref{thm:Neumann} holds.

\begin{figure}[p]
\centering
\includegraphics[width=0.3\textwidth]{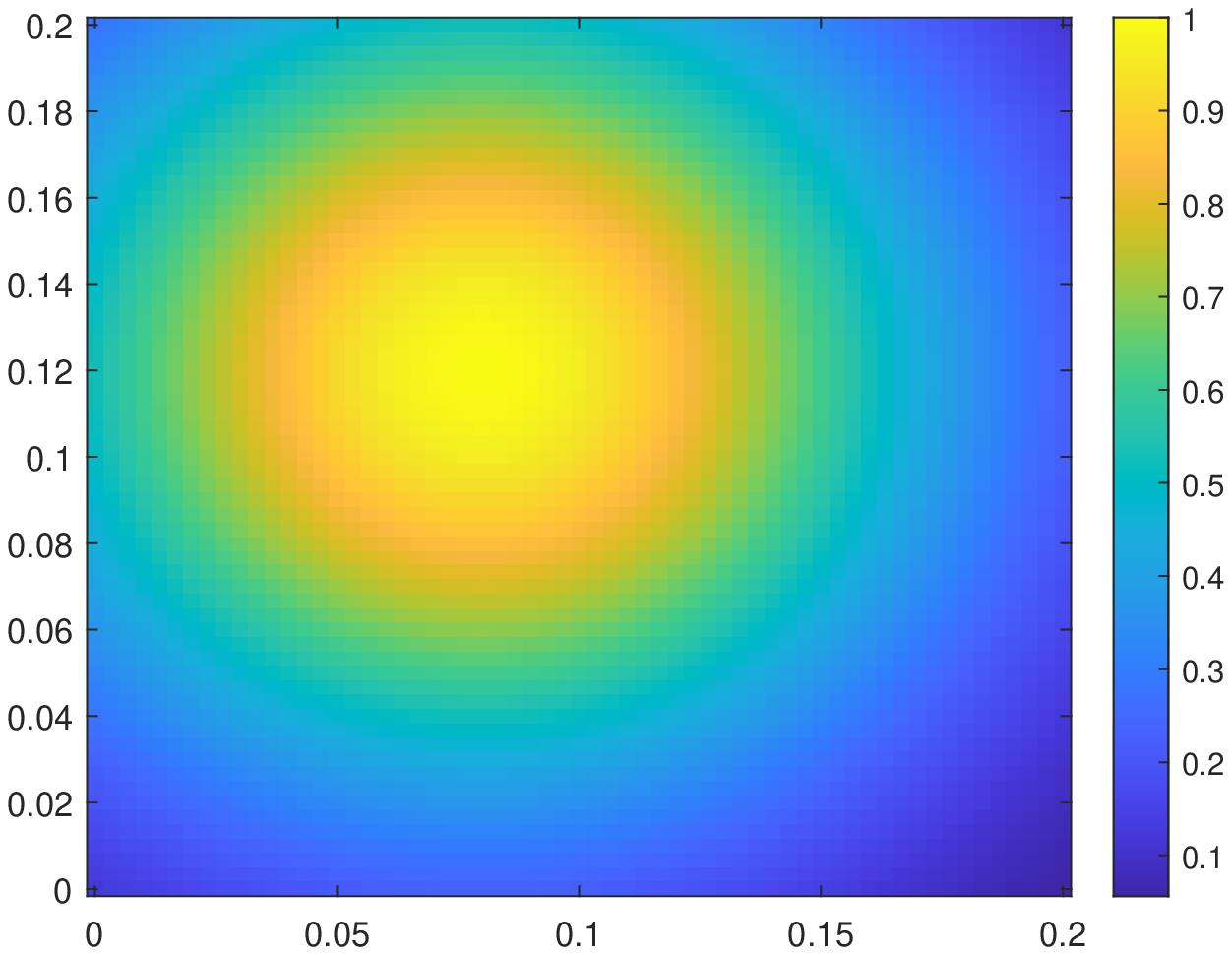}
\includegraphics[width=0.3\textwidth]{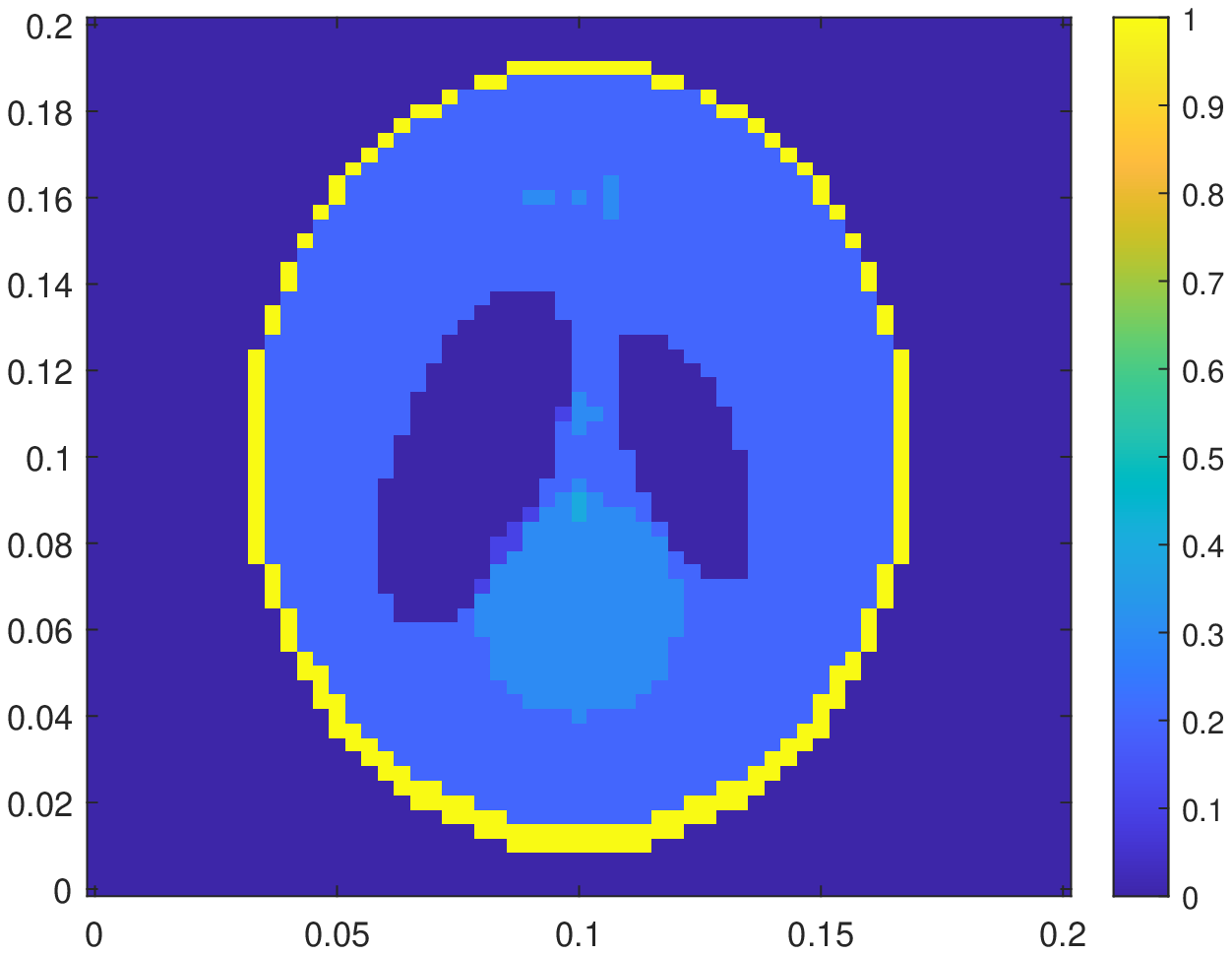}
\includegraphics[width=0.3\textwidth]{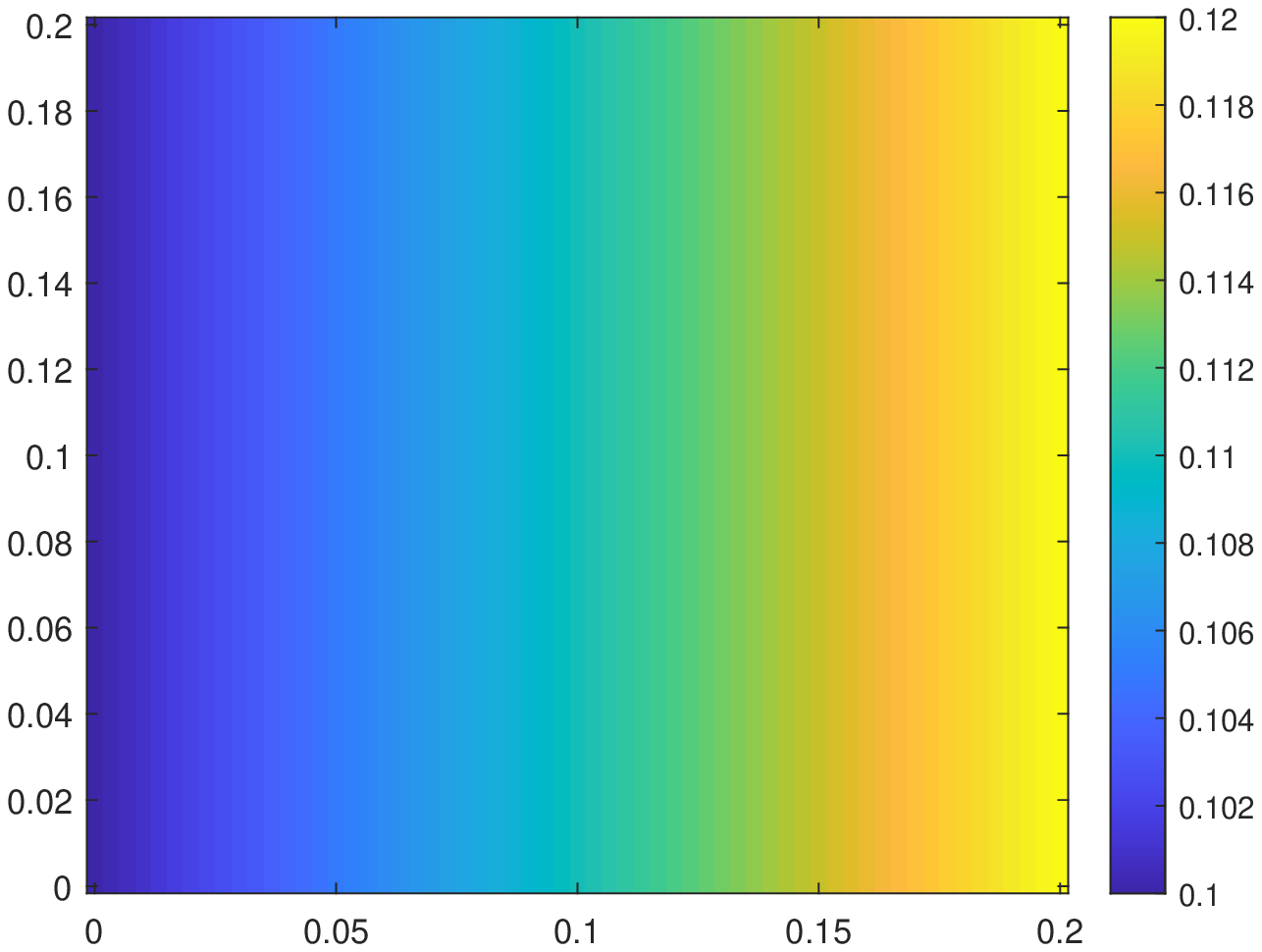}
\caption{Left: source $S_1$. Center: source $S_2$. Right: attenuation coefficient $\sigma_1$.}
\label{fig:Ncoef1}
\vspace*{\floatsep}
\includegraphics[width=0.23\textwidth]{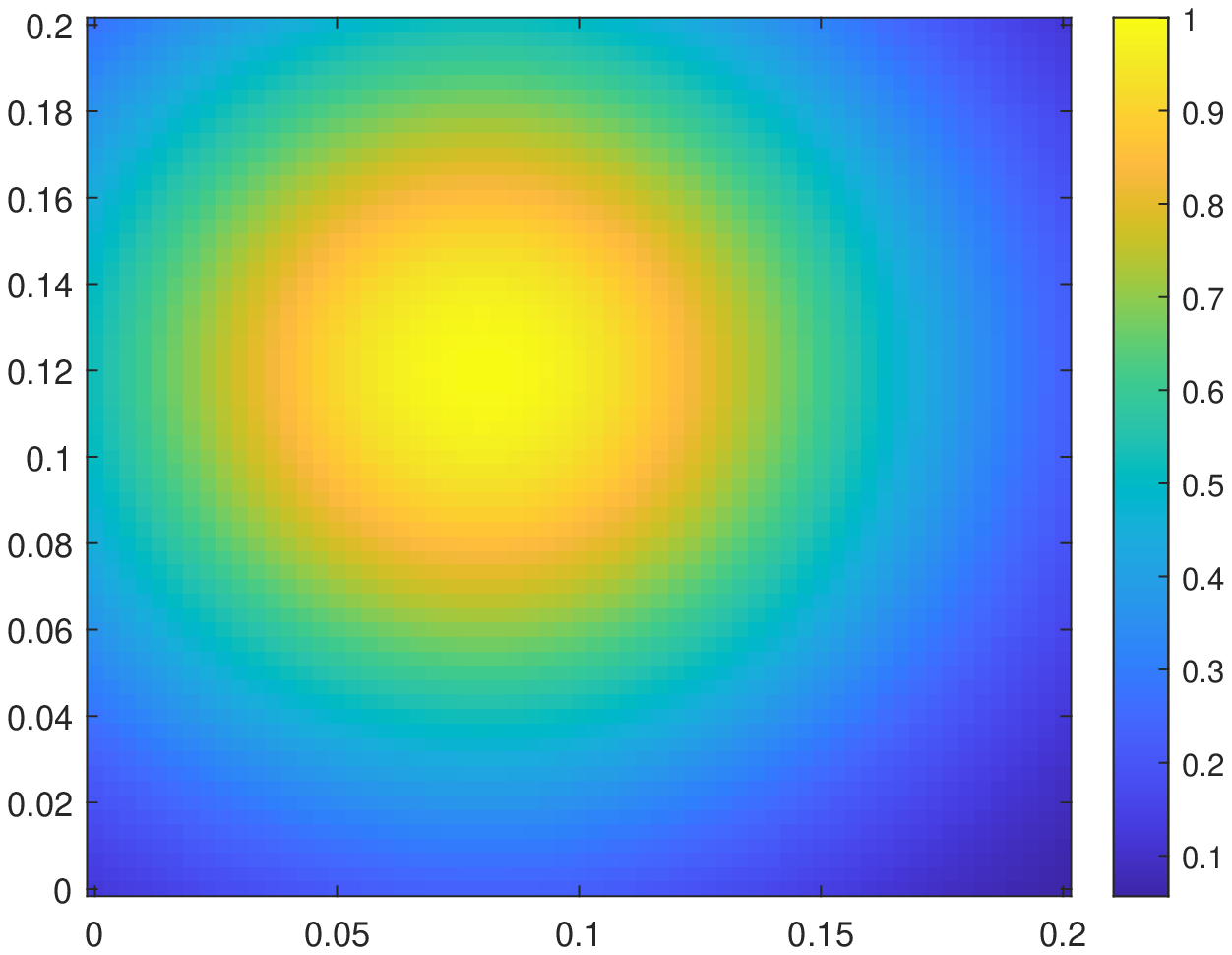}
\includegraphics[width=0.23\textwidth]{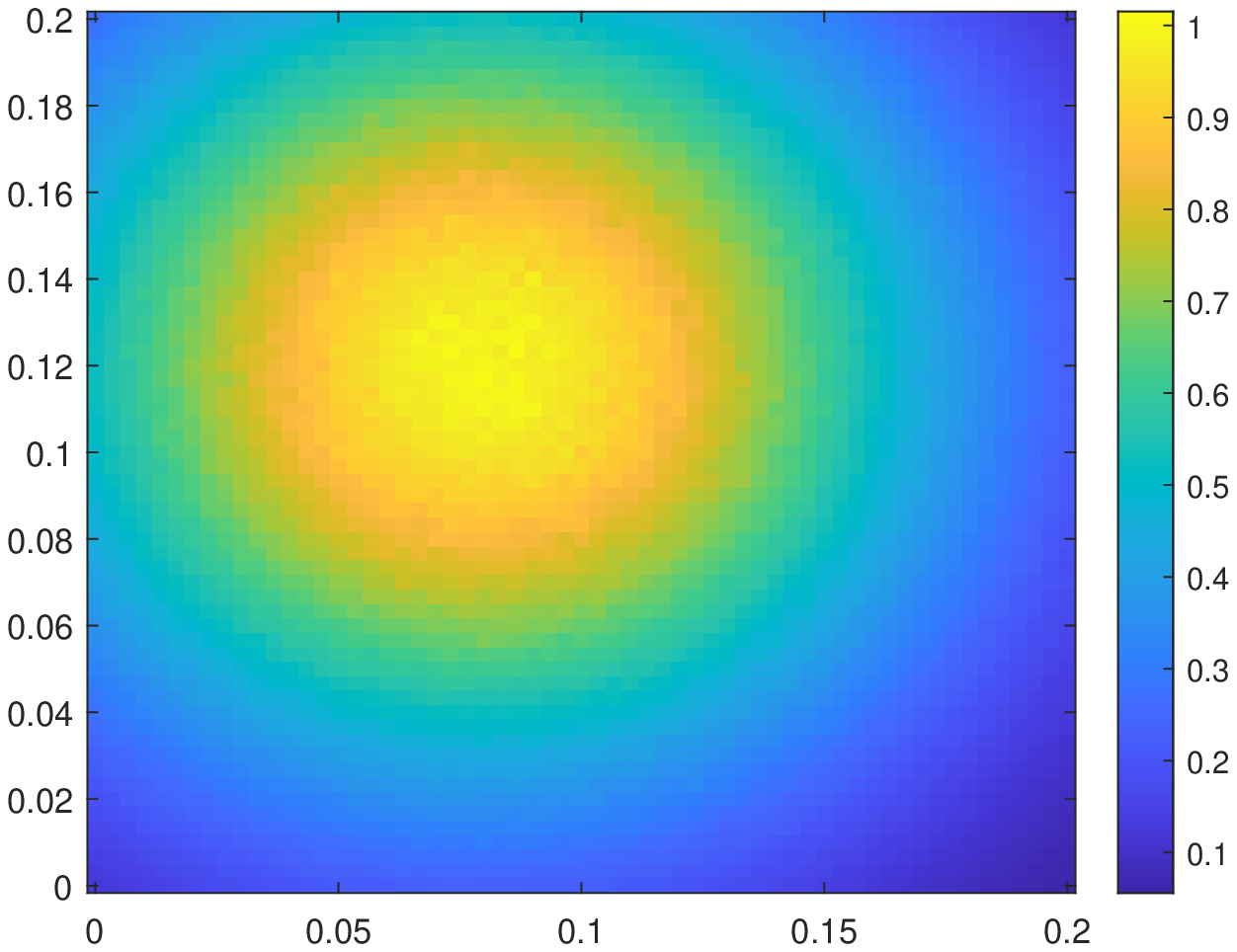}
\includegraphics[width=0.23\textwidth]{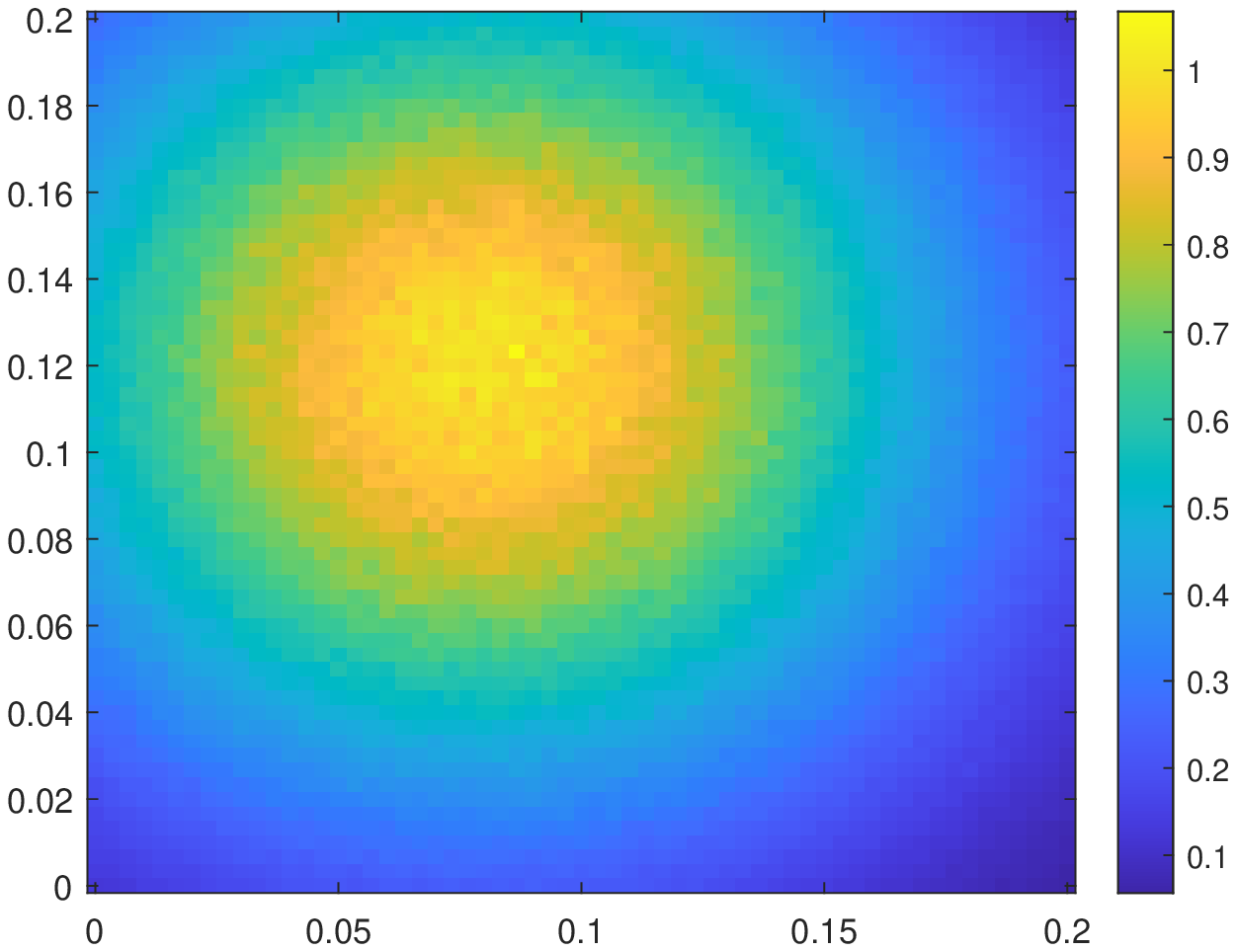}
\includegraphics[width=0.23\textwidth]{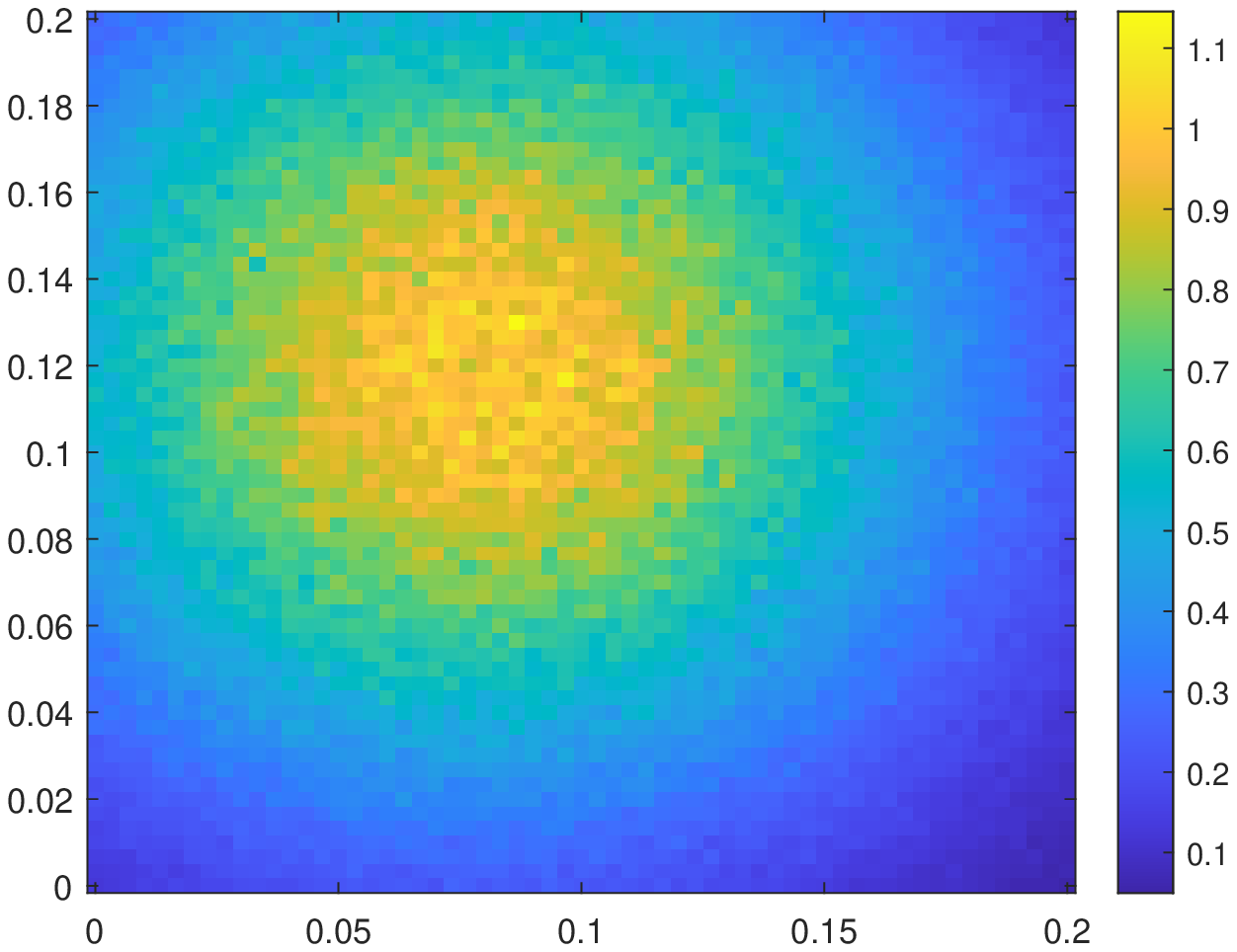}

\includegraphics[width=0.23\textwidth]{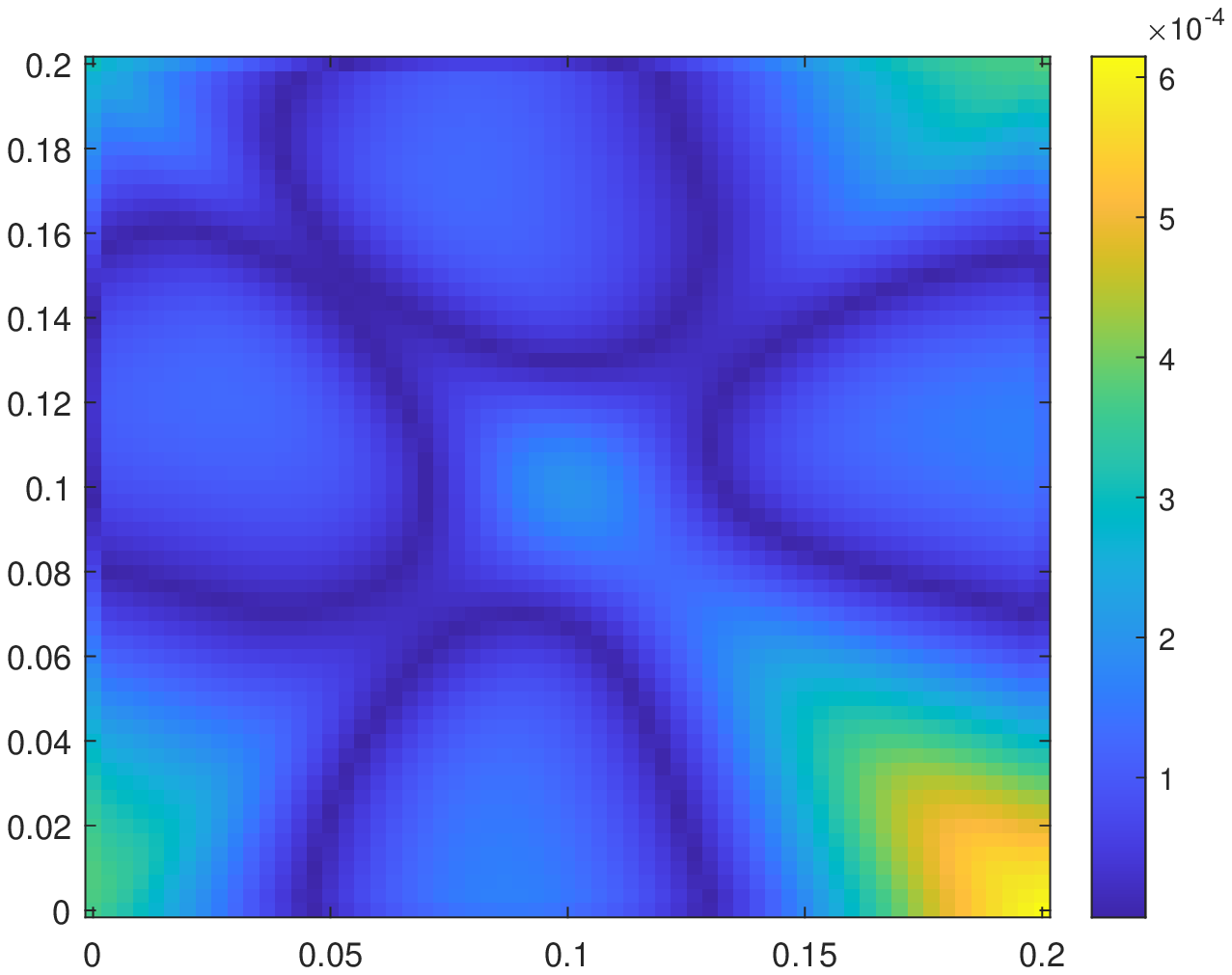}
\includegraphics[width=0.23\textwidth]{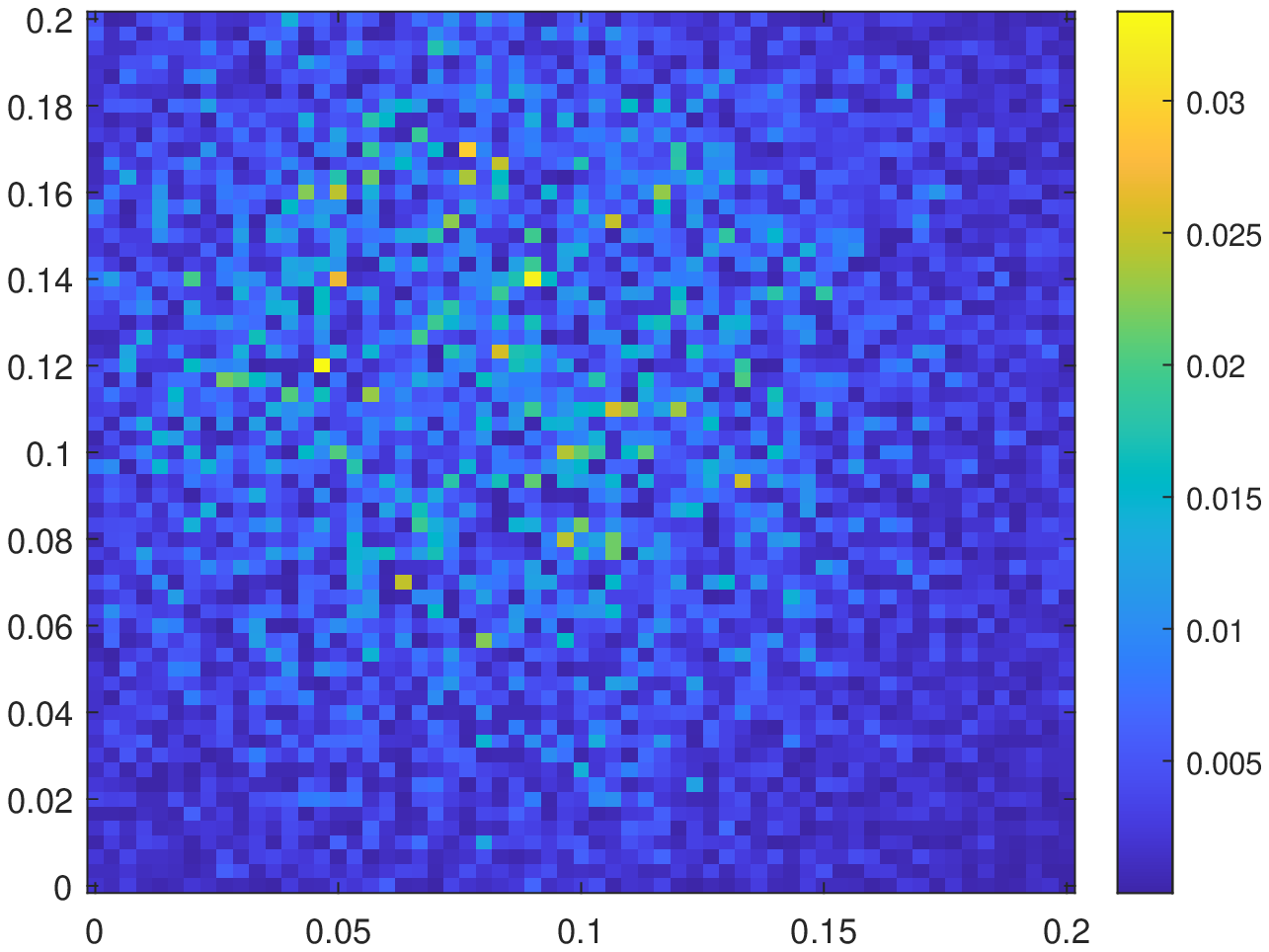}
\includegraphics[width=0.23\textwidth]{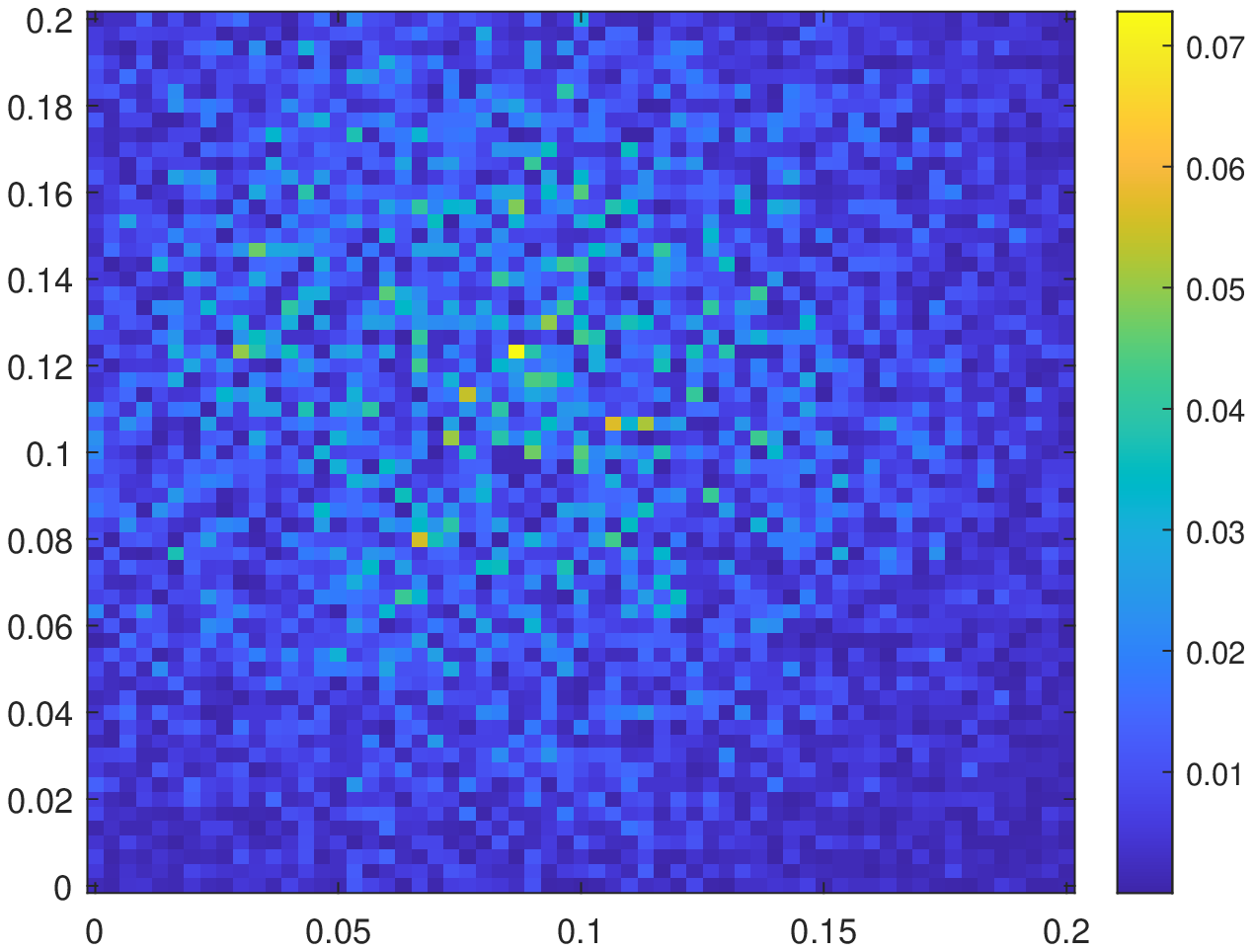}
\includegraphics[width=0.23\textwidth]{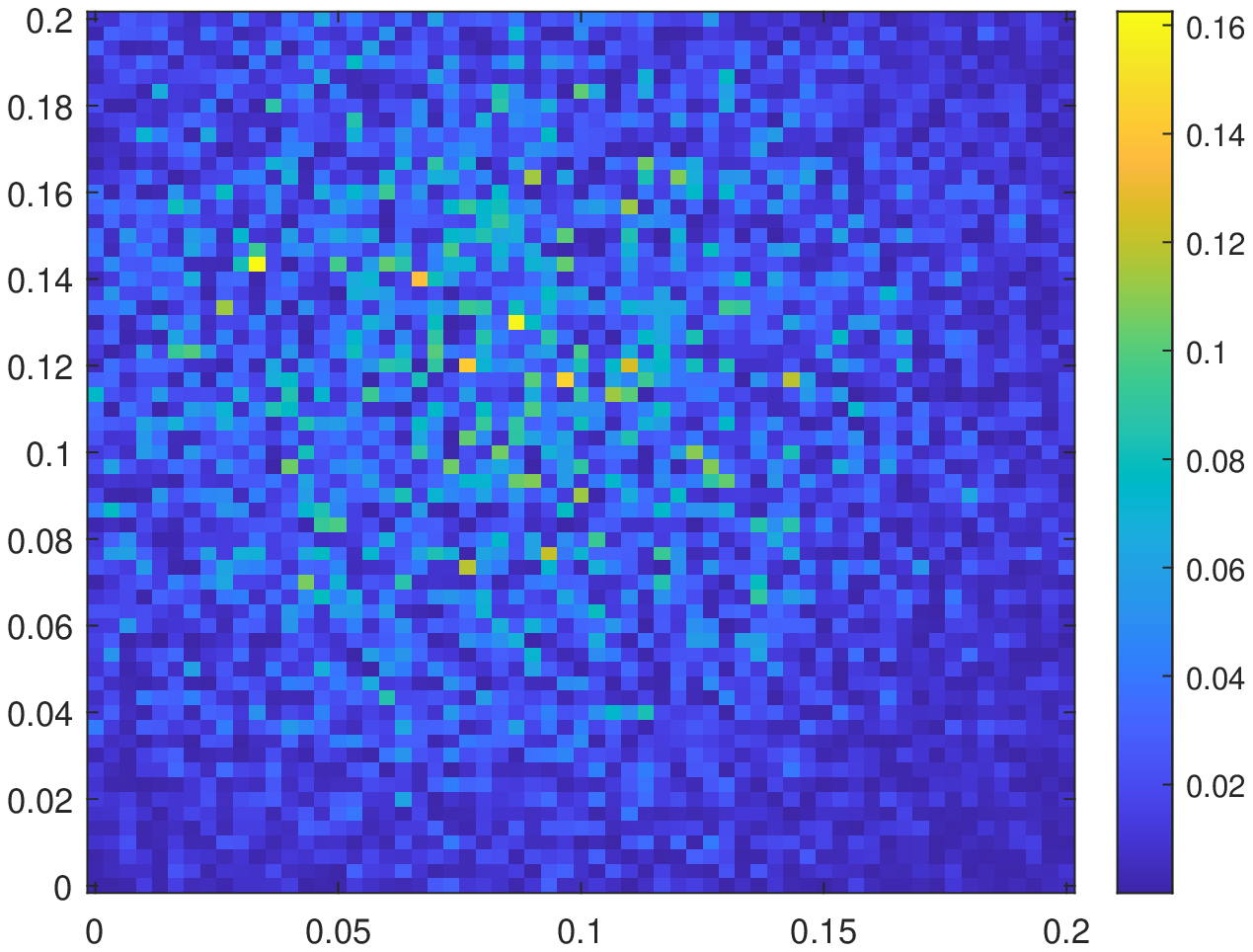}

\caption{Reconstructed $S_1$ using Neumann series. For the first row, 0\%, 1\%, 2\%, 5\% random noises are added to $H_{v_0}$. The relative $L^2$ errors of the reconstructions are 0.0268\%, 1.0682\%, 2.1759\%, 5.4680\%, respectively. The second row displays the corresponding differences between the ground truth and the reconstructions.}
\label{fig:Neumann1}
\vspace*{\floatsep}
\includegraphics[width=0.23\textwidth]{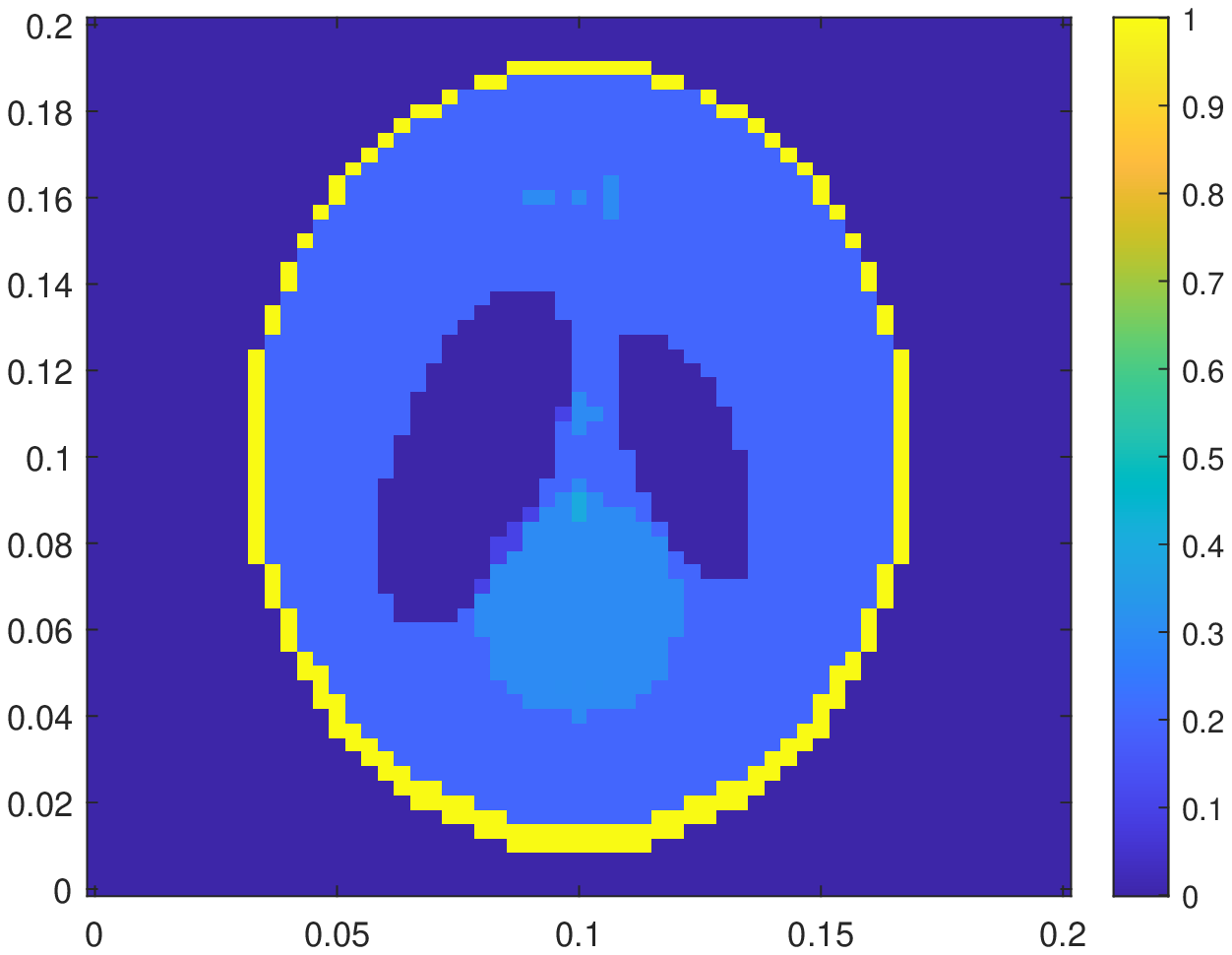}
\includegraphics[width=0.23\textwidth]{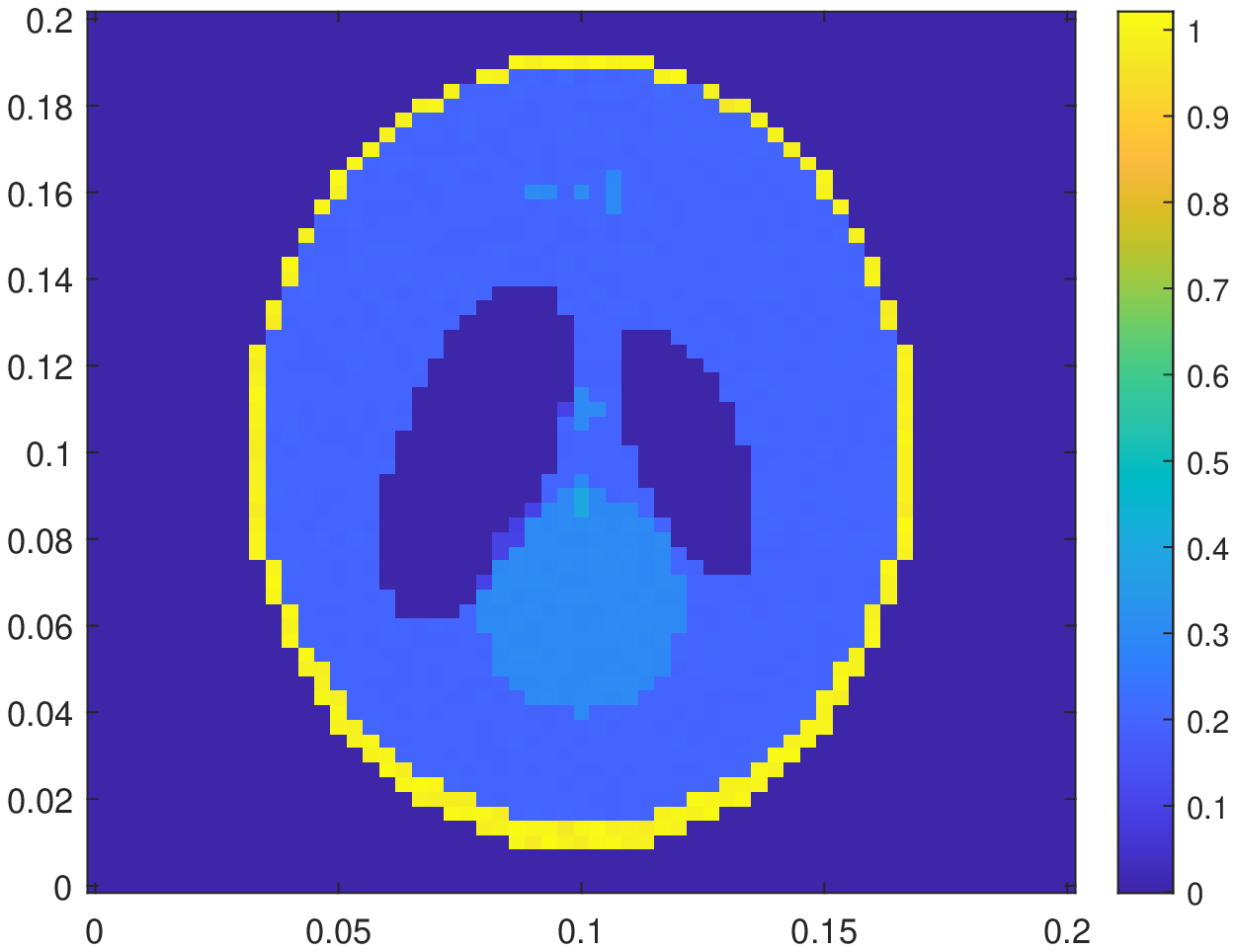}
\includegraphics[width=0.23\textwidth]{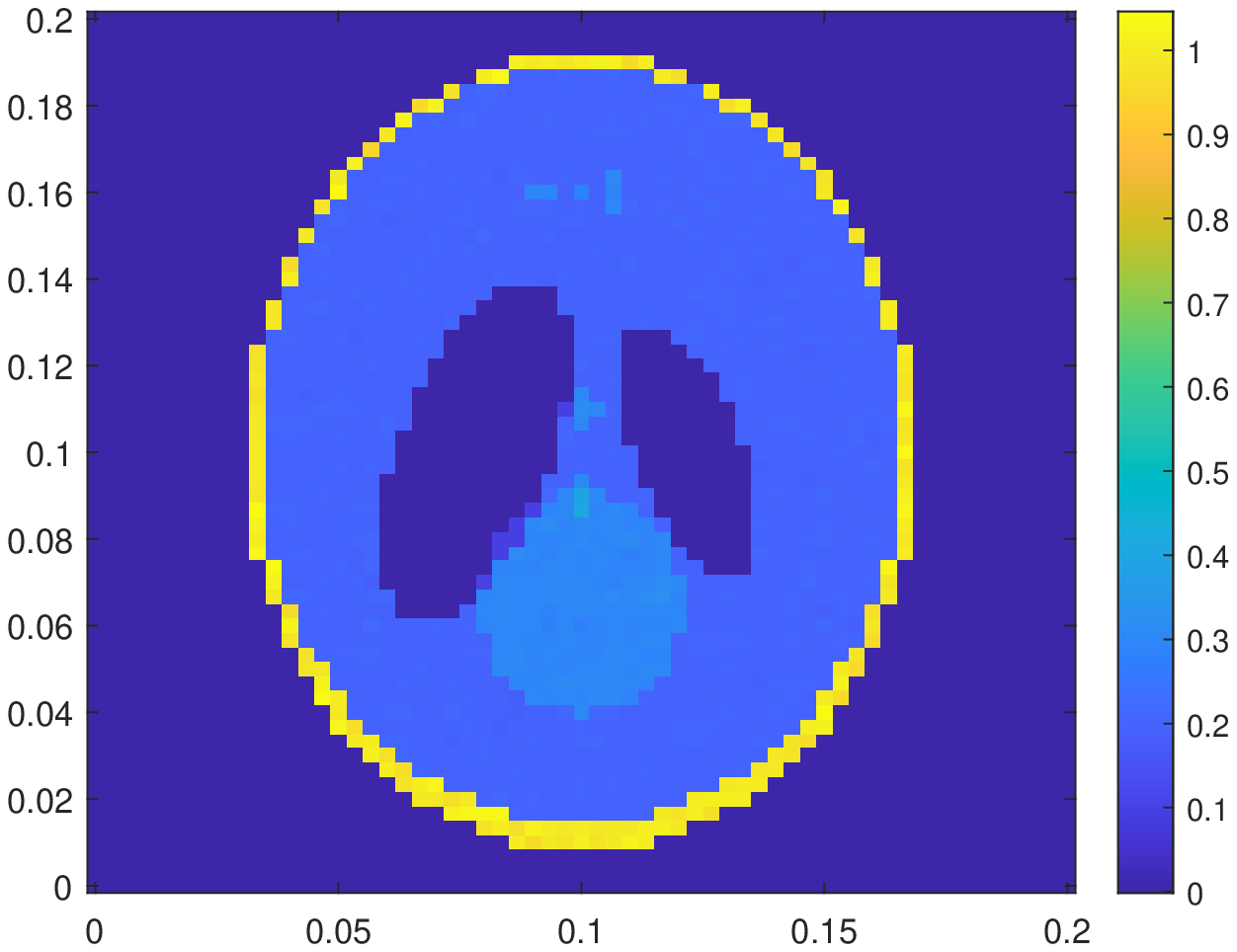}
\includegraphics[width=0.23\textwidth]{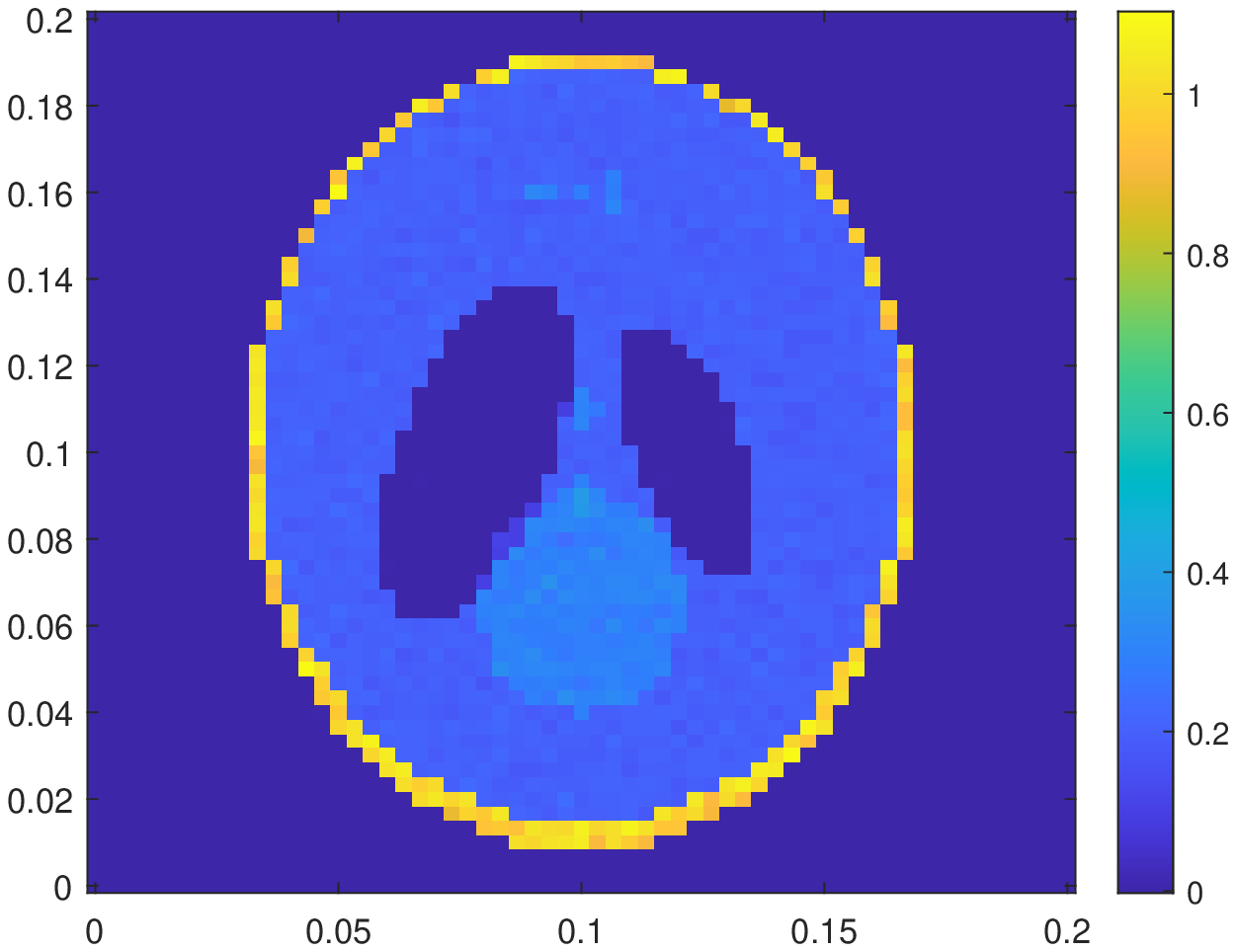}

\includegraphics[width=0.23\textwidth]{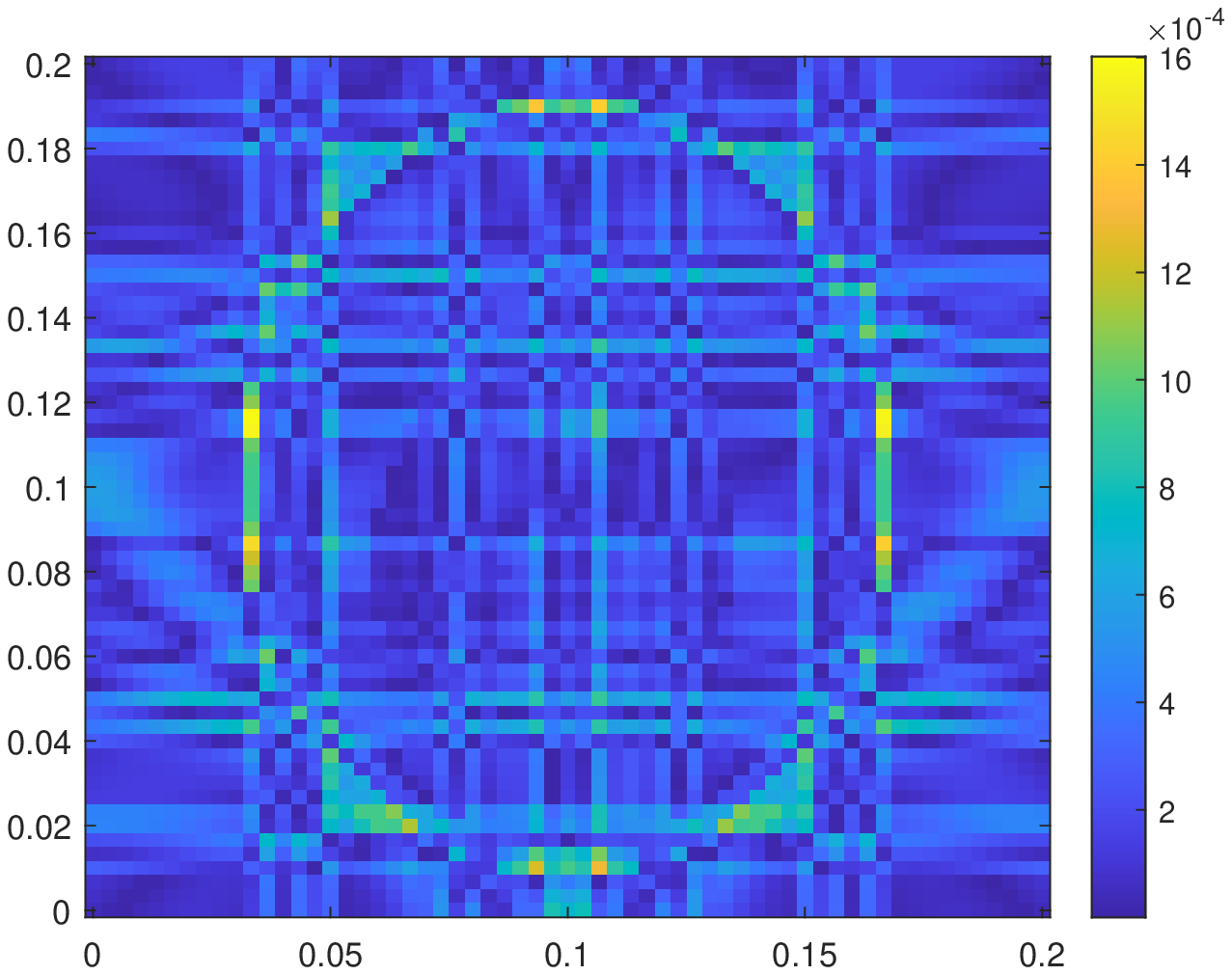}
\includegraphics[width=0.23\textwidth]{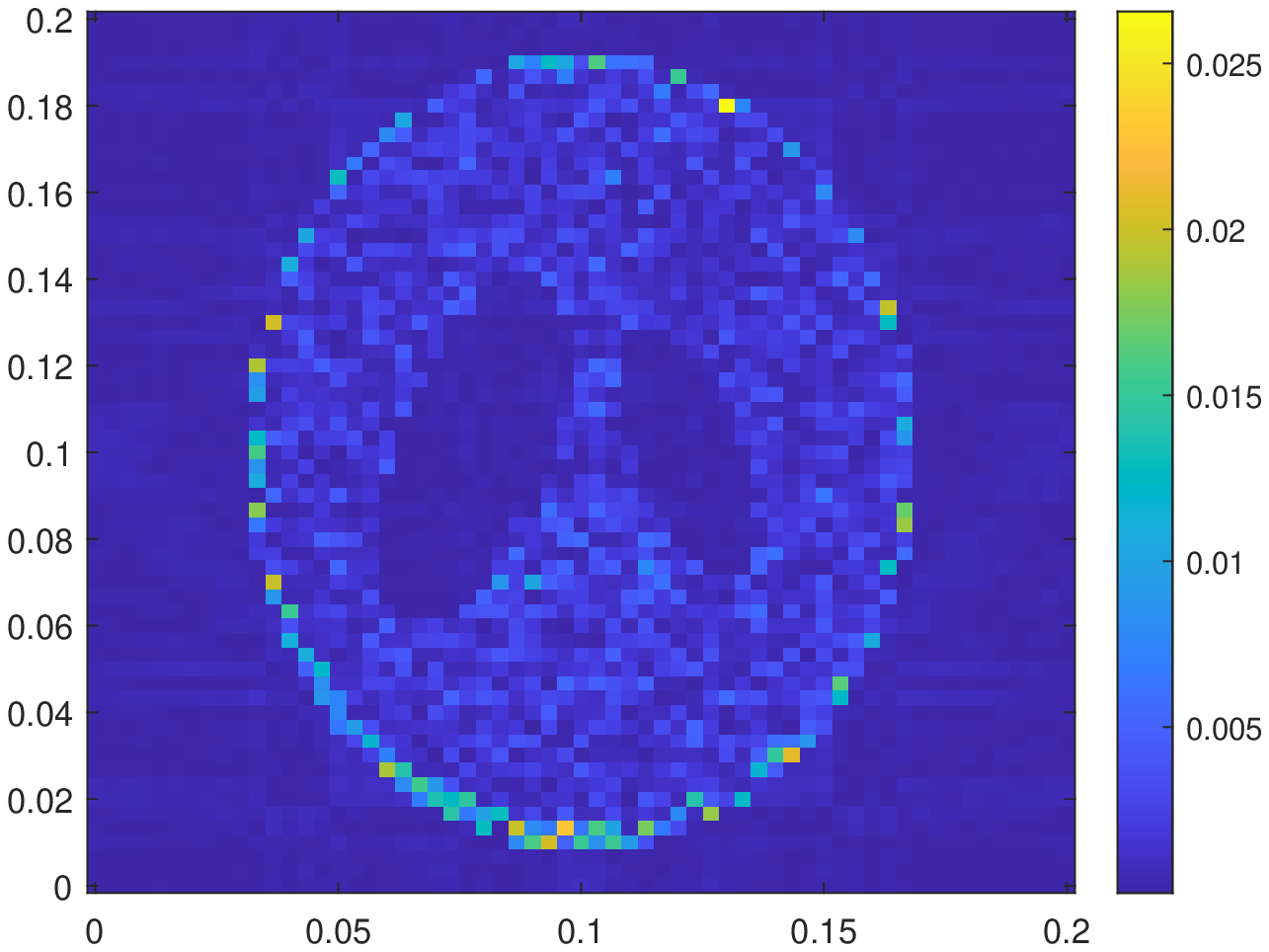}
\includegraphics[width=0.23\textwidth]{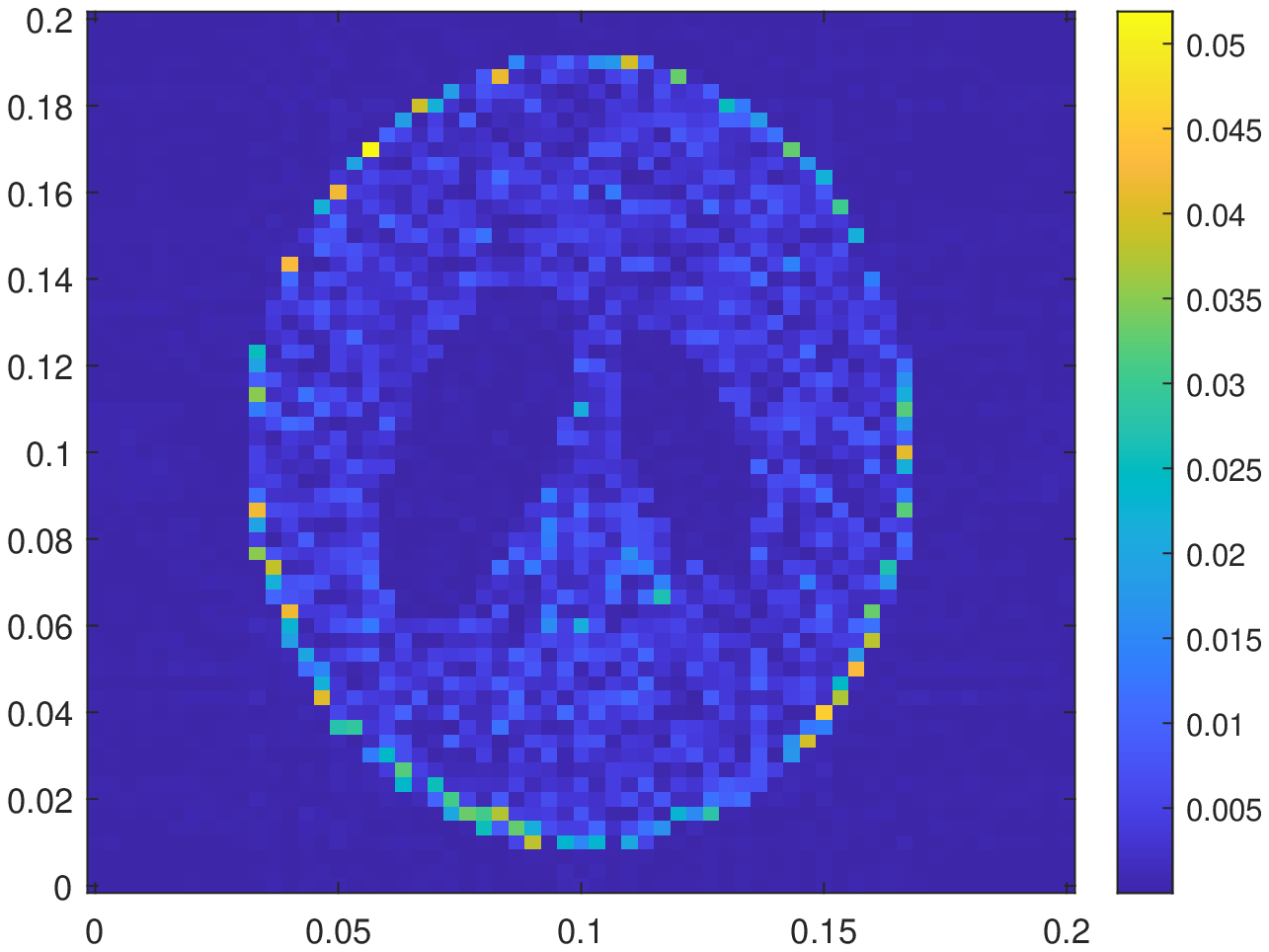}
\includegraphics[width=0.23\textwidth]{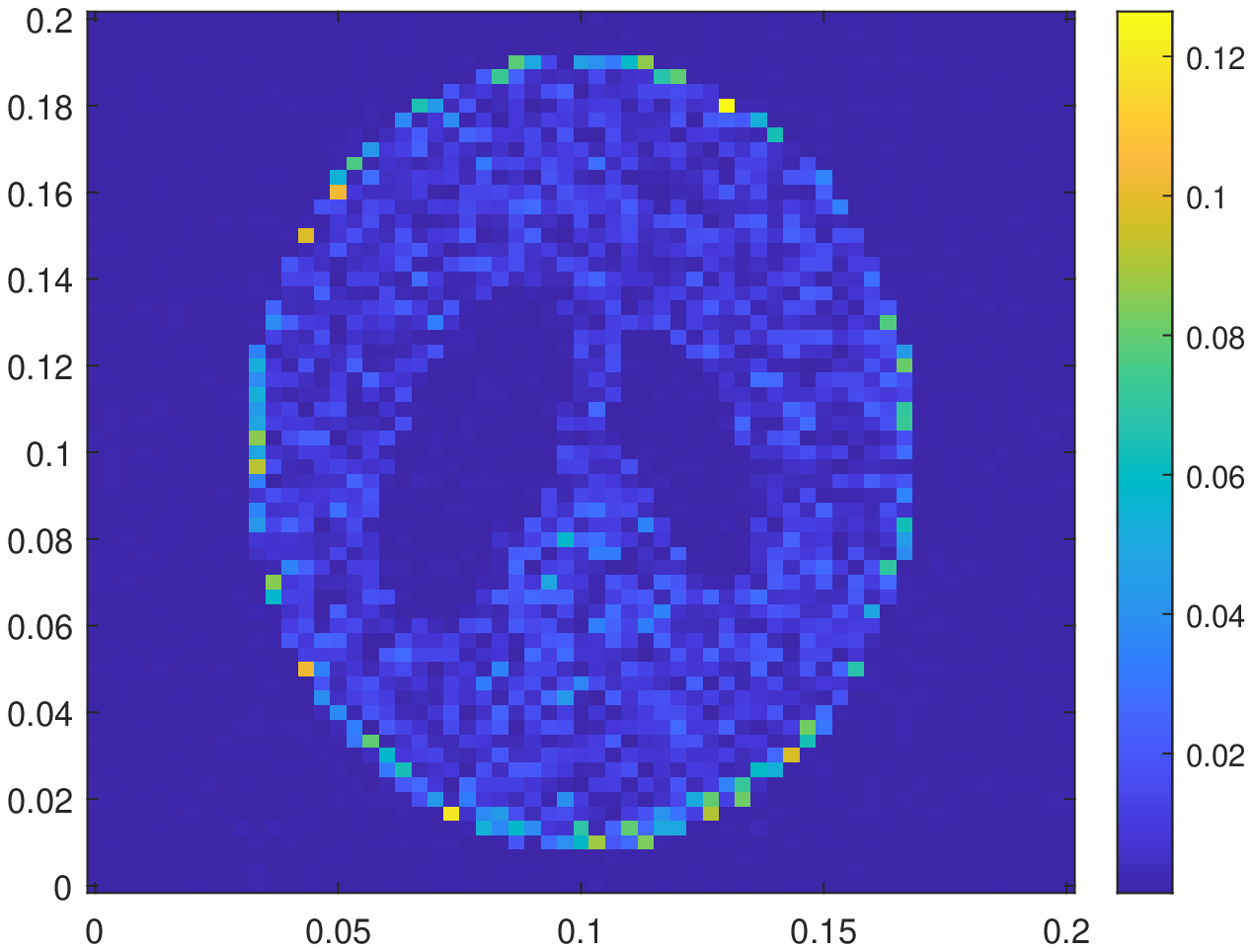}

\caption{Reconstructed $S_2$ using Neumann series. For the first row, 0\%, 1\%, 2\%, 5\% random noises are added to $H_{v_0}$. The relative $L^2$ errors of the reconstructions are 0.1383\%, 1.0152\%, 2.1301\%, 5.0305\%, respectively. The second row displays the corresponding differences between the ground truth and the reconstructions.}
\label{fig:Neumann2}
\end{figure}

We test the Neumann series inversion with a smooth source 
$$
S_1(x_1,x_2) = e^{-100[(x_1-0.08)^2+(x_2-0.12)^2]}
$$
and a discontinuous source $S_2 = $ Shepp-Logan phantom, see Figure~\ref{fig:Ncoef1}. The reconstructions with different levels of noises are illustrated in Figure~\ref{fig:Neumann1} and Figure~\ref{fig:Neumann2}, respectively.


 \begin{figure}[p]
 \centering
 \includegraphics[width=0.3\textwidth]{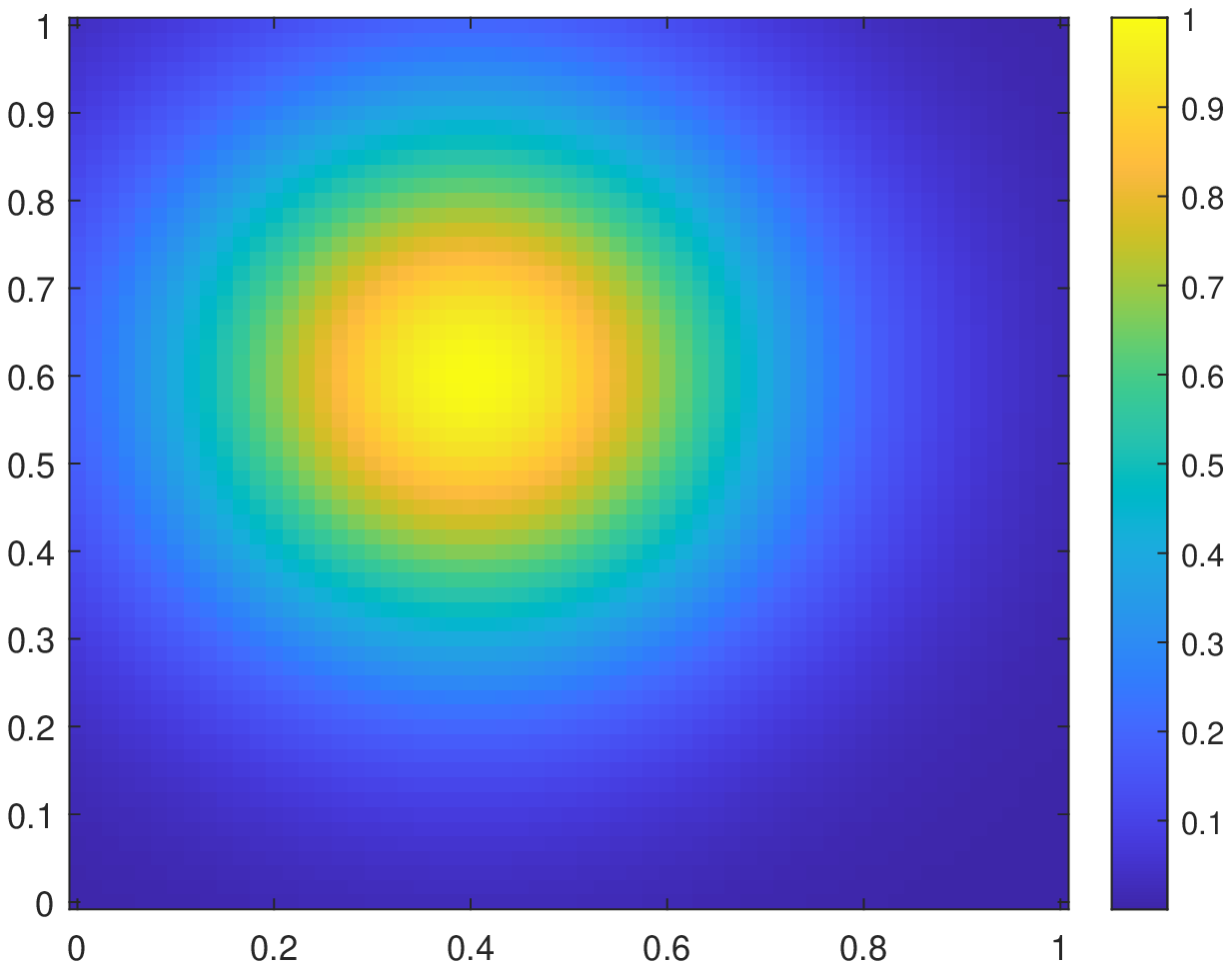}
 \includegraphics[width=0.3\textwidth]{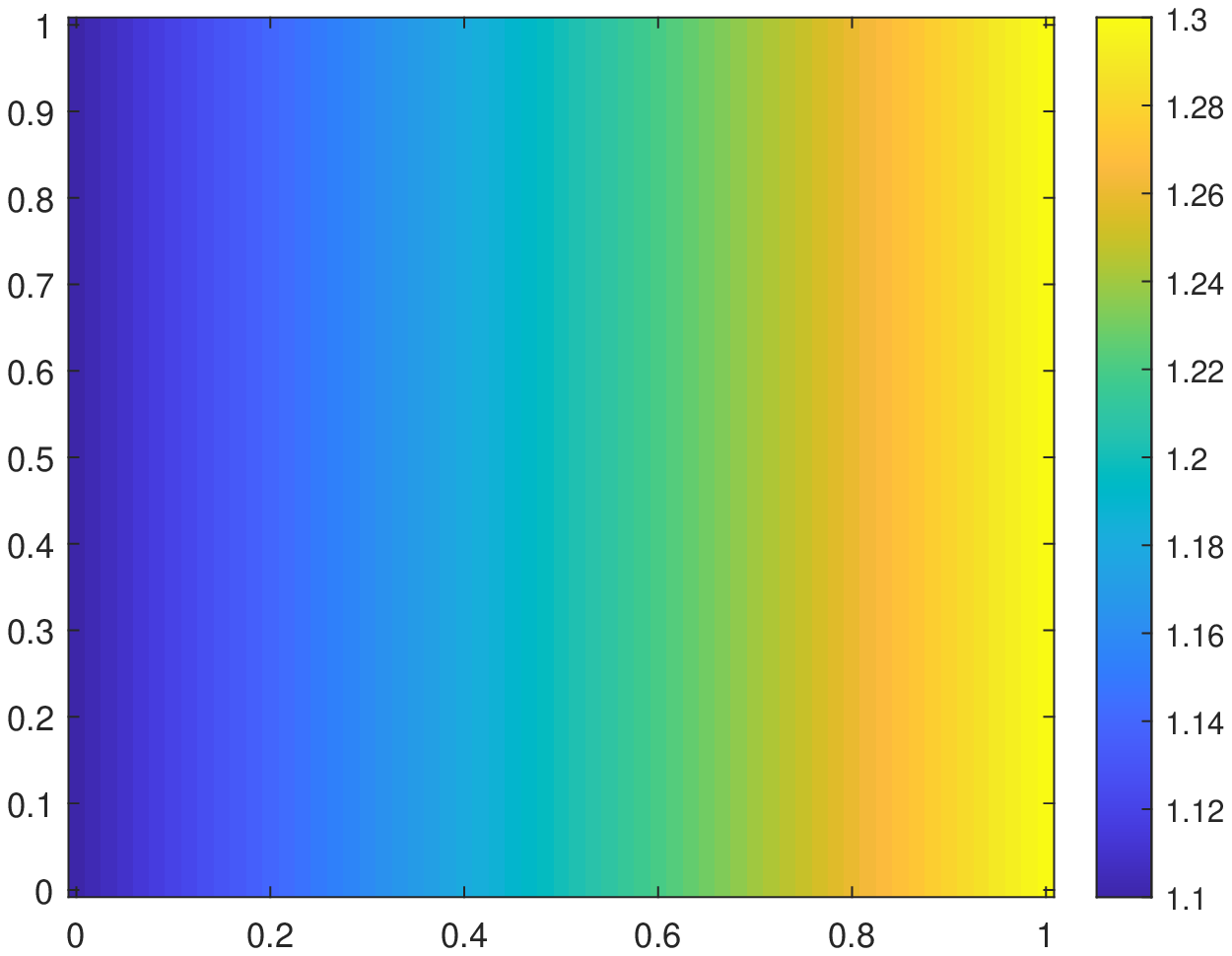}
 \caption{Left: source $S_3$. Right: attenuation coefficient $\sigma_2$}
 \label{fig:Wcoef1}
\vspace*{\floatsep}
 \includegraphics[width=0.23\textwidth]{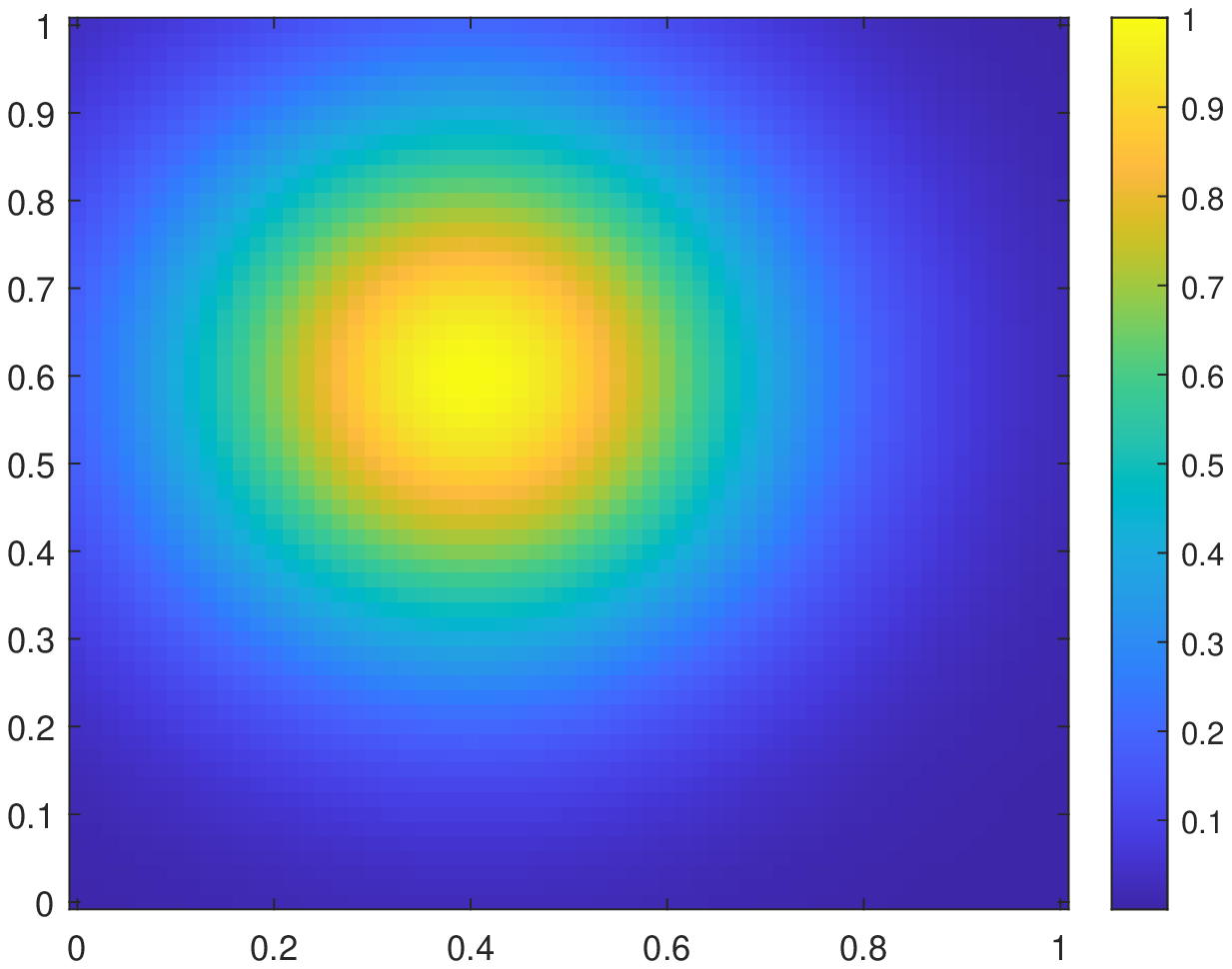}
 \includegraphics[width=0.23\textwidth]{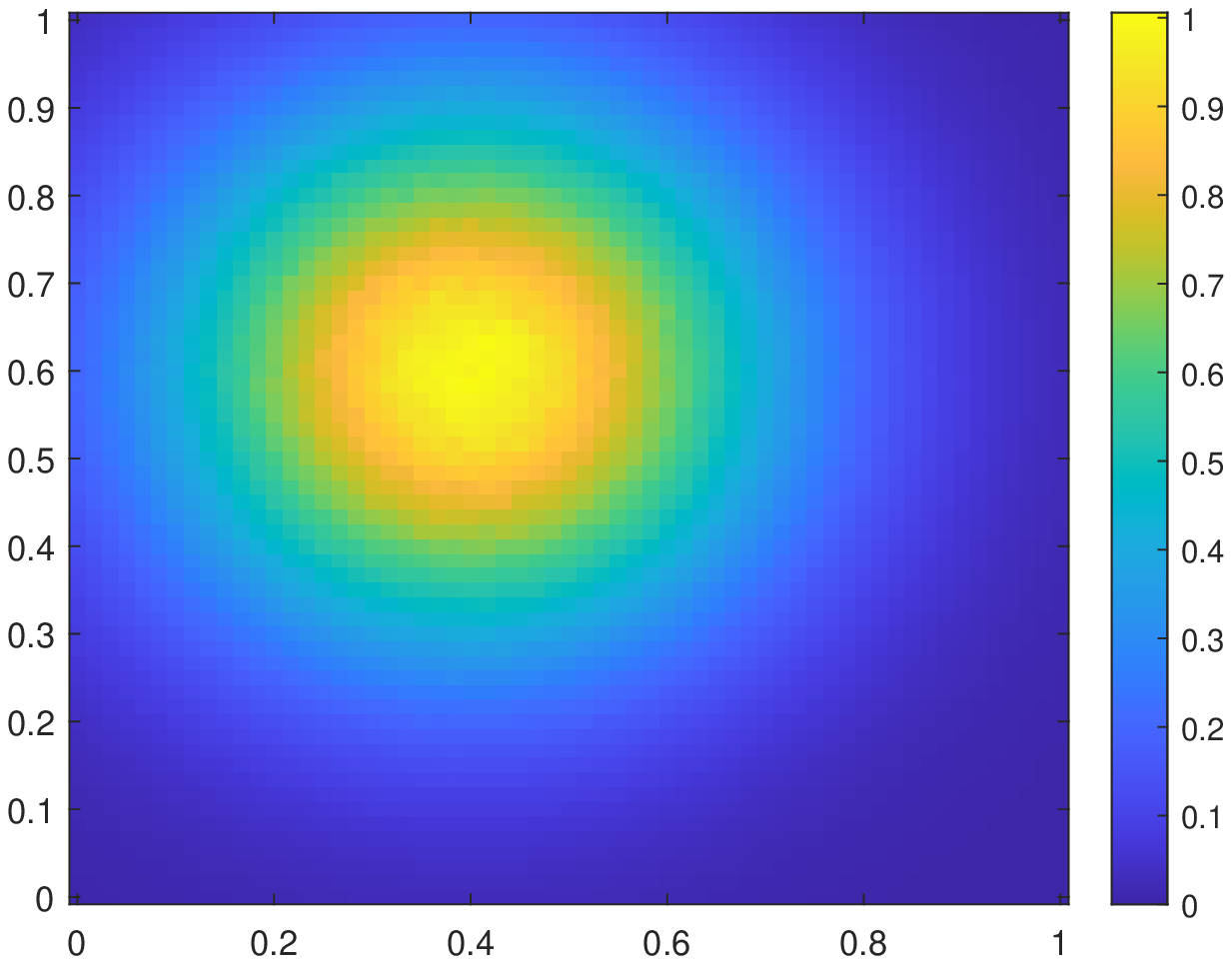}
 \includegraphics[width=0.23\textwidth]{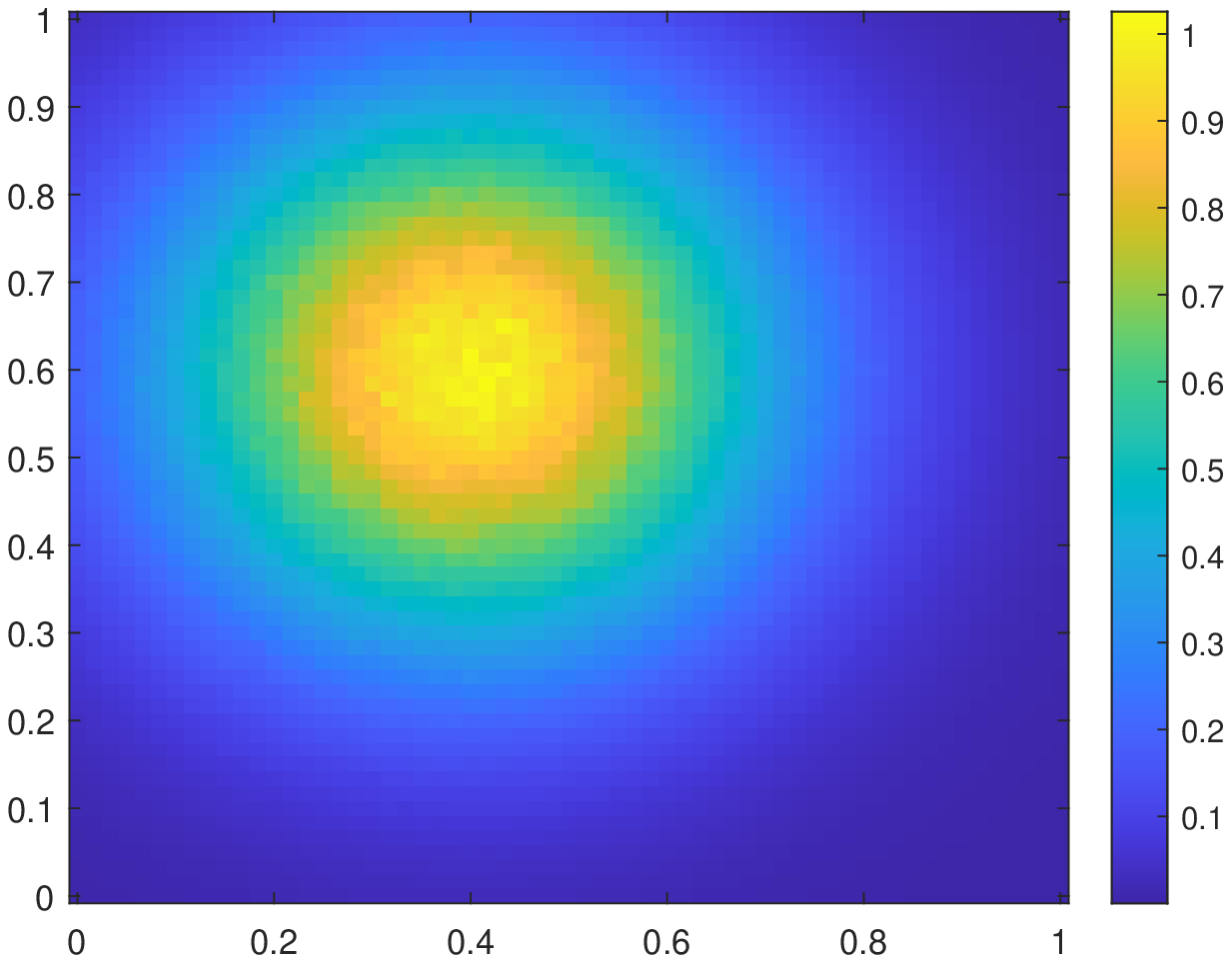}
 \includegraphics[width=0.23\textwidth]{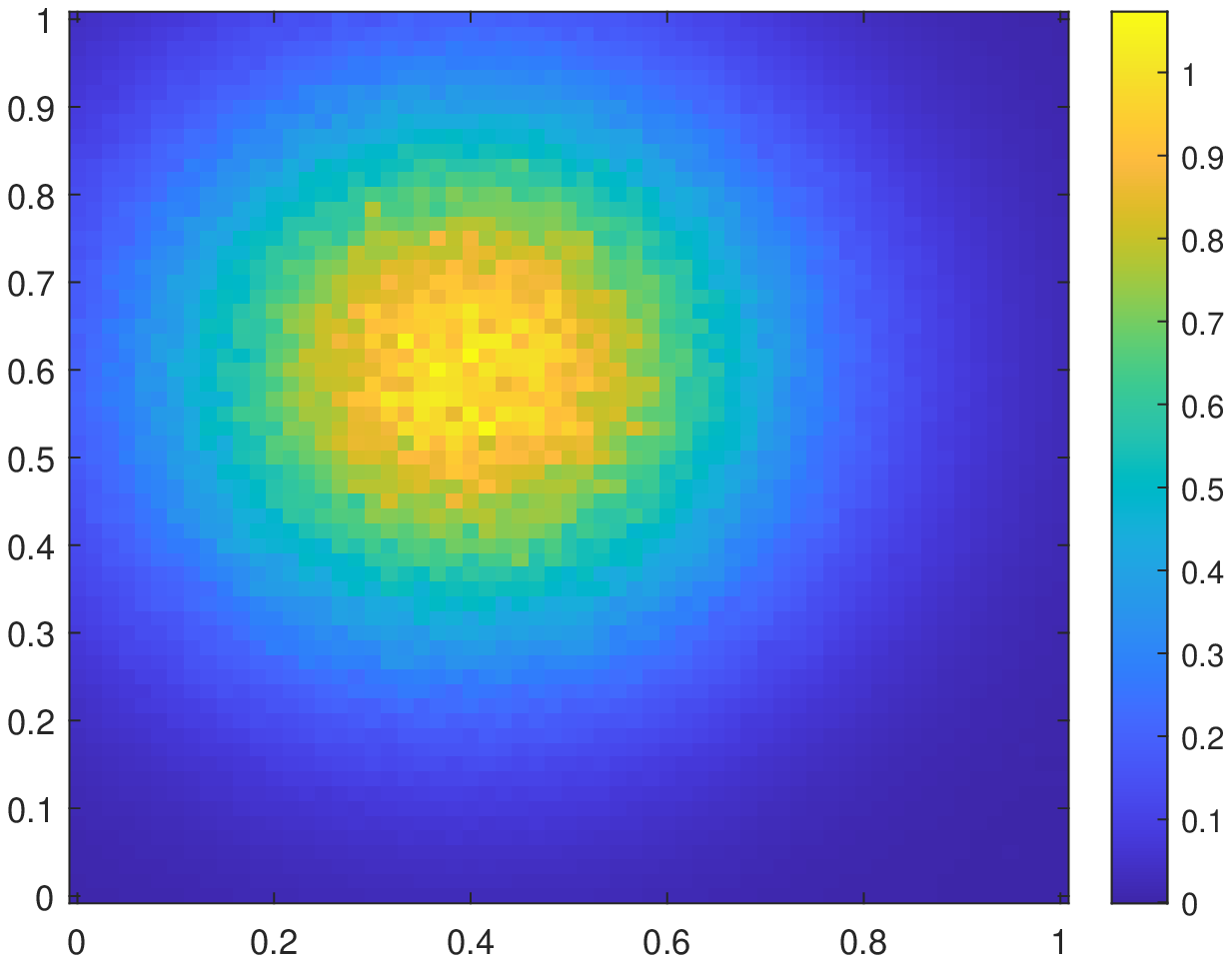}

 \includegraphics[width=0.23\textwidth]{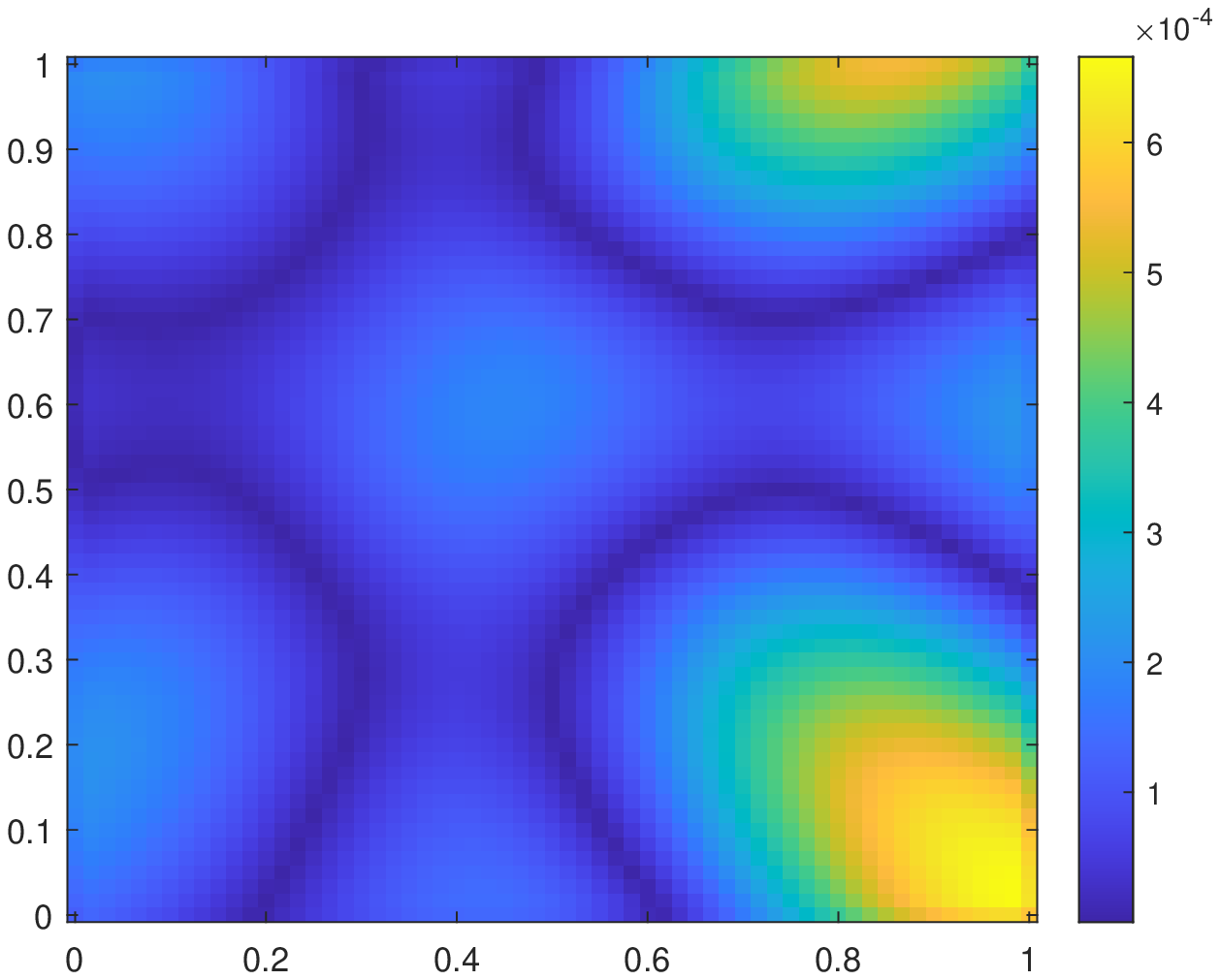}
 \includegraphics[width=0.23\textwidth]{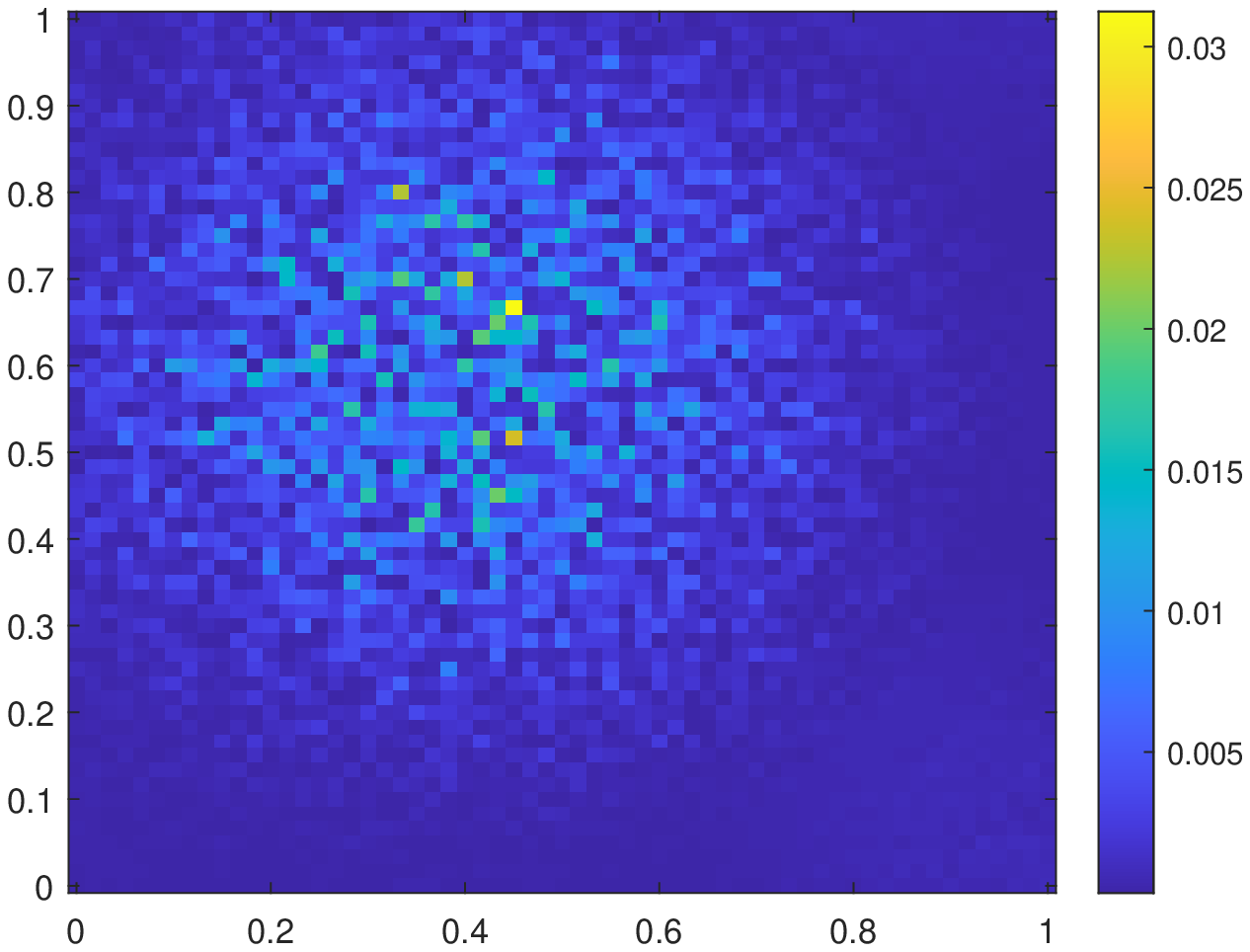}
 \includegraphics[width=0.23\textwidth]{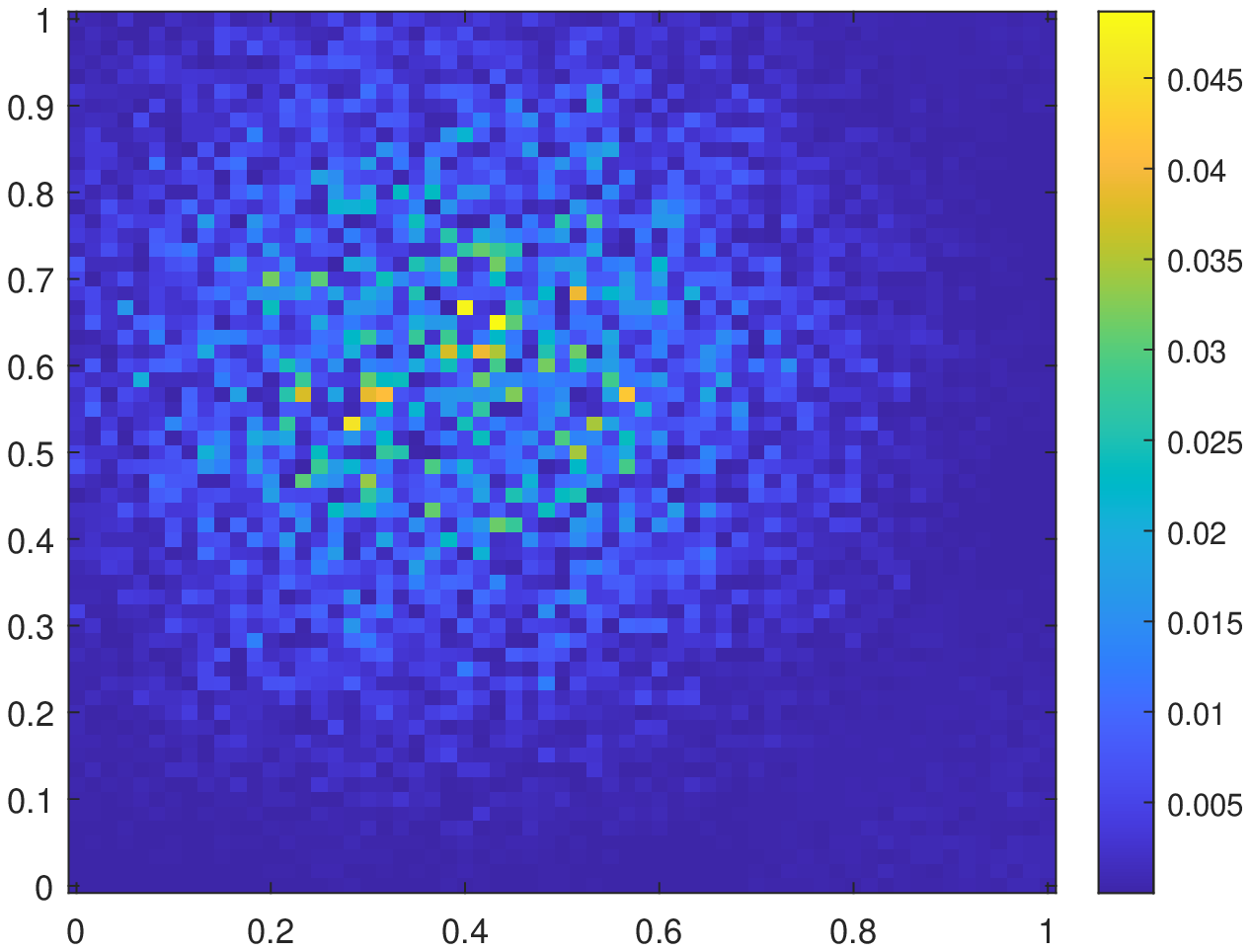}
 \includegraphics[width=0.23\textwidth]{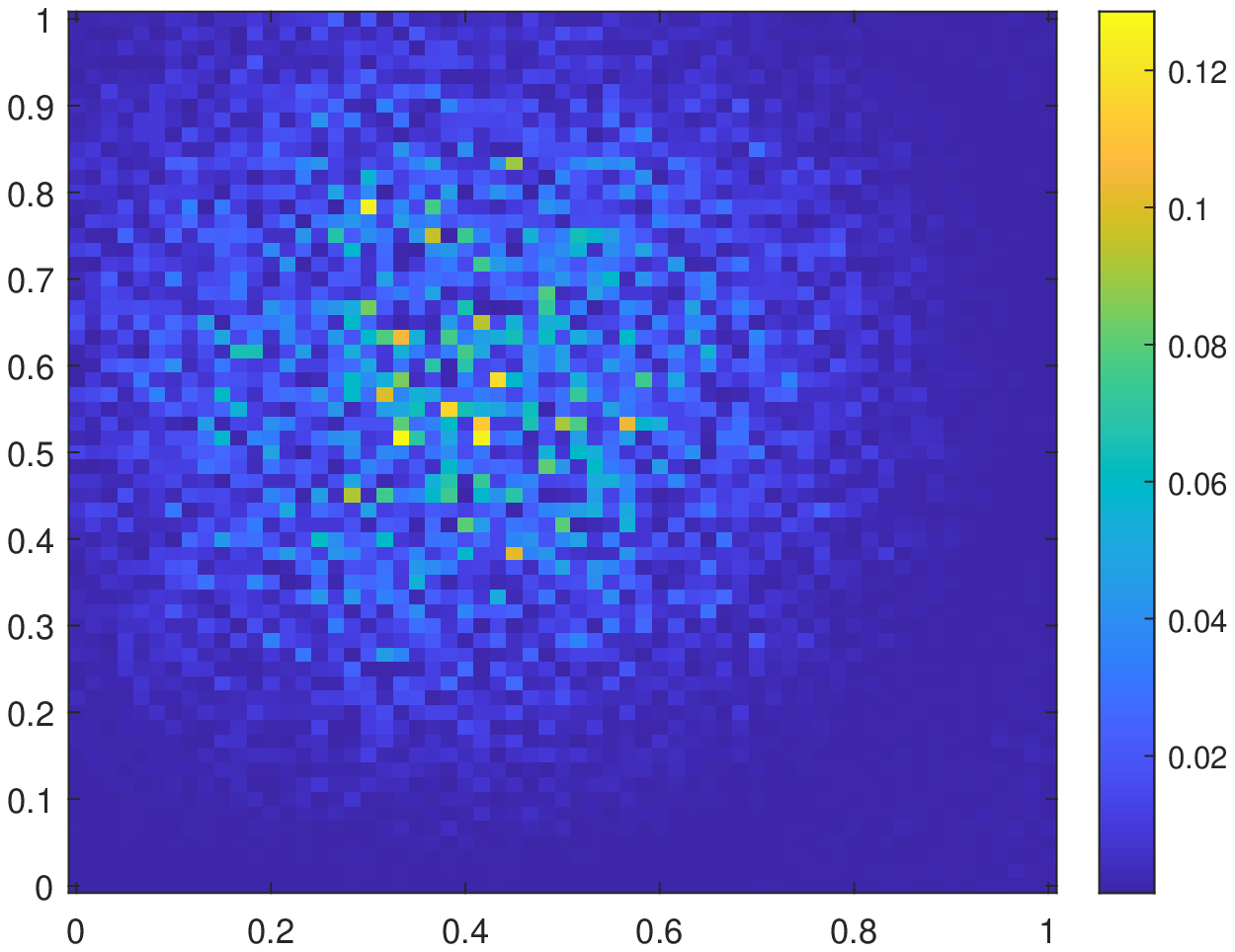}

 \caption{Reconstructed $S_3$ beyond the assumption of Theorem~\ref{thm:Neumann}. For the first row, 0\%, 1\%, 2\%, 5\% random noises are added to $H_{v_0}$. The relative $L^2$ errors of the reconstructions are 0.0526\%, 0.9410\%, 1.8096\%, 4.5969\%, respectively. 
 The second row  displays  the  corresponding  differences  between  the  ground  truth  and  the reconstructions.
 }
 \label{fig:Wrong1}
 \vspace*{\floatsep}
 \includegraphics[width=0.23\textwidth]{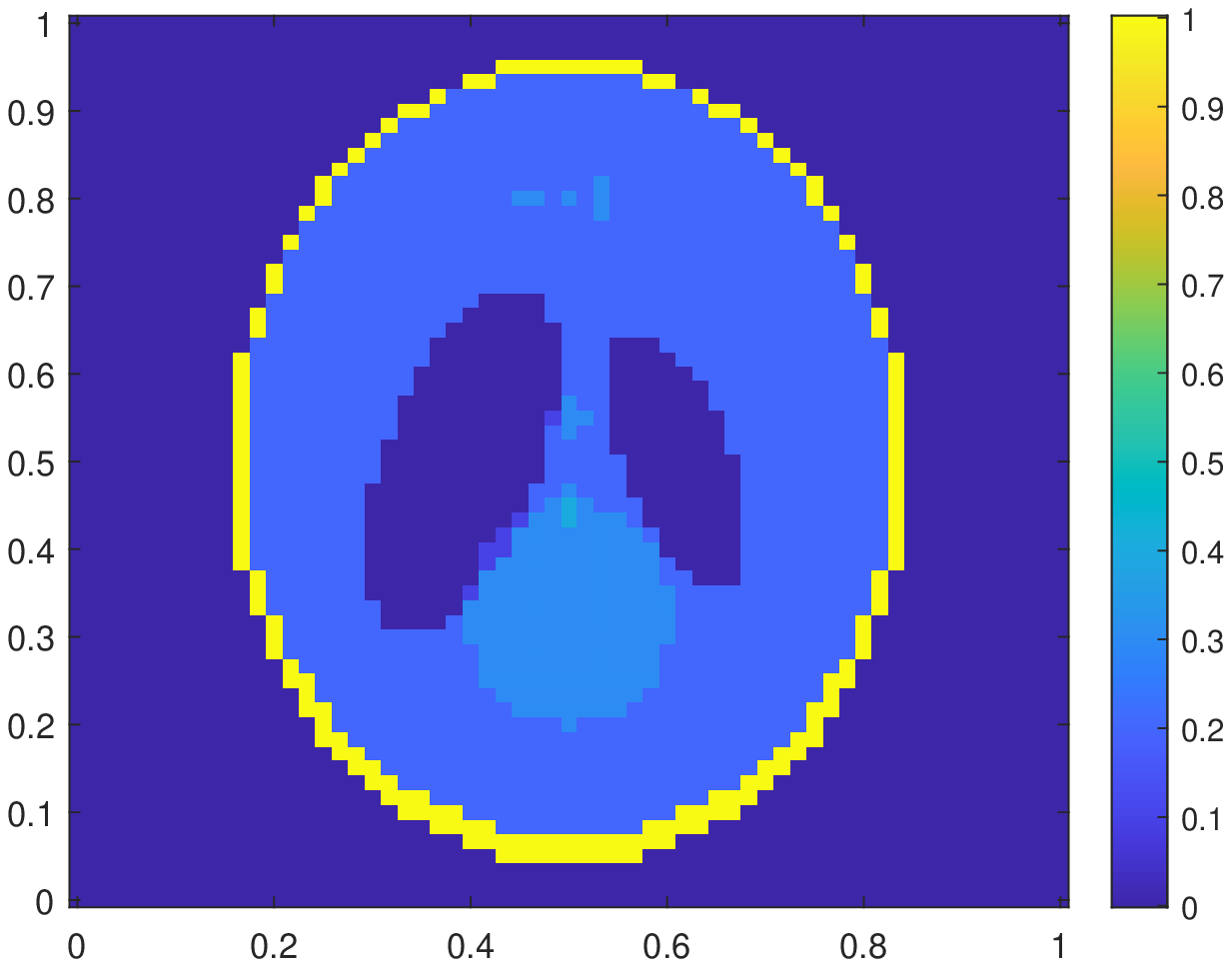}
 \includegraphics[width=0.23\textwidth]{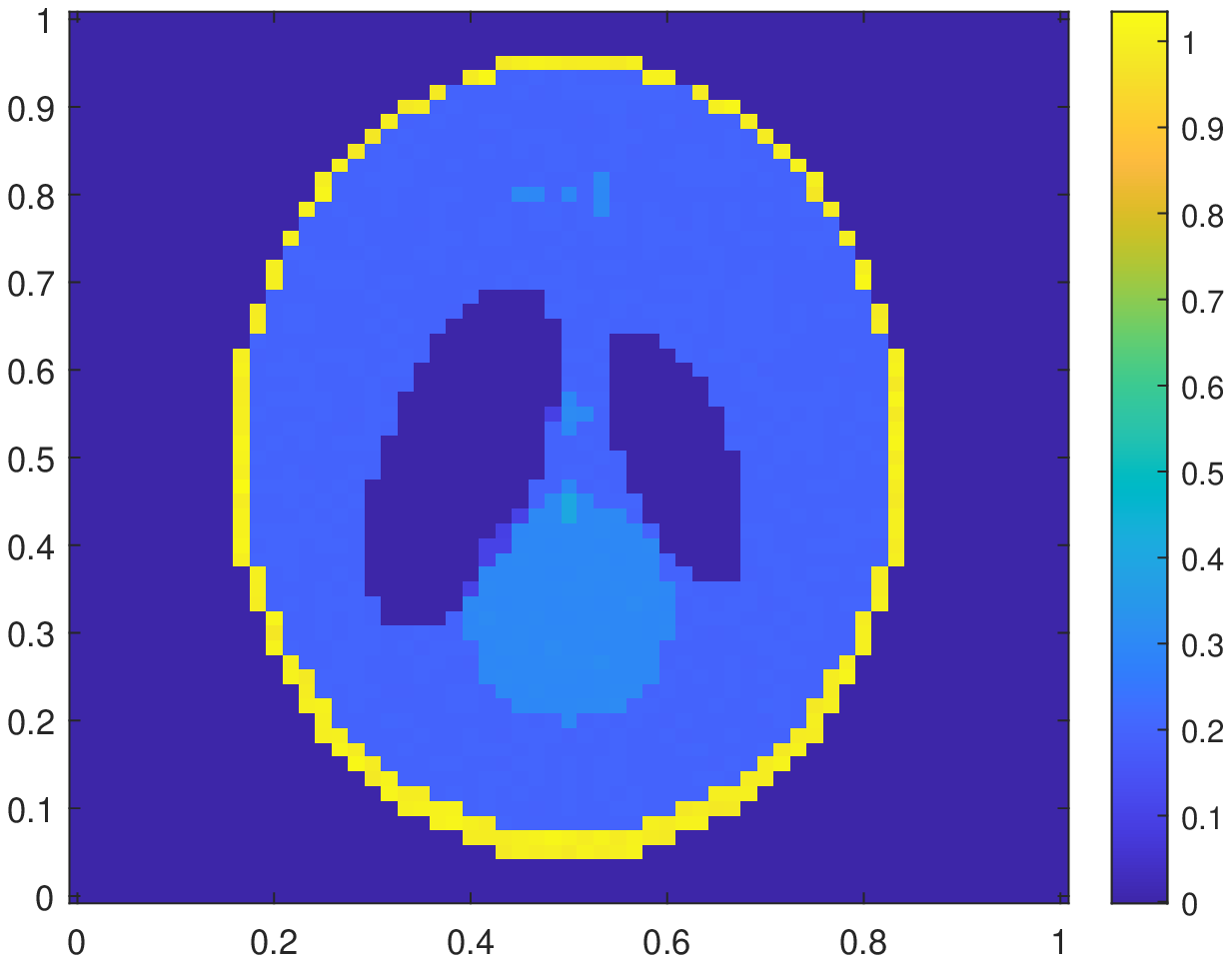}
 \includegraphics[width=0.23\textwidth]{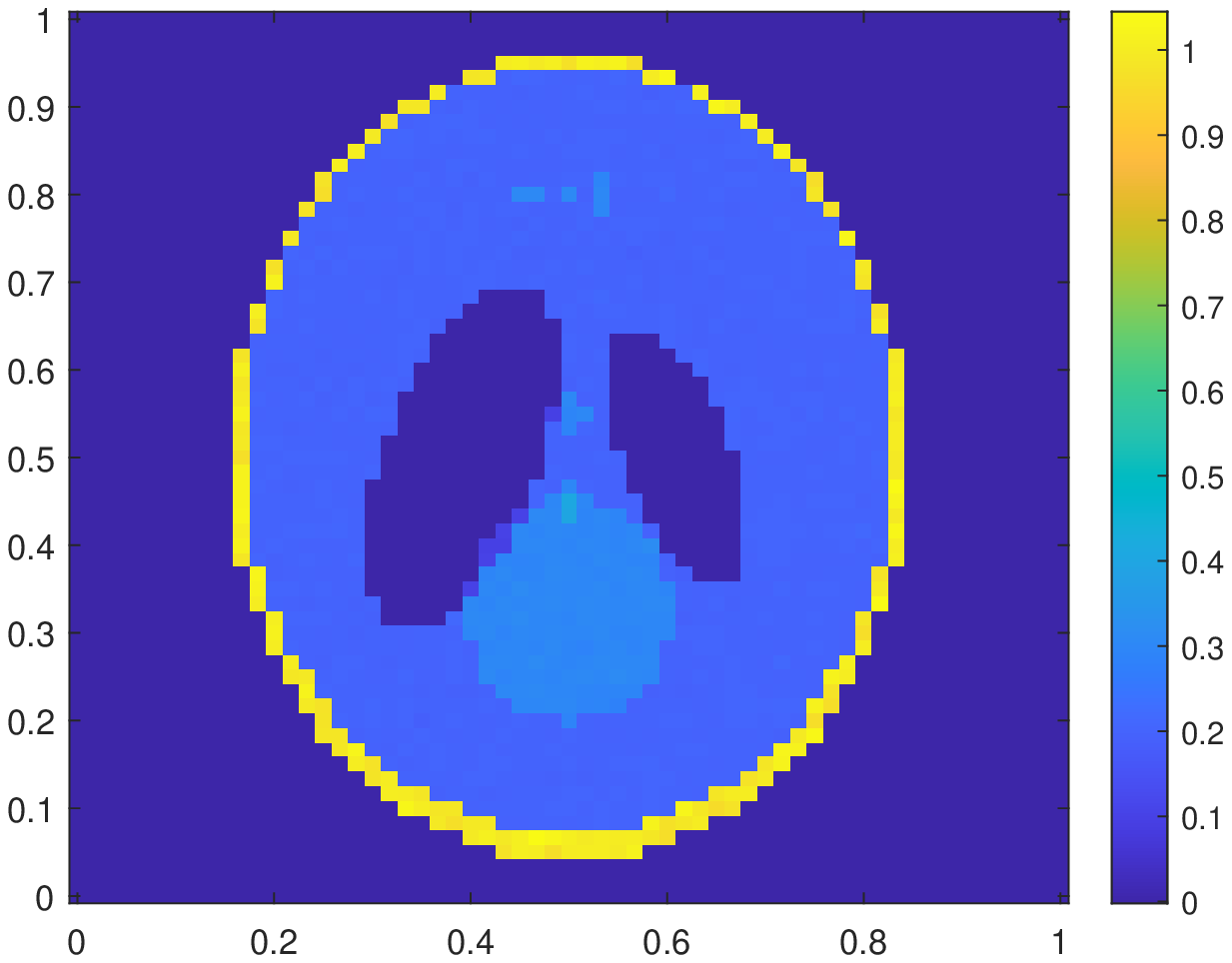}
 \includegraphics[width=0.23\textwidth]{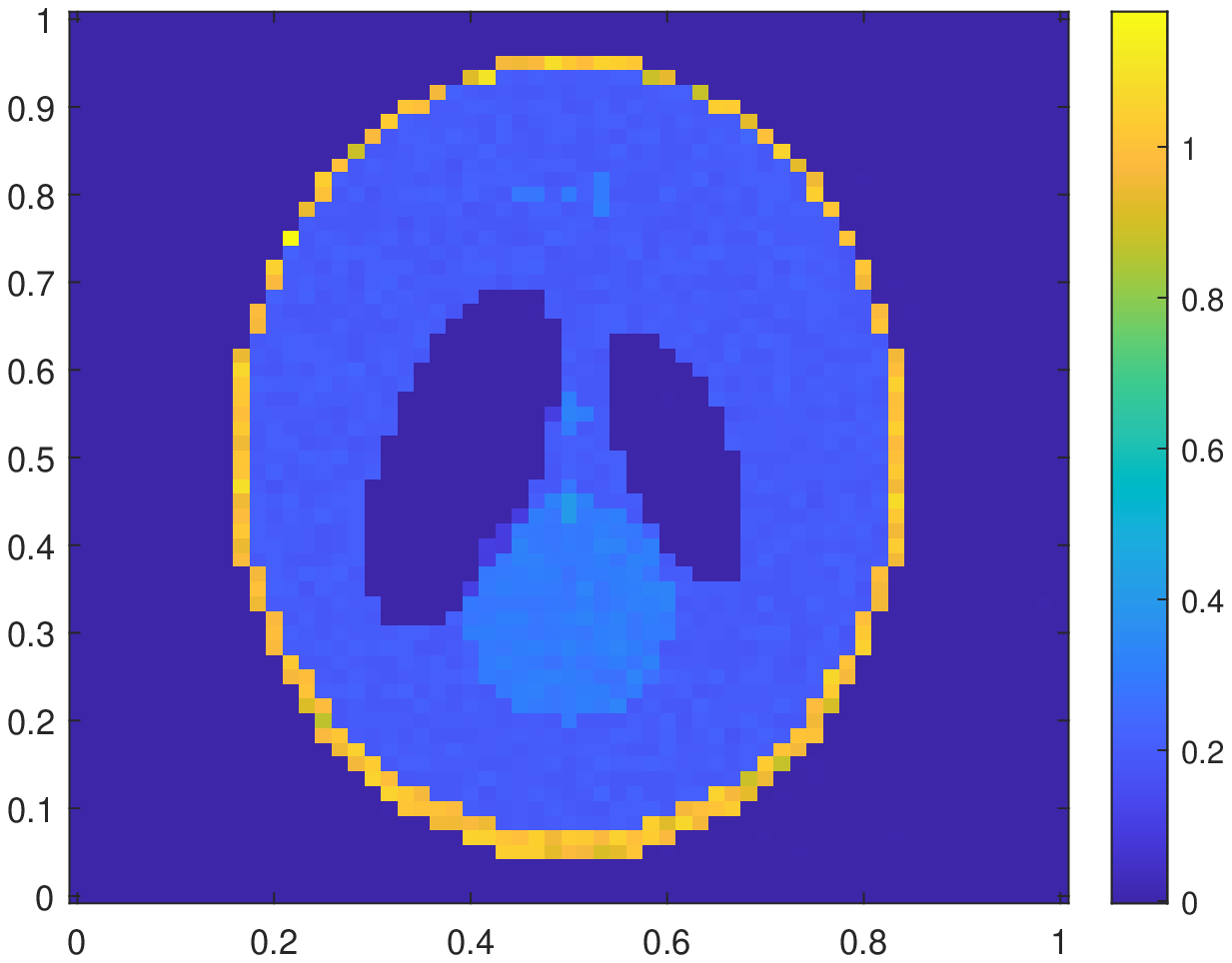}

 \includegraphics[width=0.23\textwidth]{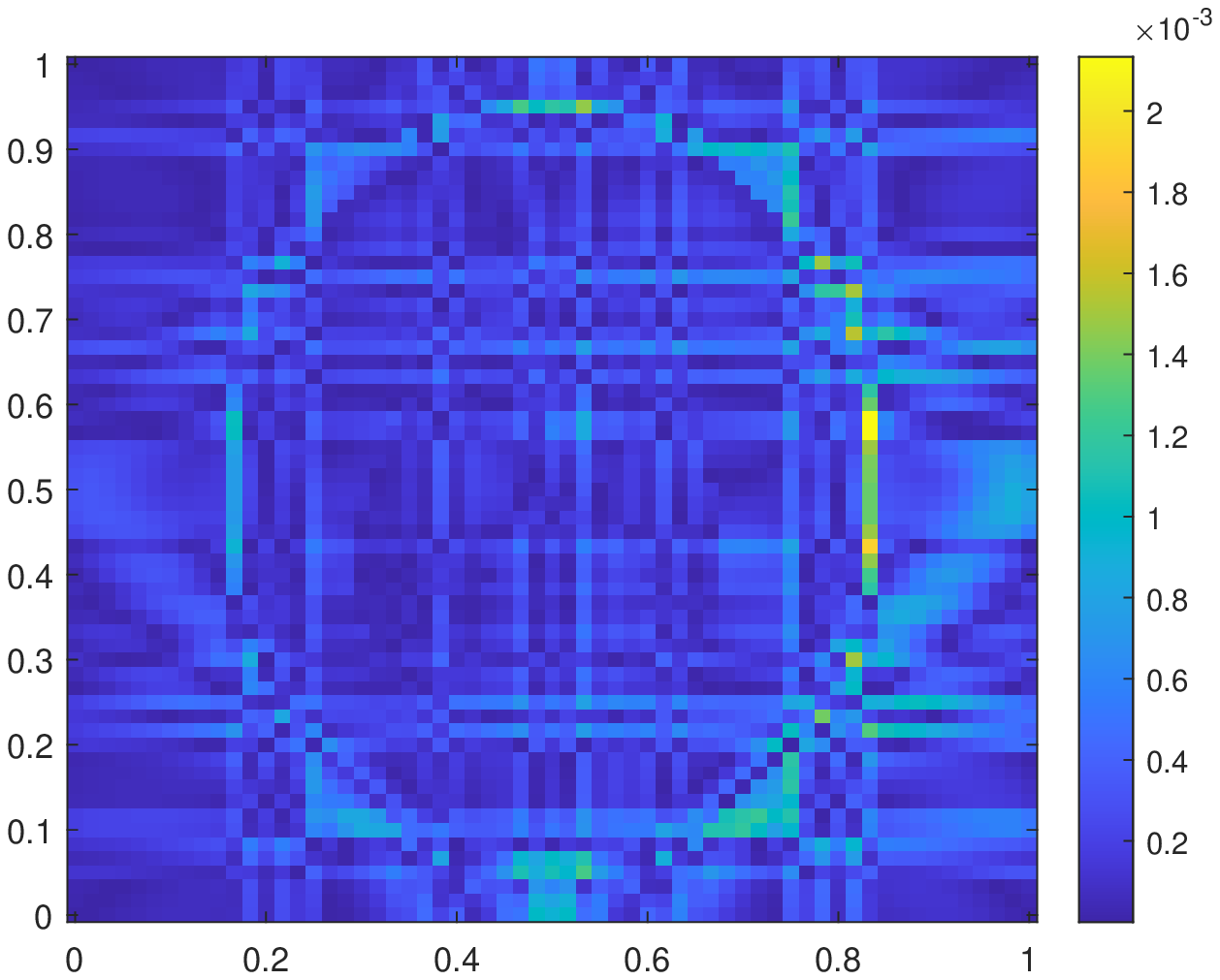}
 \includegraphics[width=0.23\textwidth]{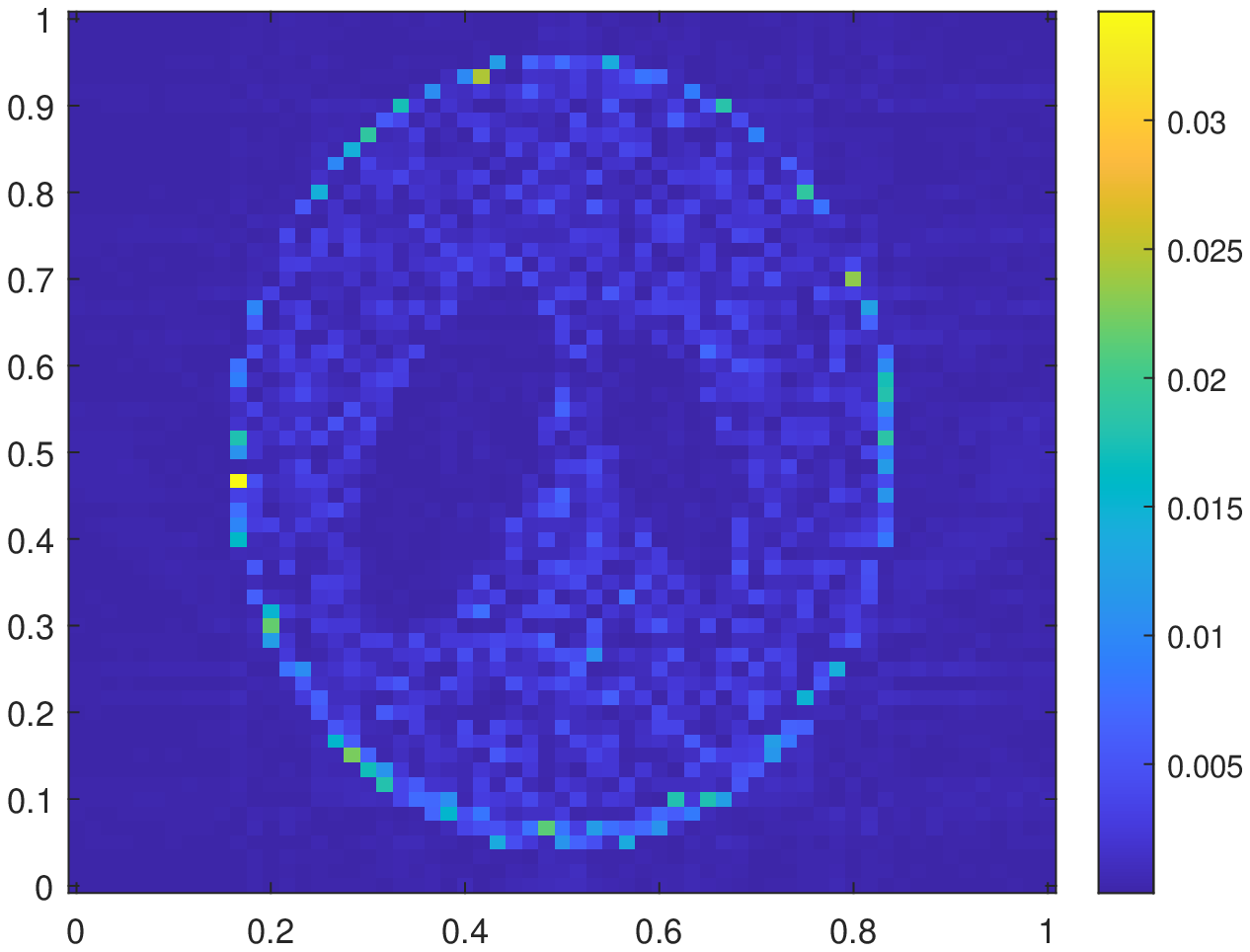}
 \includegraphics[width=0.23\textwidth]{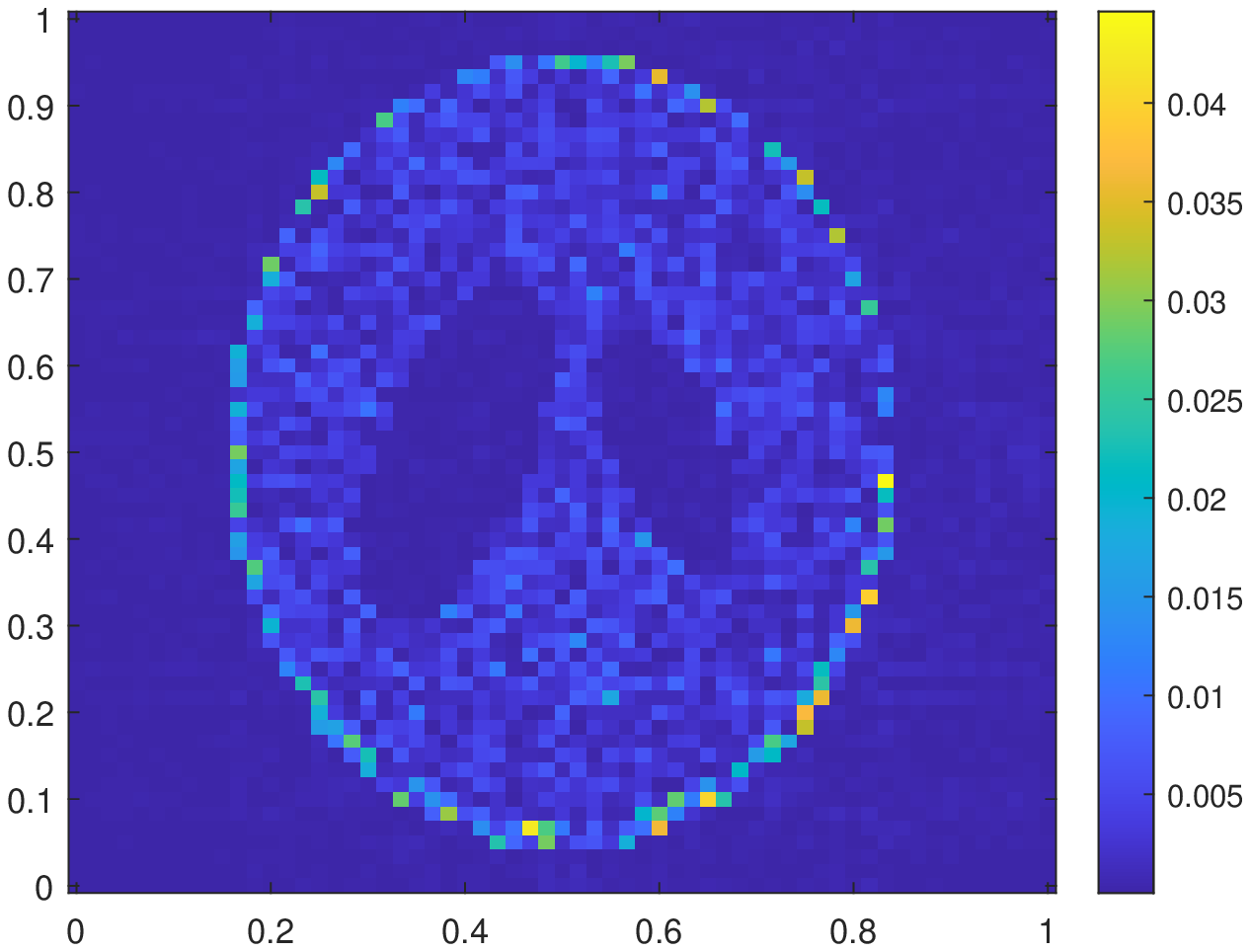}
 \includegraphics[width=0.23\textwidth]{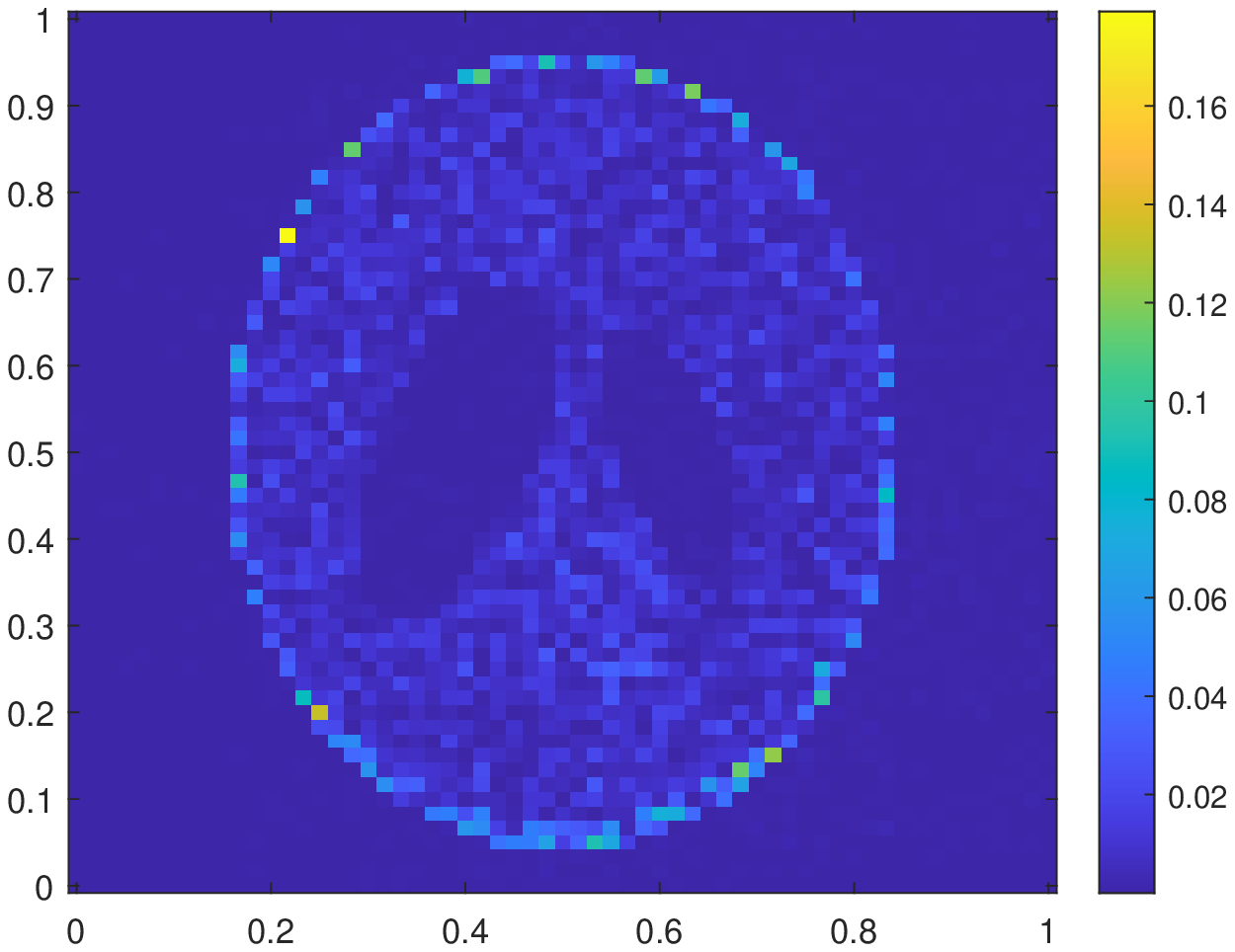}

 \caption{Reconstructed $S_2$ beyond the assumption of Theorem~\ref{thm:Neumann}. For the first row, 0\%, 1\%, 2\%, 5\% random noises are added to $H_{v_0}$. The relative $L^2$ errors of the reconstructions are 0.1408\%, 0.9864\%, 1.8437\%, 4.8237\%, respectively. 
 The second row displays the corresponding differences between the ground truth and the reconstructions.}
 \label{fig:Wrong2}
  \end{figure}










\bigskip
\textbf{Experiment 2: Inversion beyond the Assumption of Theorem~\ref{thm:Neumann}.}
The Neumann series in Theorem~\ref{thm:Neumann} was proved convergent under the sufficient condition~\eqref{eq:combound}. Here we also test the case when this condition fails. The experiment shows the series still converges in certain circumstances when the condition is violated.

We choose the computational domain $X=[0,1]\times [0,1]$, the attenuation coefficient $\sigma_2(x_1,x_2) = 1.1+0.2 x_1$. The constant $\rho=1$, the anisotropy parameter $g=0.5$, and the adjoint solution $v_0$ with $v_0|_{\Gamma_+}=1$ remain the same as in Experiment 1. Notice that $\textup{diam}(X) \rho = \sqrt{2} > 1 $ so the assumption~\eqref{eq:combound} does not hold. In this case, the well-posedness of the forward RTE is ensured by~\eqref{eqn:X1} but not by~\eqref{eqn:X2}. This is because
$$
\left(\inf_{x\in\overline{X}}\sigma\right) - \rho = 1.1 - 1 = 0.1 > 0.
$$
However, the first bound we obtained in~\eqref{eq:twobounds} is never less than 1, as was explained before Theorem~\ref{thm:Neumann}. This numerical experiment is therefore not covered by the proposed theorem.

We test the Neumann series inversion with a smooth source 
$$
S_3(x_1,x_2) = e^{-10[(x_1-0.4)^2+(x_2-0.6)^2]}
$$
and the discontinuous Shepp-Logan phantom $S_2$ (see Figure~\ref{fig:Wcoef1}). The reconstructions with different levels of noises are illustrated in Figure~\ref{fig:Wrong1} and Figure~\ref{fig:Wrong2}, respectively.

\bigskip
\textbf{Experiment 3: Inversion with Theorem~\ref{thm:Fredholm}.}
We test the Fredholm inversion Theorem~\ref{thm:Fredholm} in this experiment. Choose the computational domain $X=[0,1]\times [0,1]$, the attenuation coefficient $\sigma_1(x_1,x_2) = 0.1+0.1x_1$ (see Figure~\ref{fig:Ncoef1}), $g=0.5$ in the scattering kernel~\eqref{eq:HGfunc}, the adjoint RTE solution $v_0$ with $v_0|_{\Gamma_+}=1$.

\begin{figure}[!htb]
\includegraphics[width=0.2\textwidth]{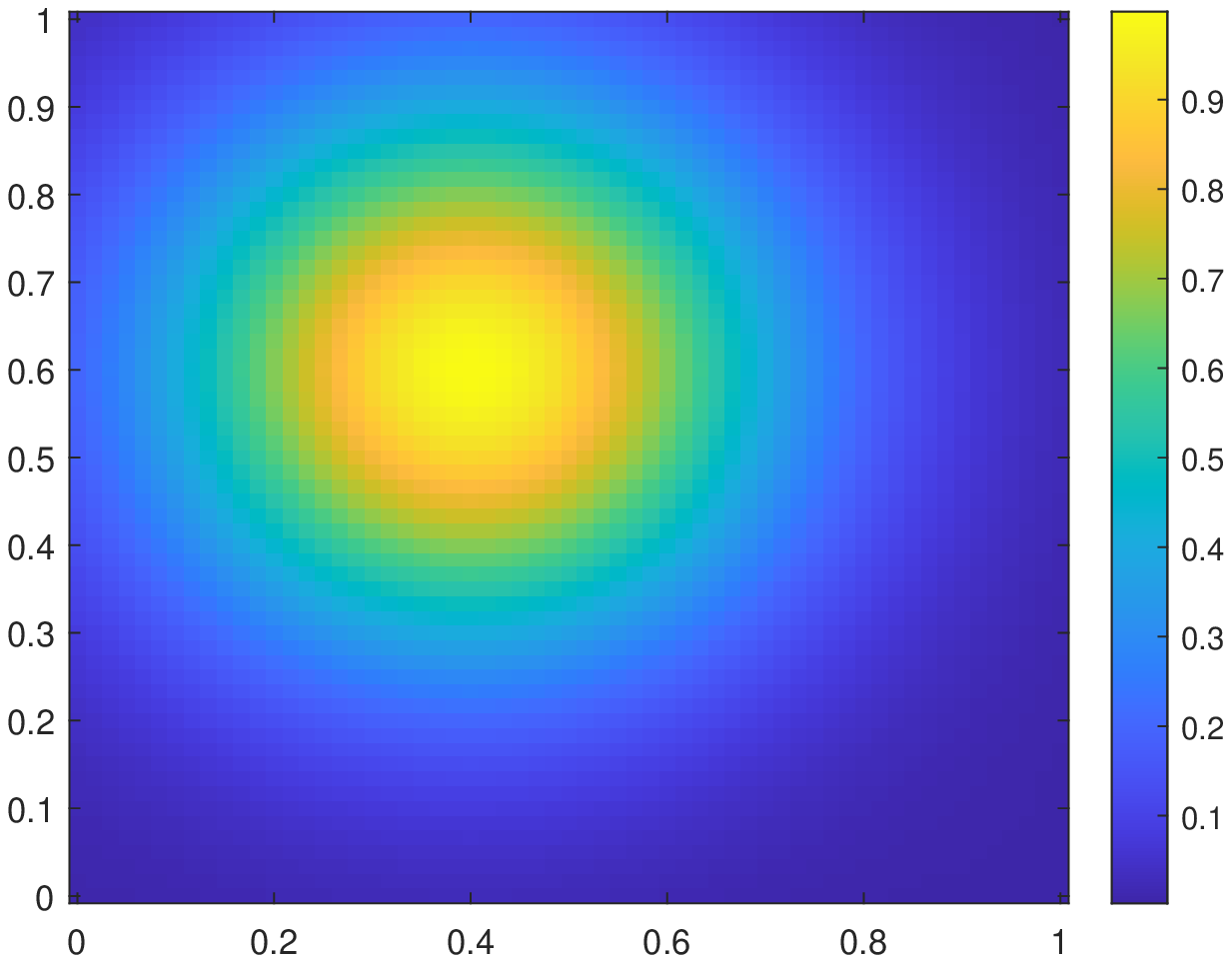}
\includegraphics[width=0.2\textwidth]{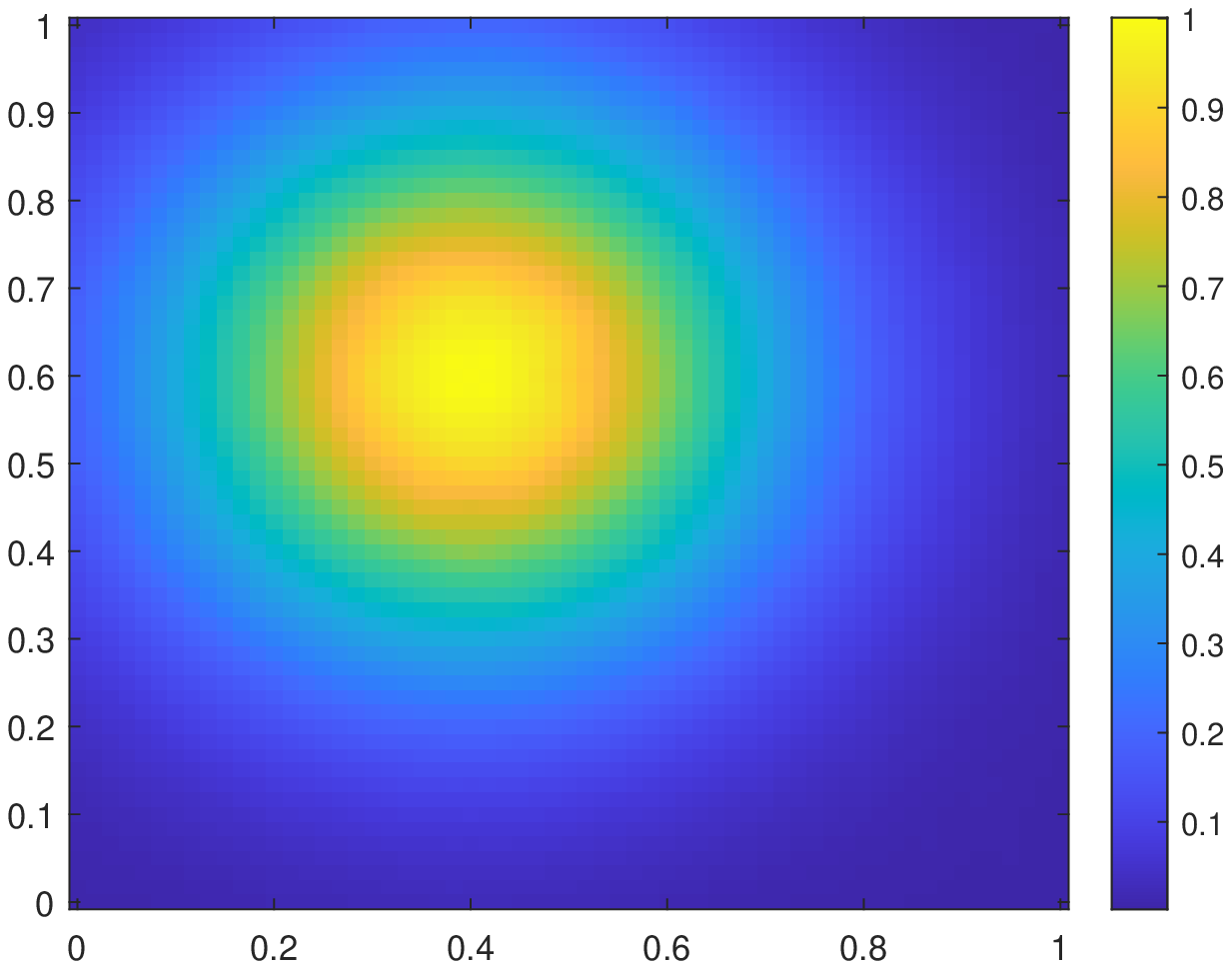}
\includegraphics[width=0.2\textwidth]{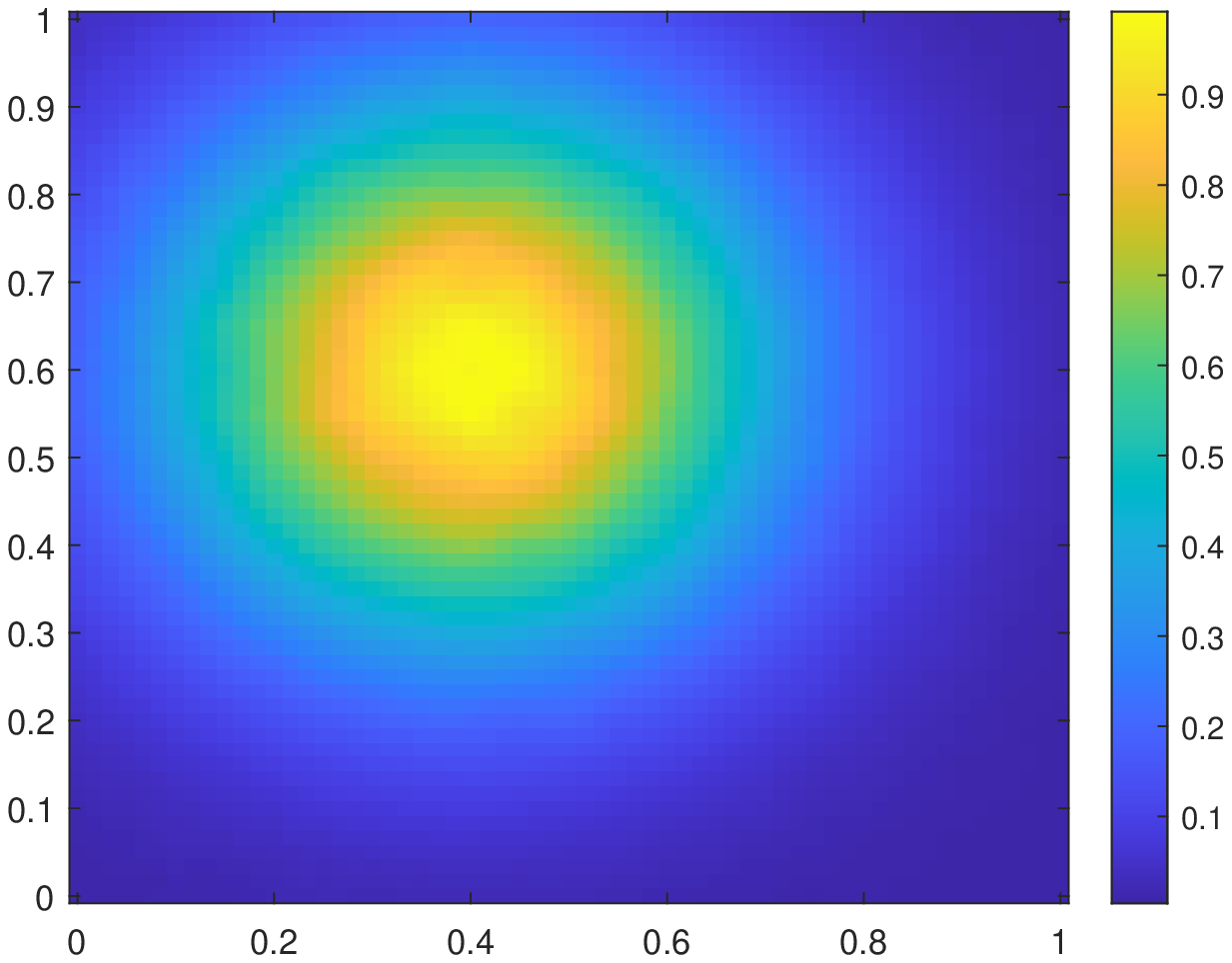}
\includegraphics[width=0.2\textwidth]{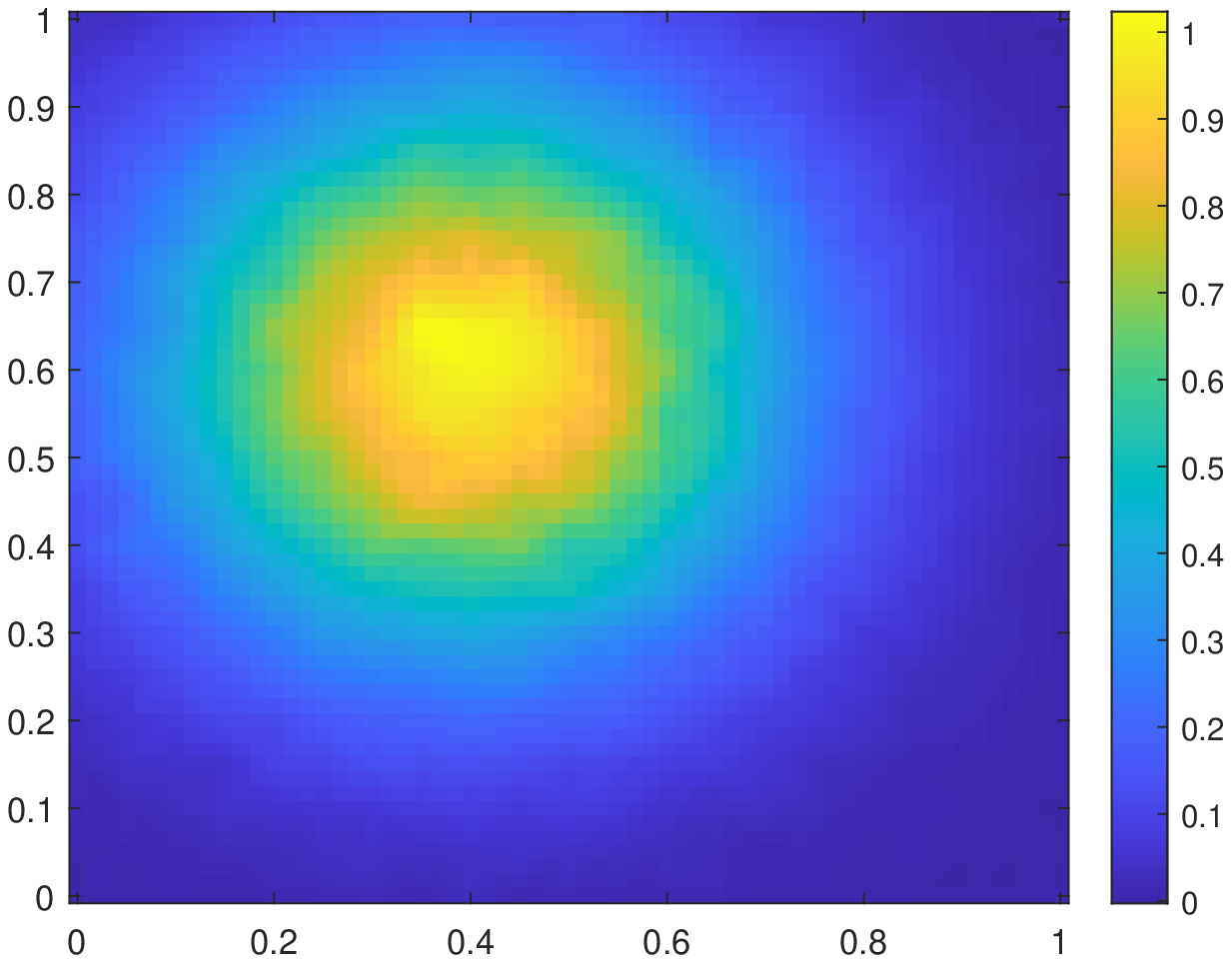}

\includegraphics[width=0.2\textwidth]{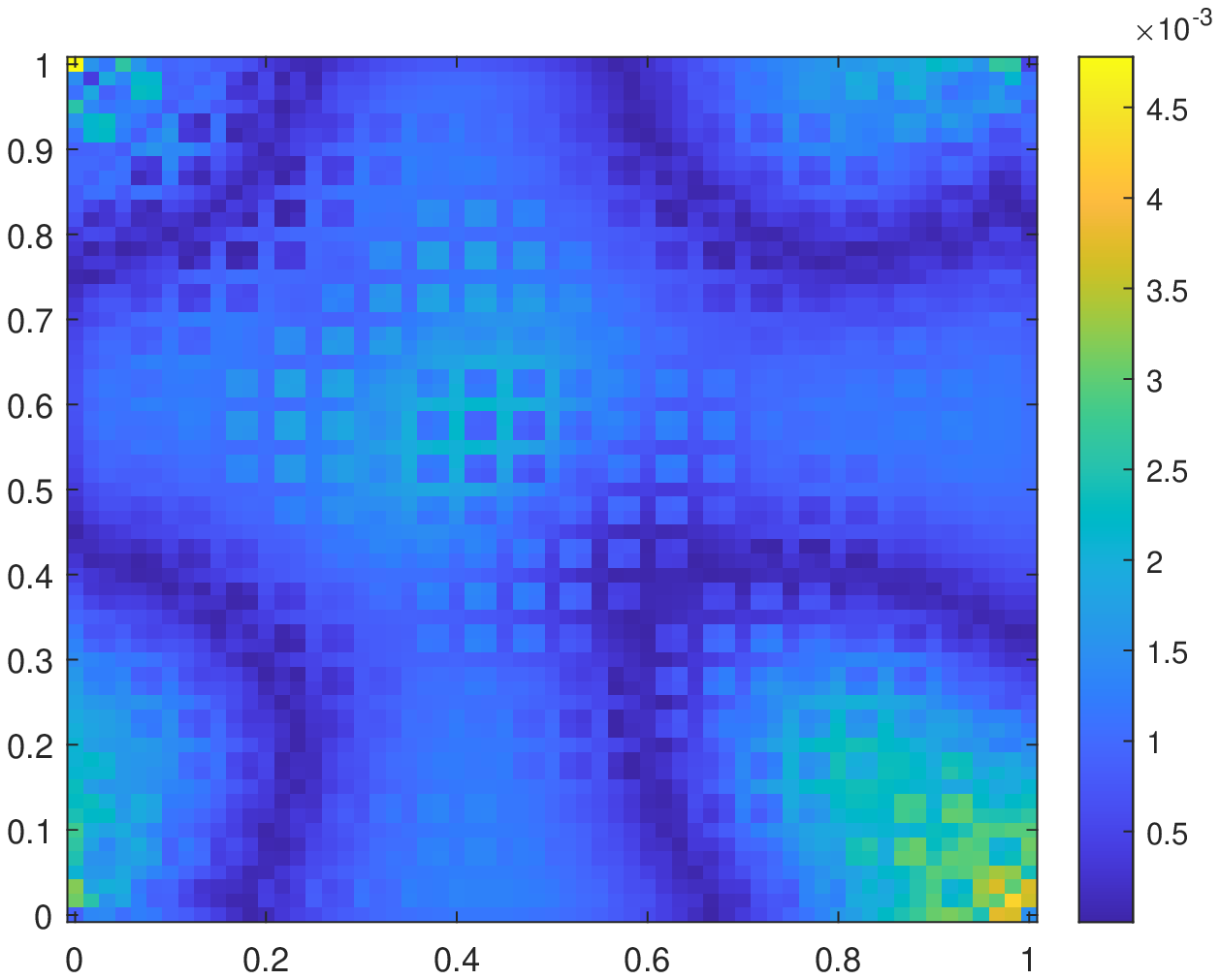}
\includegraphics[width=0.2\textwidth]{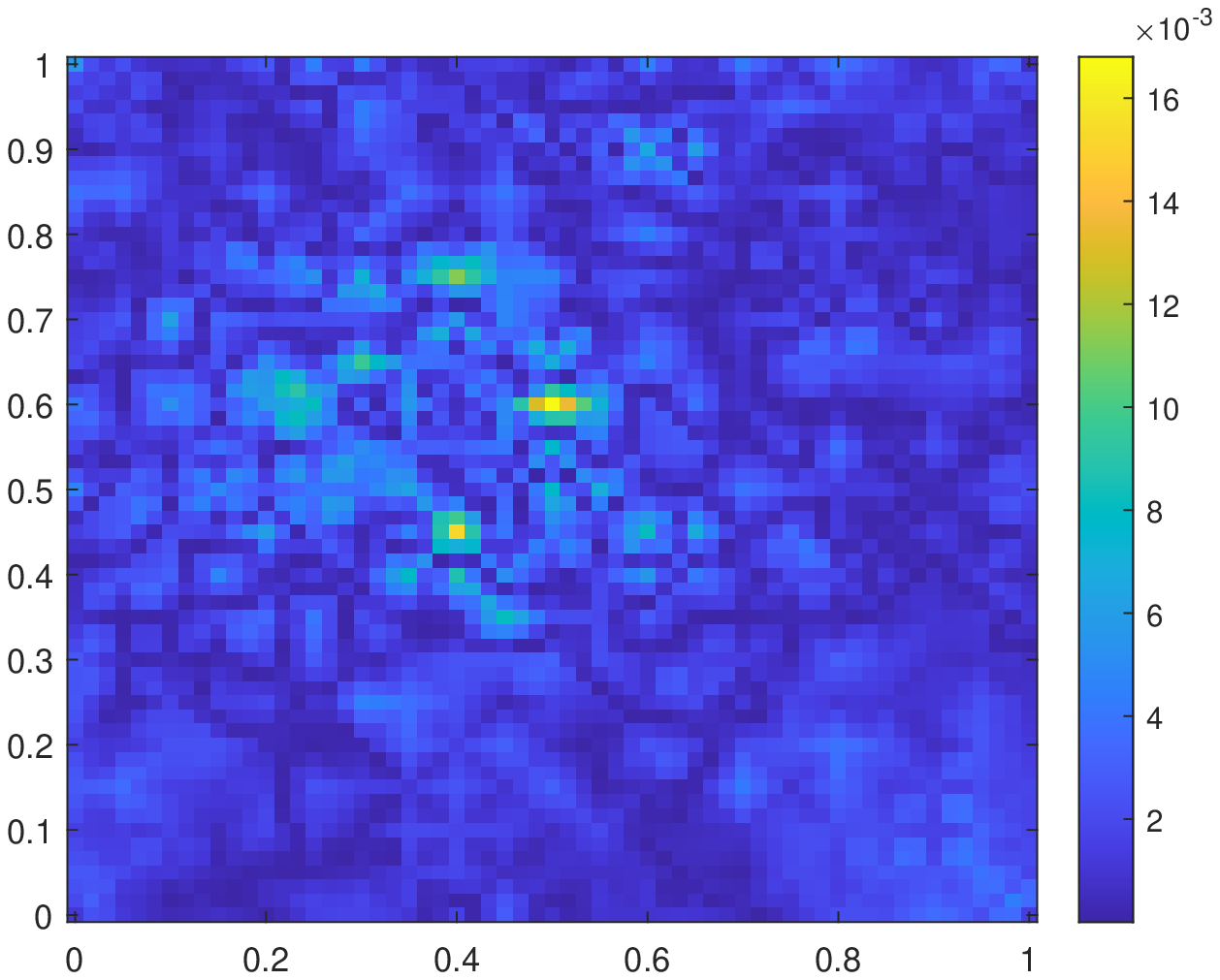}
\includegraphics[width=0.2\textwidth]{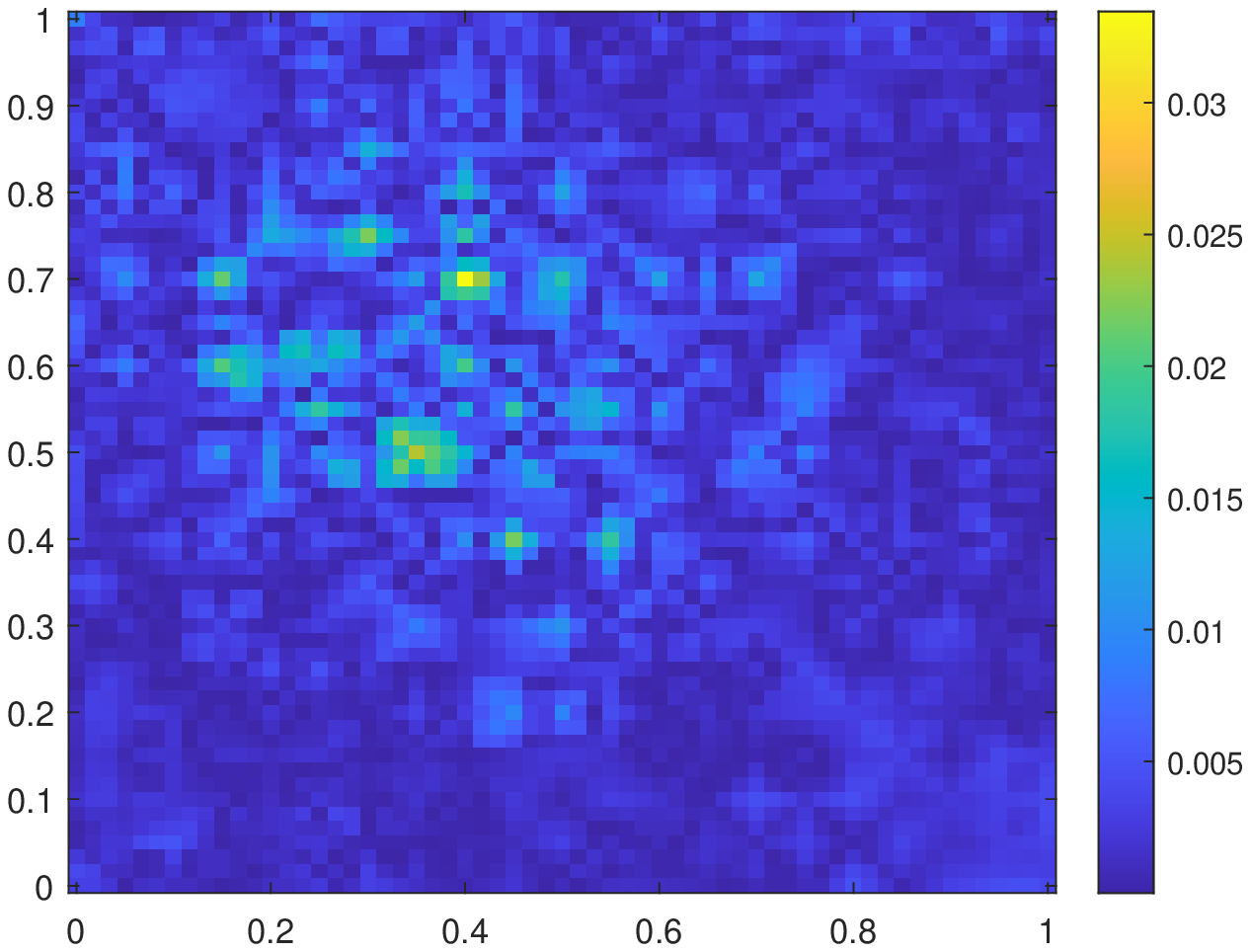}
\includegraphics[width=0.2\textwidth]{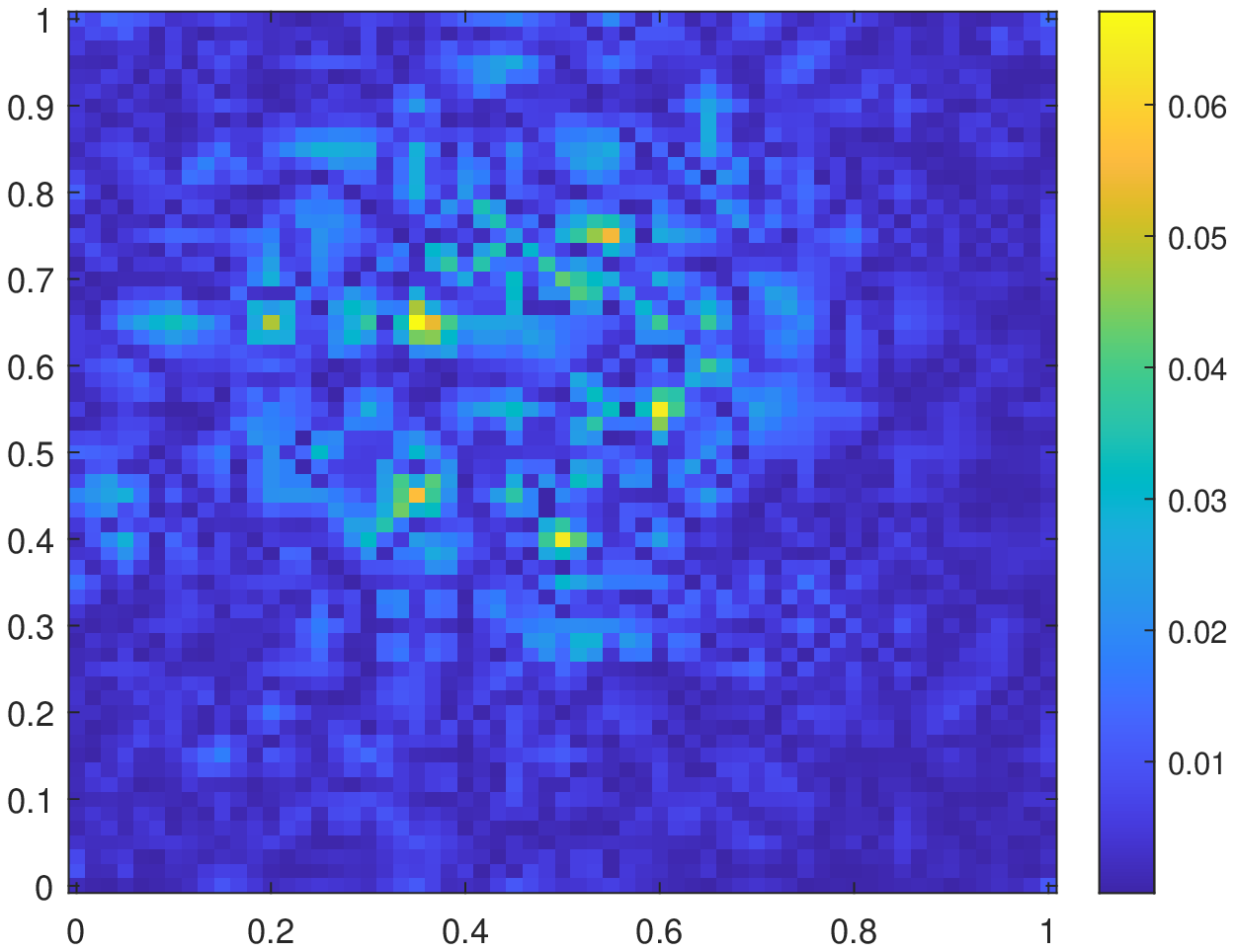}
\caption{Reconstructed $S_3$ using Fredholm inversion. For the first row, 0\%, 1\%, 2\%, 5\% random noises are added to $H_{v_0}$. The relative $L^2$ errors of the reconstructions are 0.2878\%, 0.5784\%, 1.0792\%, 2.6614\%, respectively. 
The second row displays the corresponding differences between the ground truth and the reconstructions.}
\label{fig:Fredholm1}
\vspace*{\floatsep}
\includegraphics[width=0.2\textwidth]{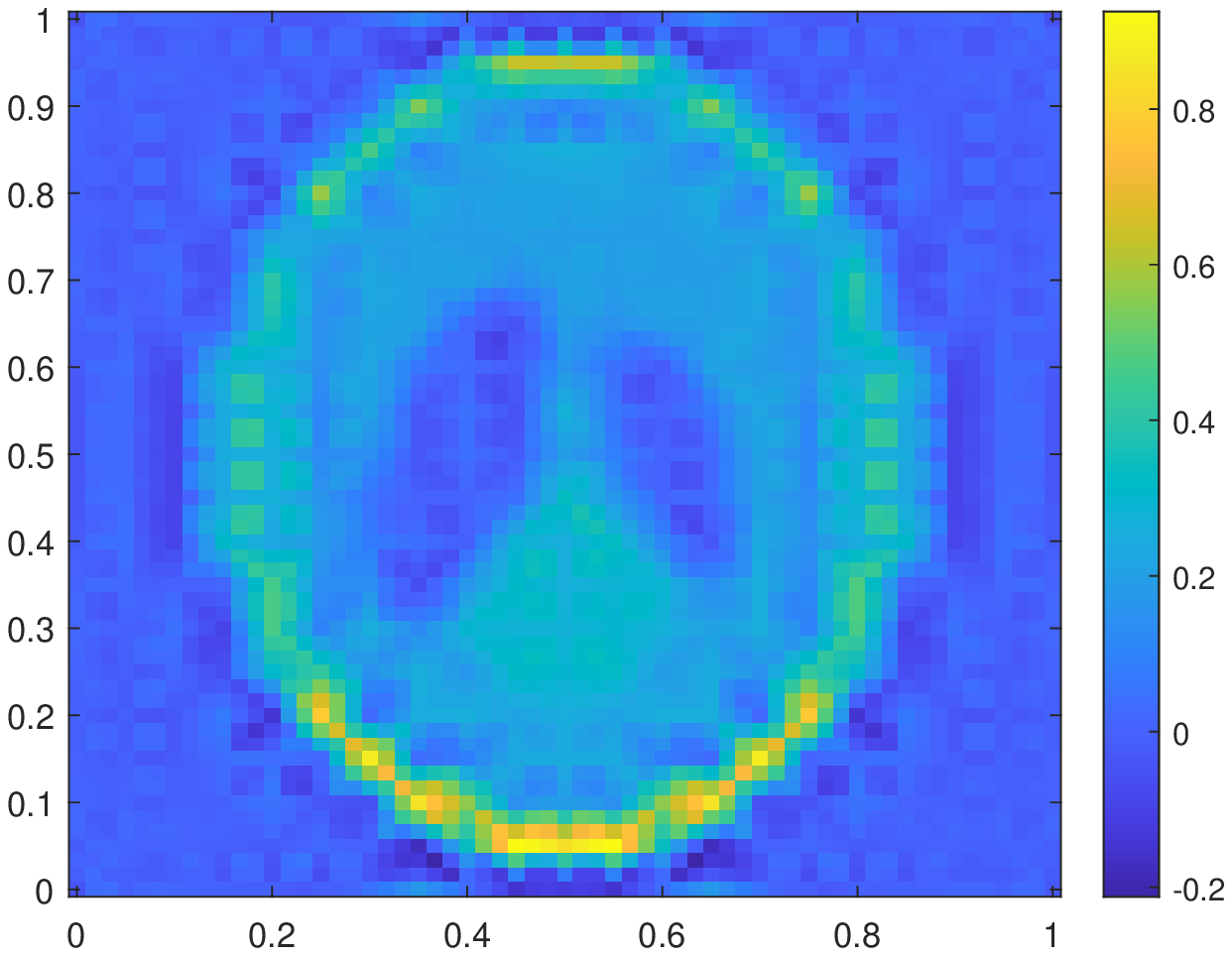}
\includegraphics[width=0.2\textwidth]{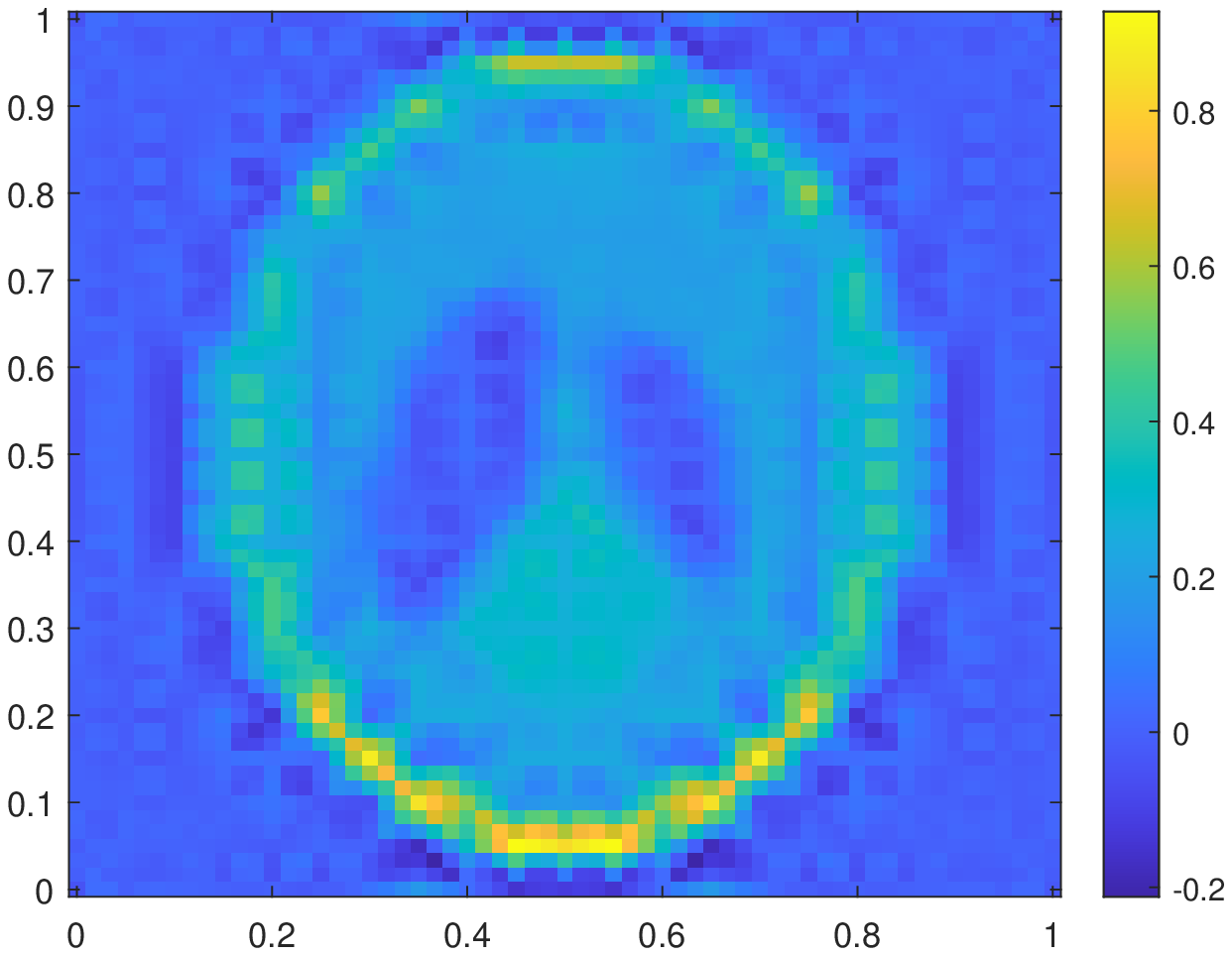}
\includegraphics[width=0.2\textwidth]{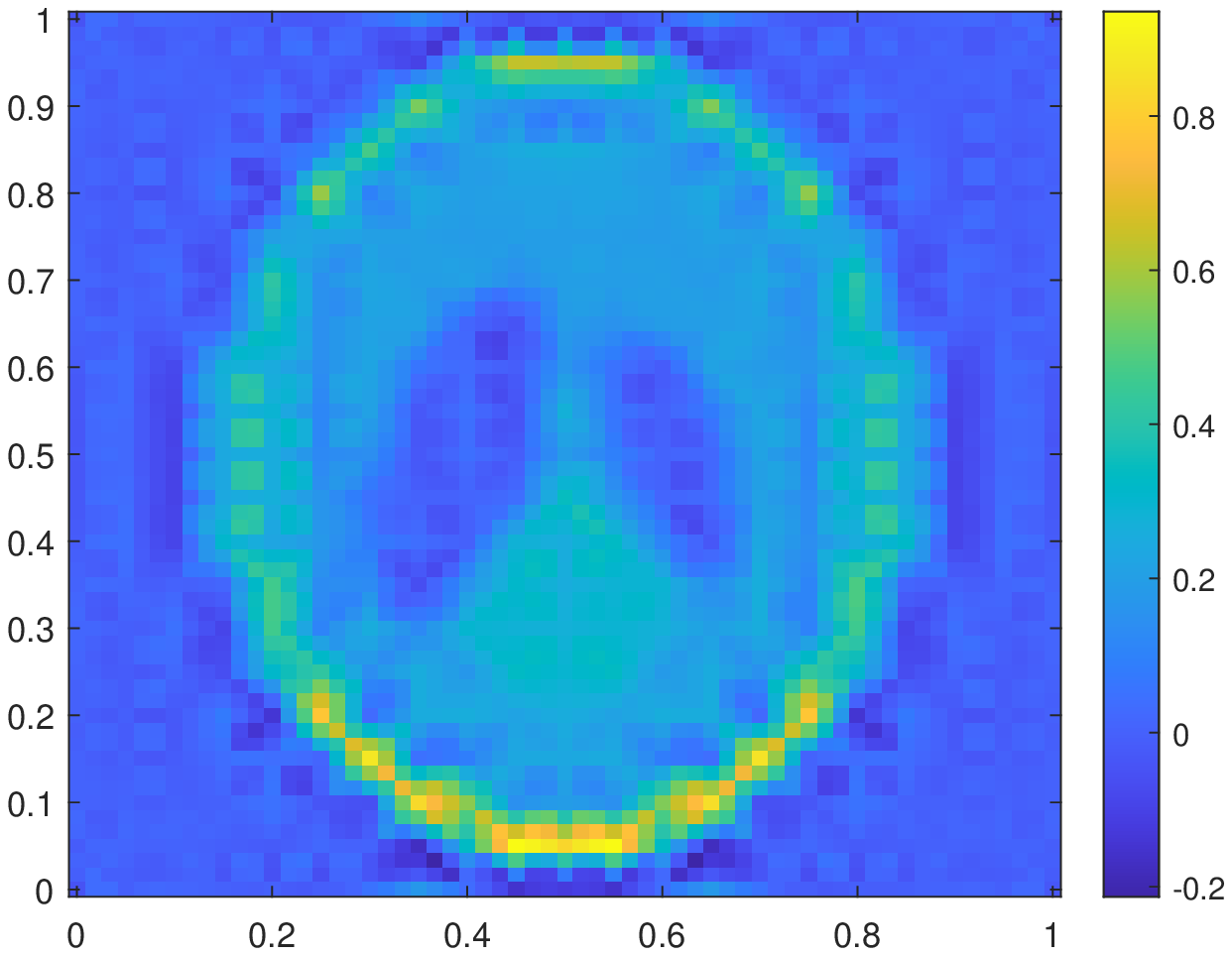}
\includegraphics[width=0.2\textwidth]{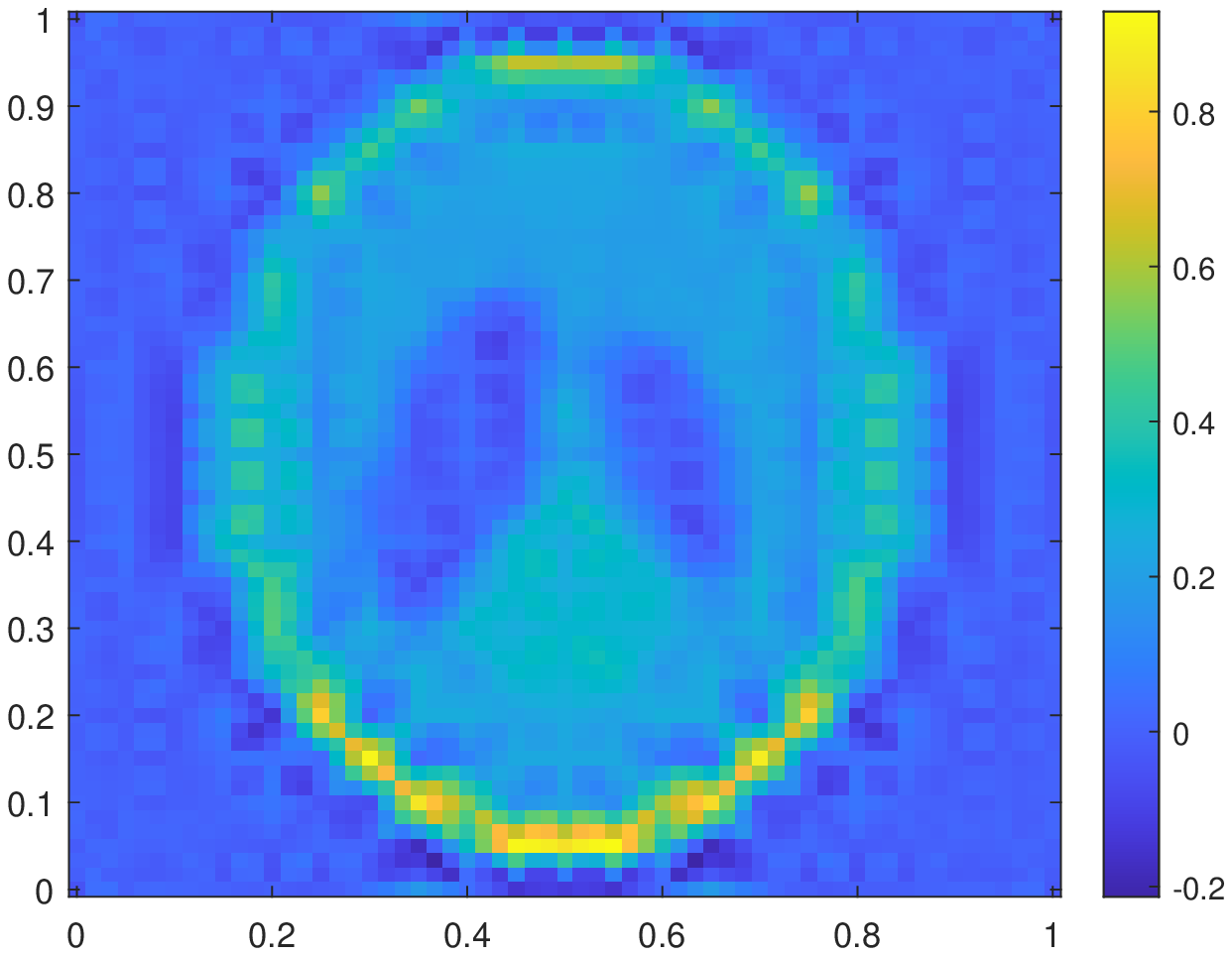}

\includegraphics[width=0.2\textwidth]{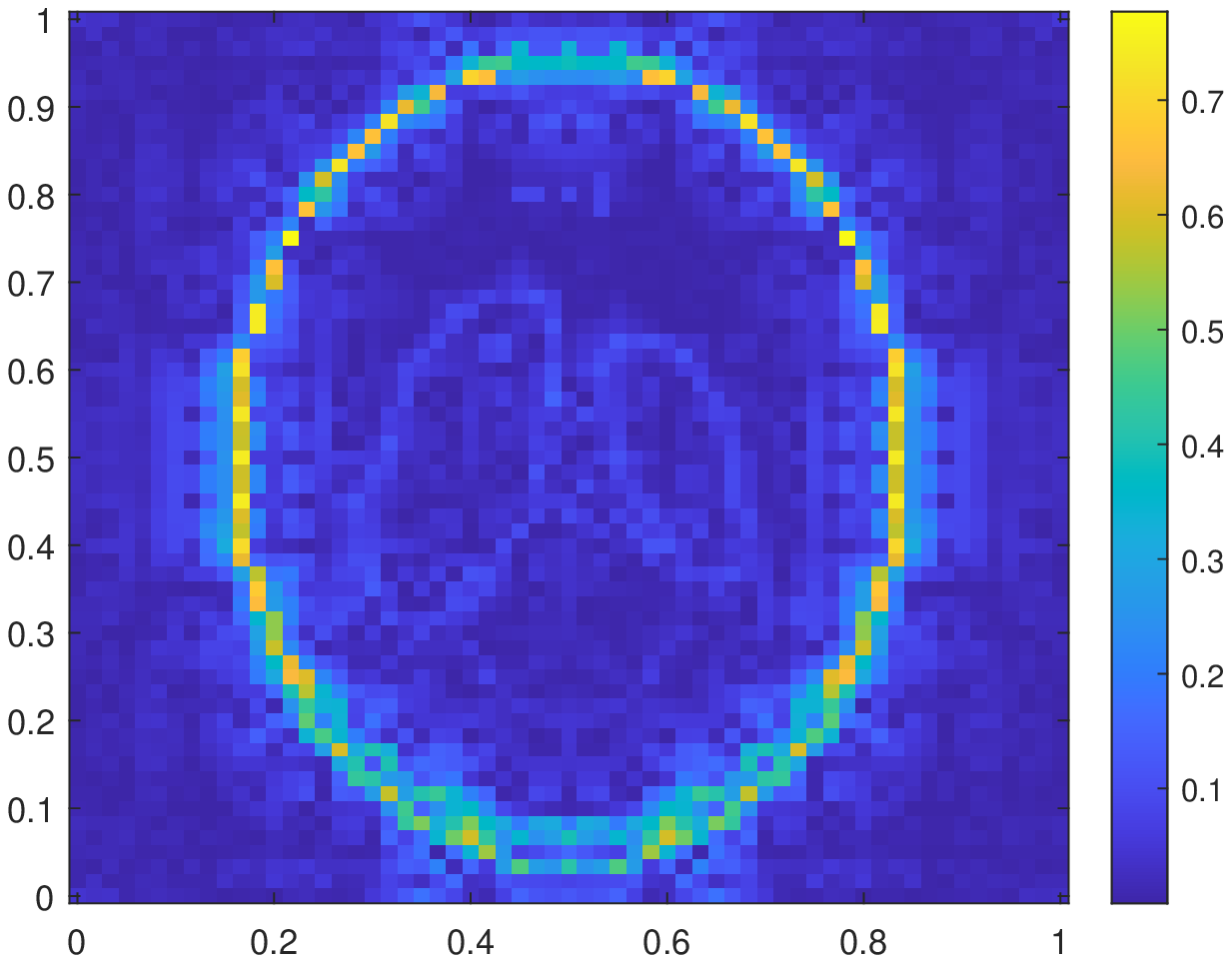}
\includegraphics[width=0.2\textwidth]{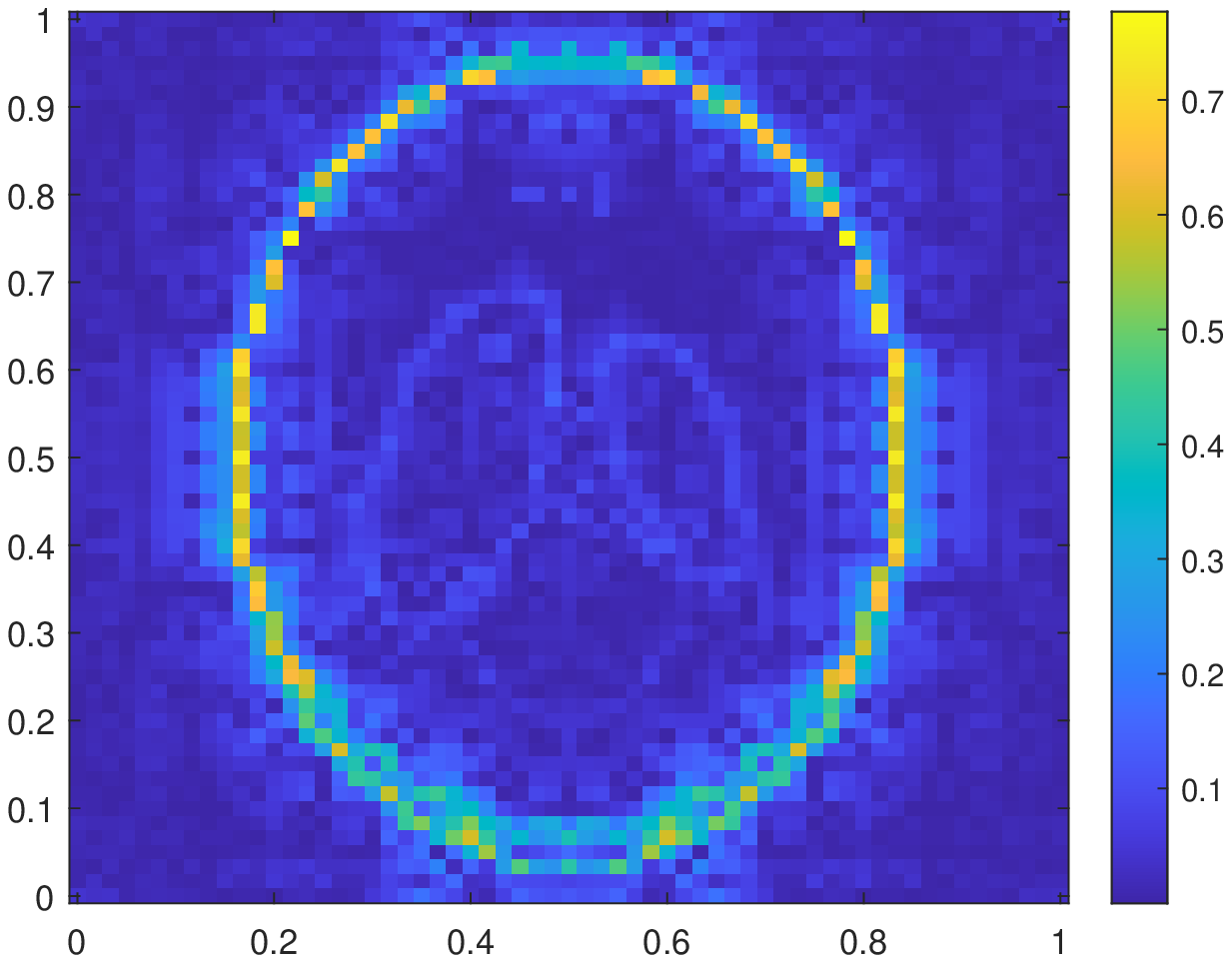}
\includegraphics[width=0.2\textwidth]{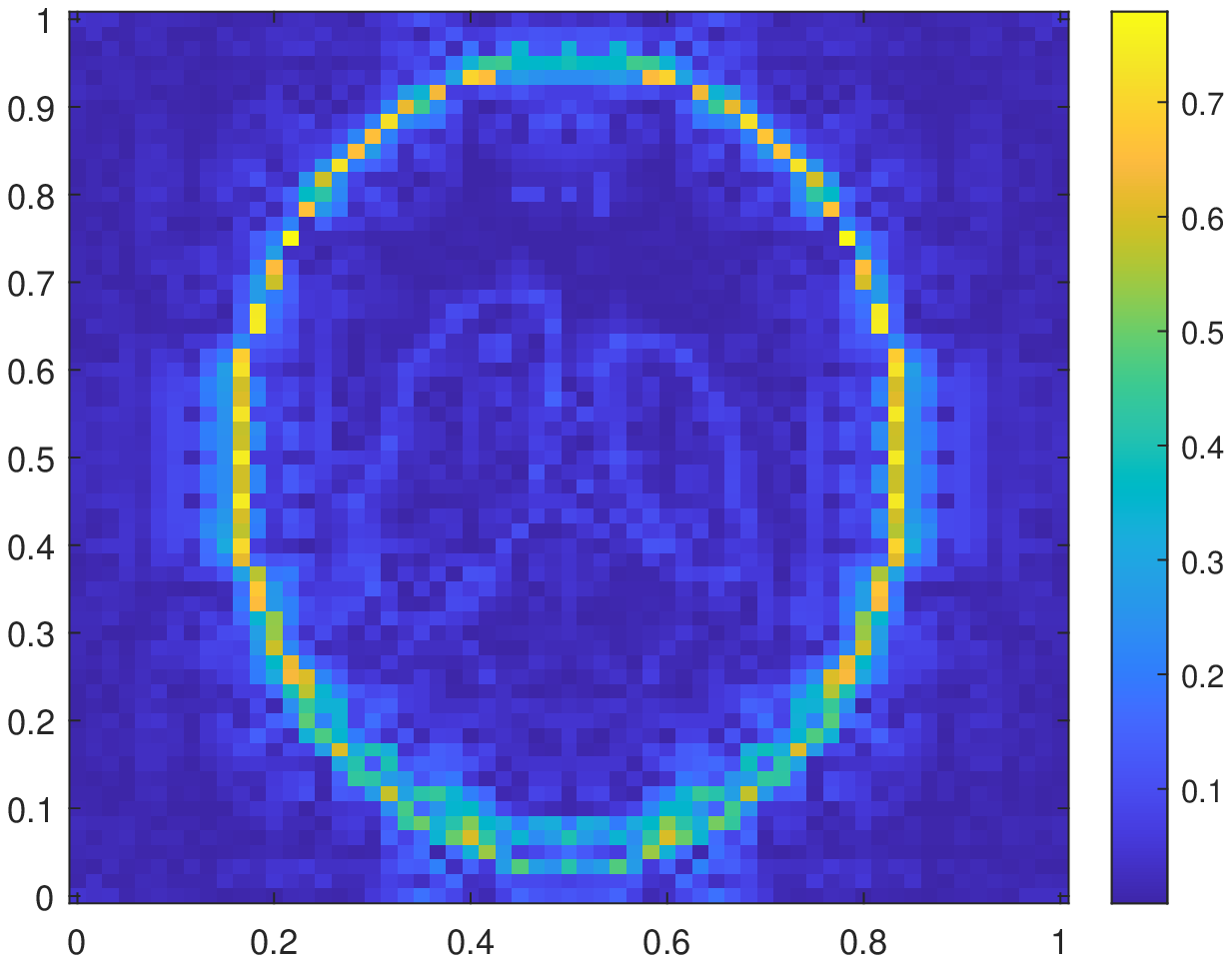}
\includegraphics[width=0.2\textwidth]{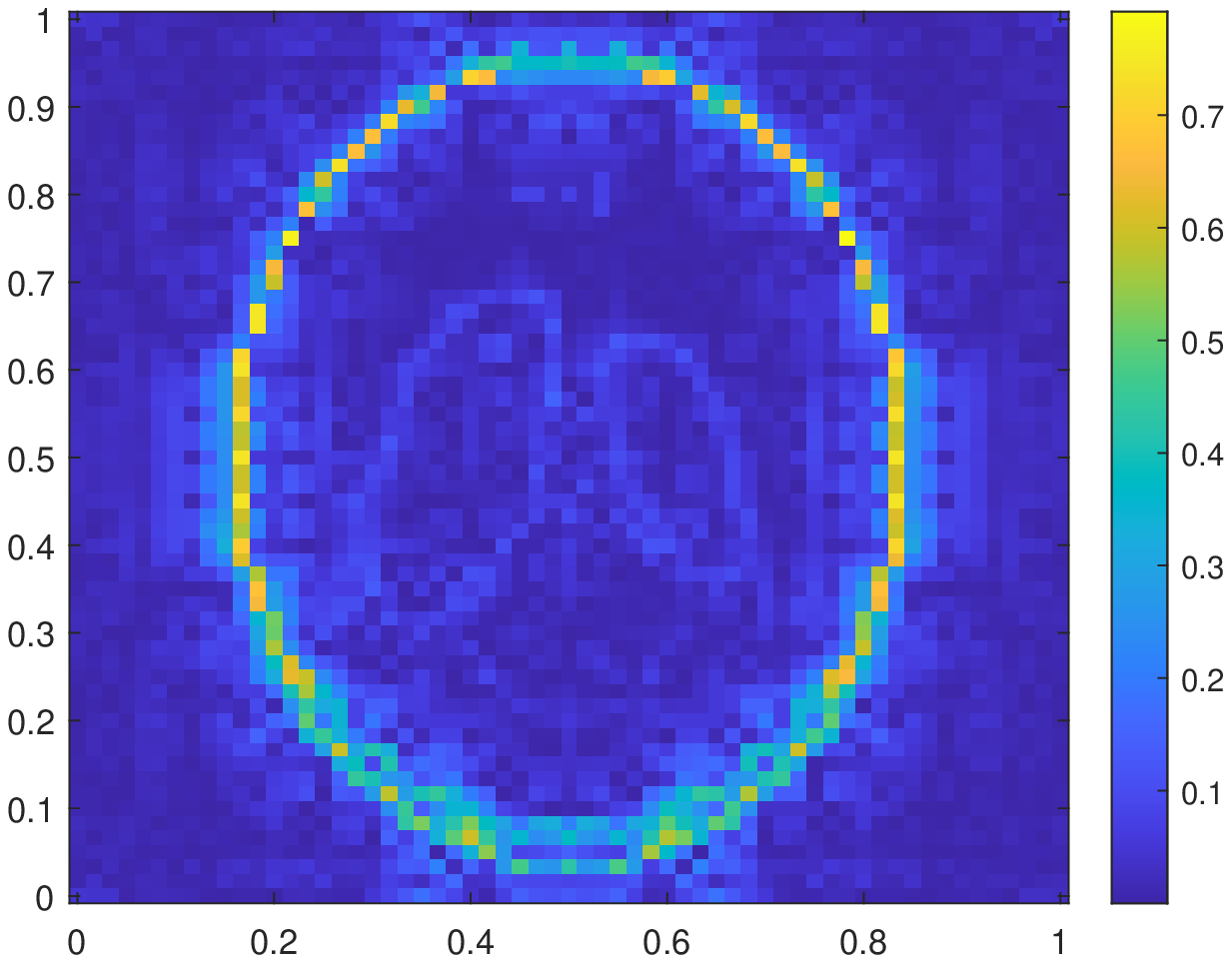}
\caption{Reconstructed $S_2$ using Fredholm inversion. For the first row, 0\%, 1\%, 2\%, 5\% random noises are added to $H_{v_0}$. The relative $L^2$ errors of the reconstructions are 57.5806\%, 57.5818\%, 57.5880\%, 57.6199\%, respectively. 
The second row displays the corresponding differences between the ground truth and the reconstructions.}
\label{fig:Fredholm2}
\end{figure}

We test the Fredholm inversion with the smooth source $S_3$ 
(see Figure~\ref{fig:Wcoef1}) and the discontinuous Shepp-Logan phantom $S_2$ (see Figure~\ref{fig:Ncoef1}). The reconstructions with different levels of noises are illustrated in Figure~\ref{fig:Fredholm1} and Figure~\ref{fig:Fredholm2}, respectively.

\begin{figure}[!htb]
\centering
\includegraphics[width=0.45\textwidth]{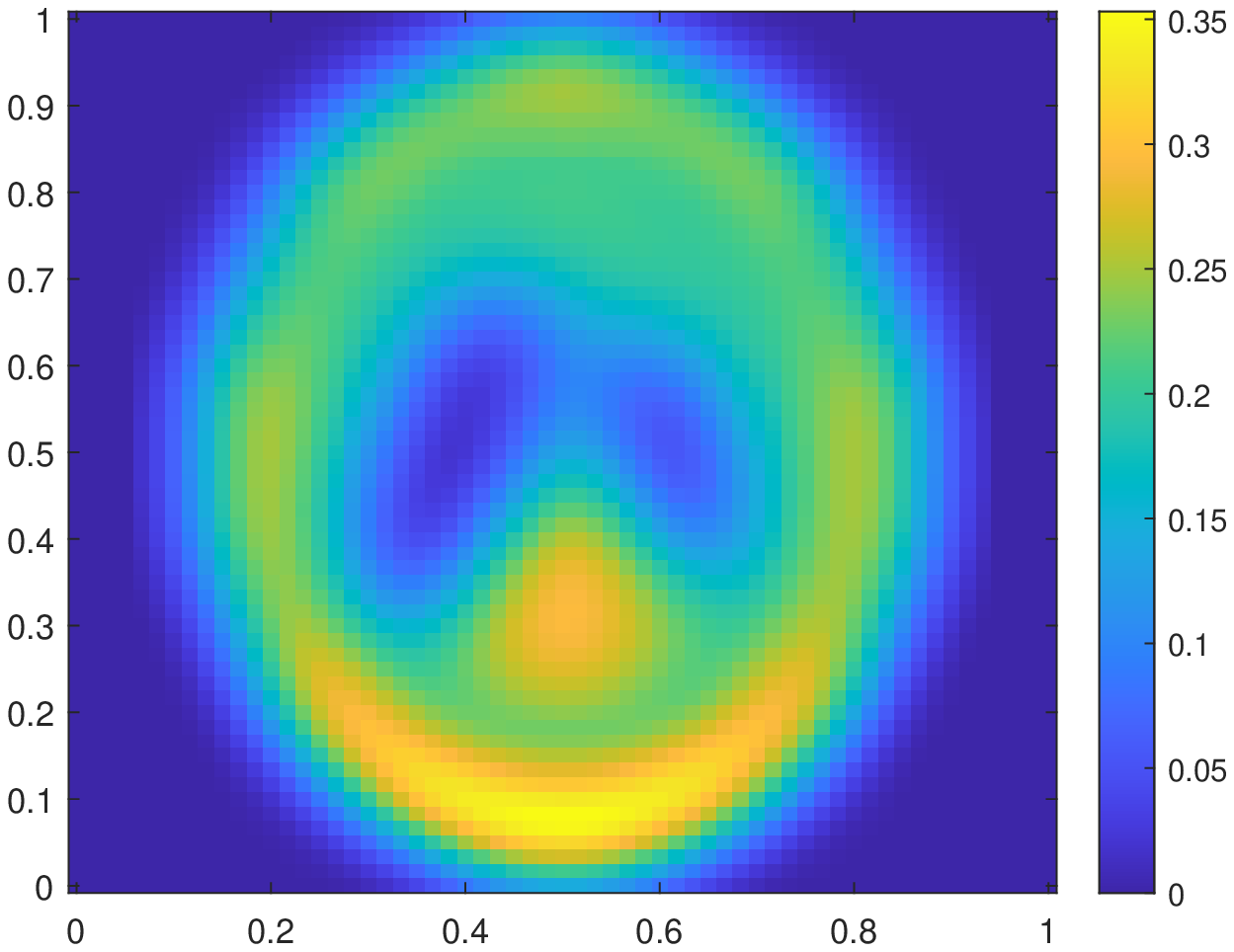}
\caption{Smoothed Shepp-Logan phantom $S_4$.}
\label{fig:Fcoef3}
\vspace*{\floatsep}
\includegraphics[width=0.2\textwidth]{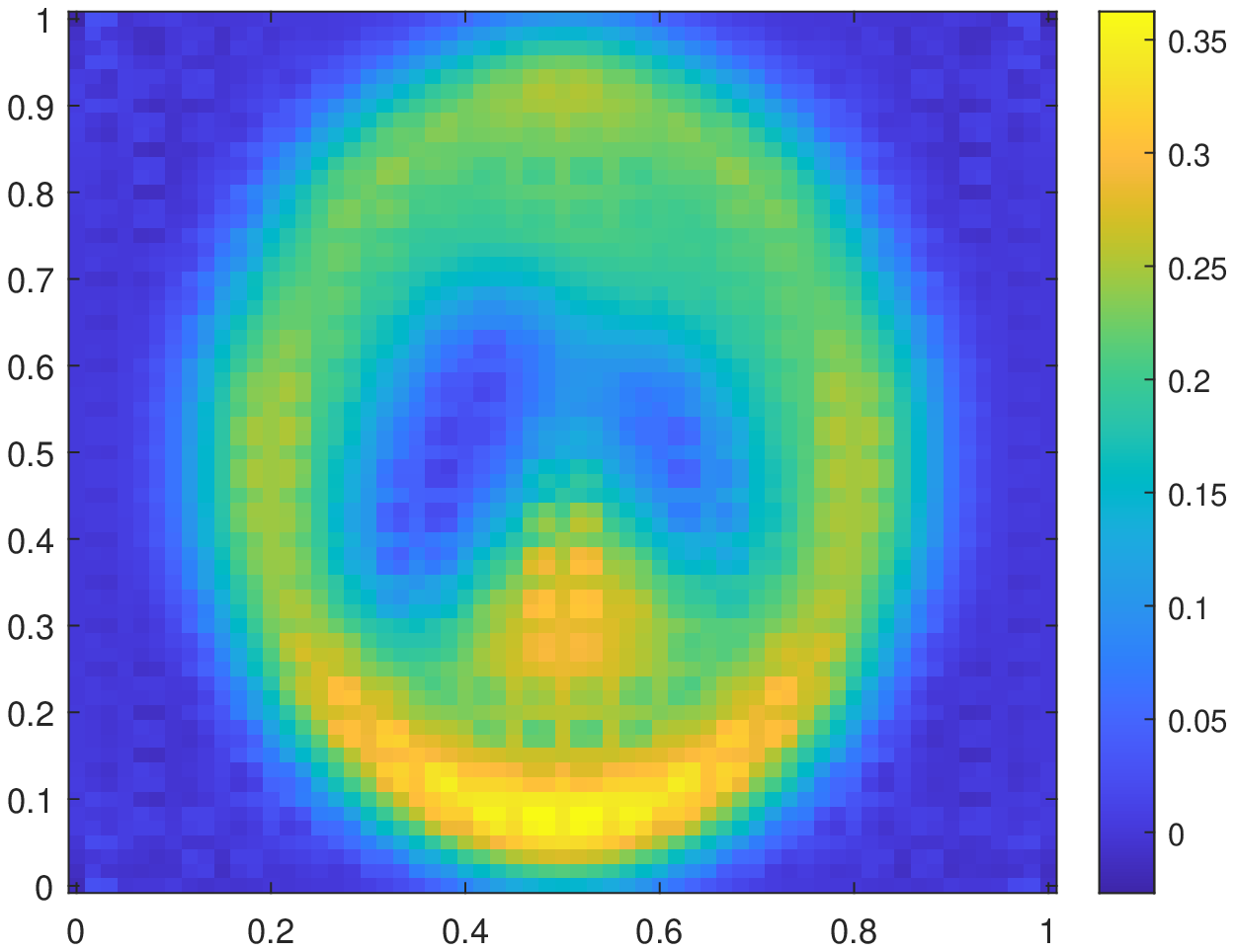}
\includegraphics[width=0.2\textwidth]{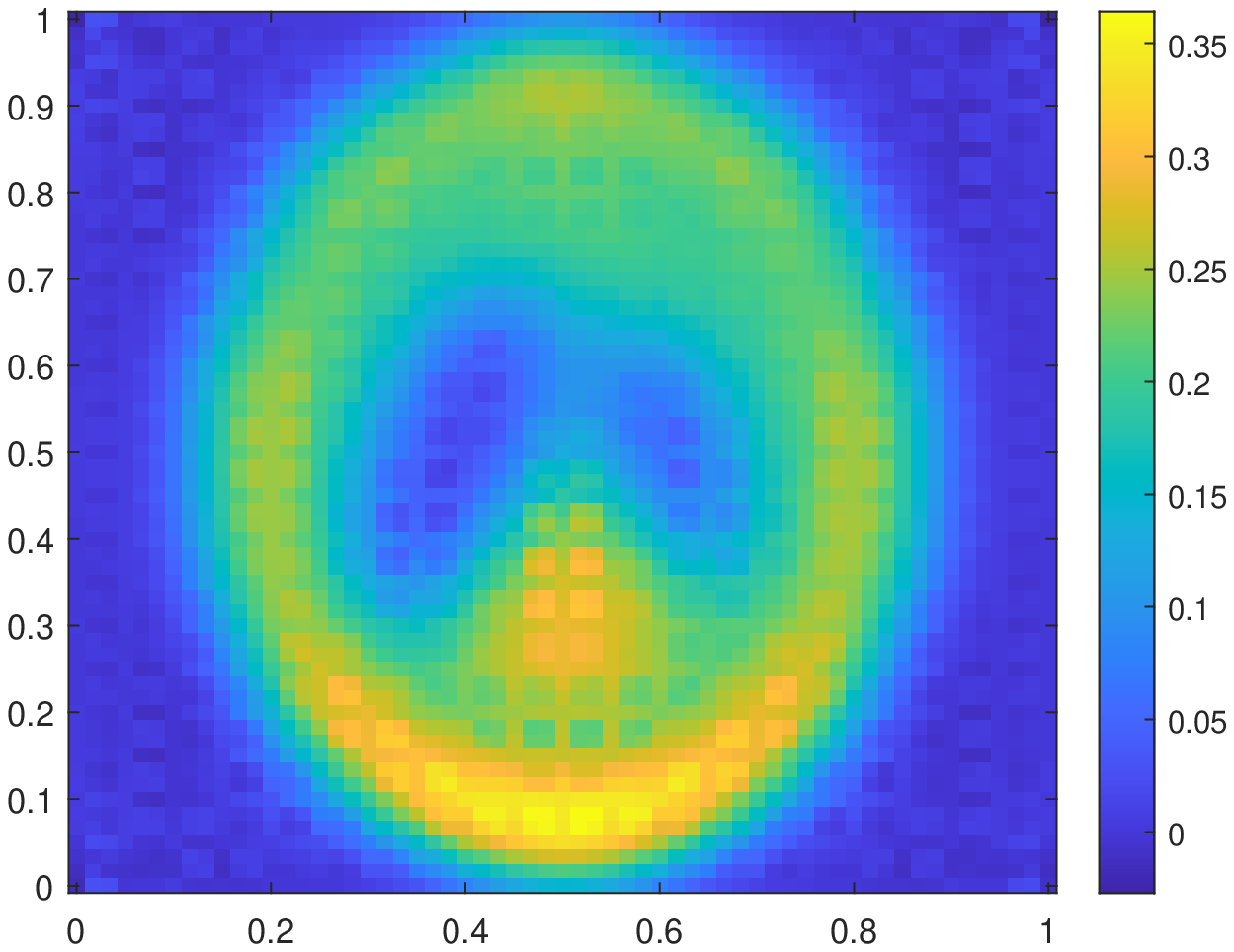}
\includegraphics[width=0.2\textwidth]{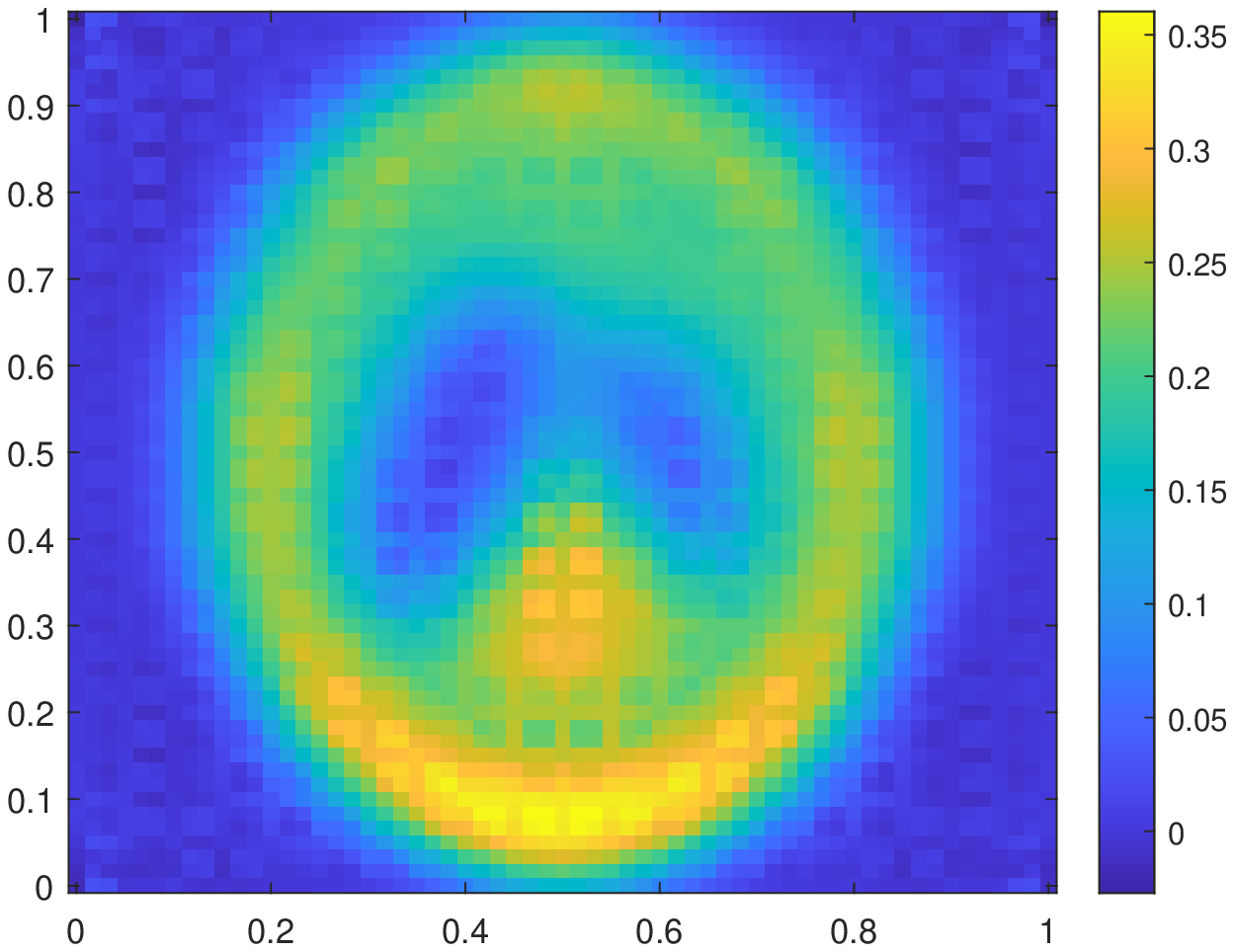}
\includegraphics[width=0.2\textwidth]{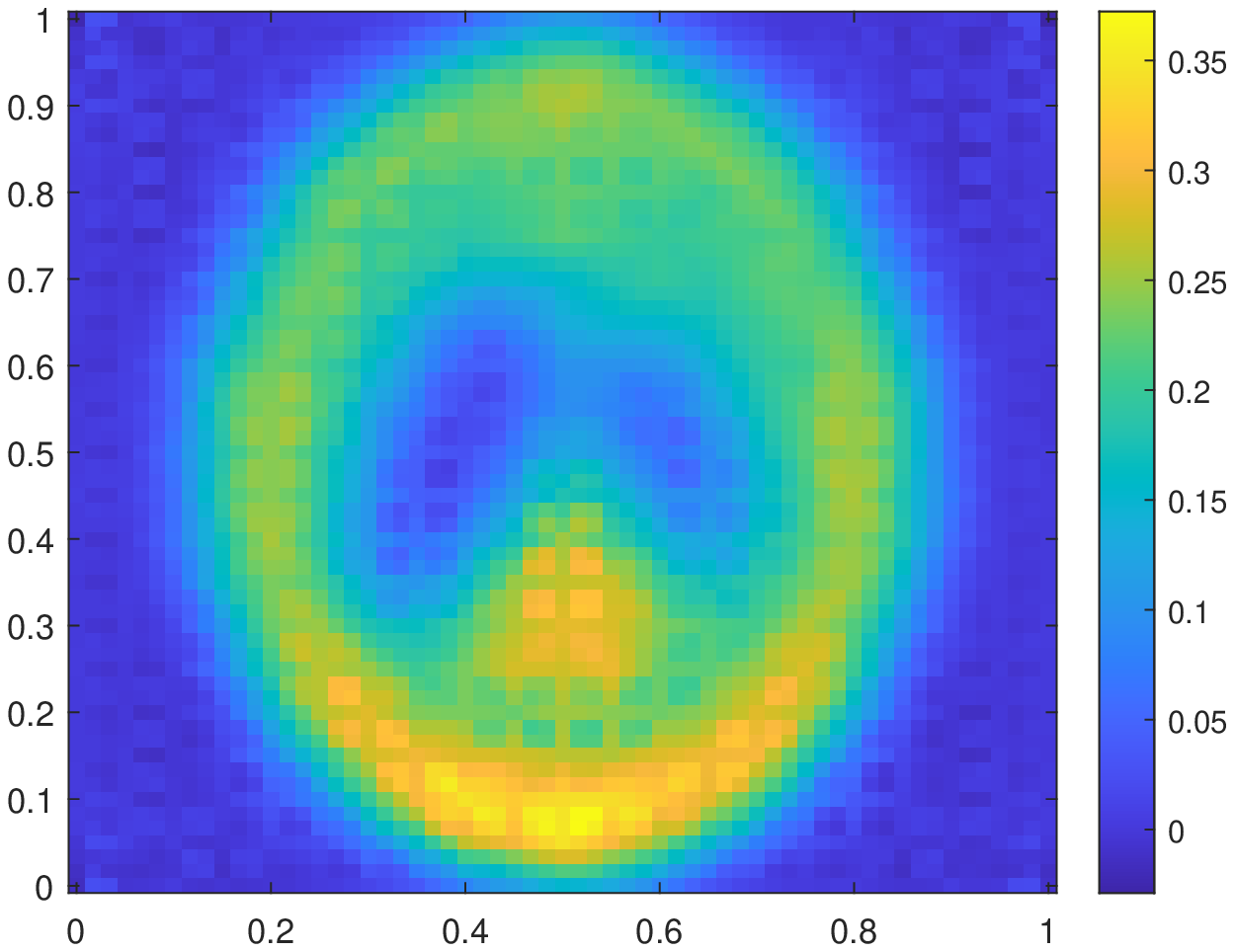}

\includegraphics[width=0.2\textwidth]{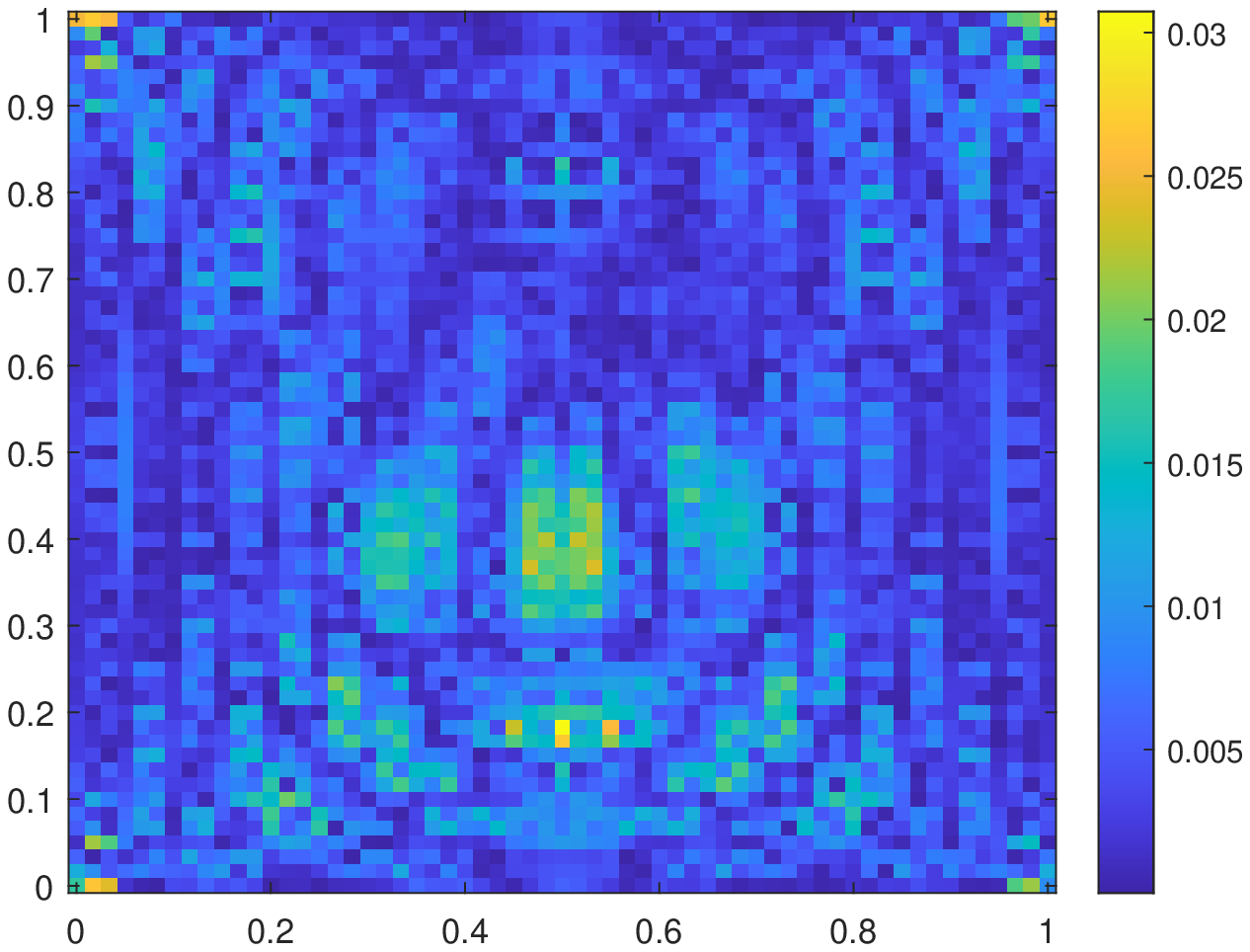}
\includegraphics[width=0.2\textwidth]{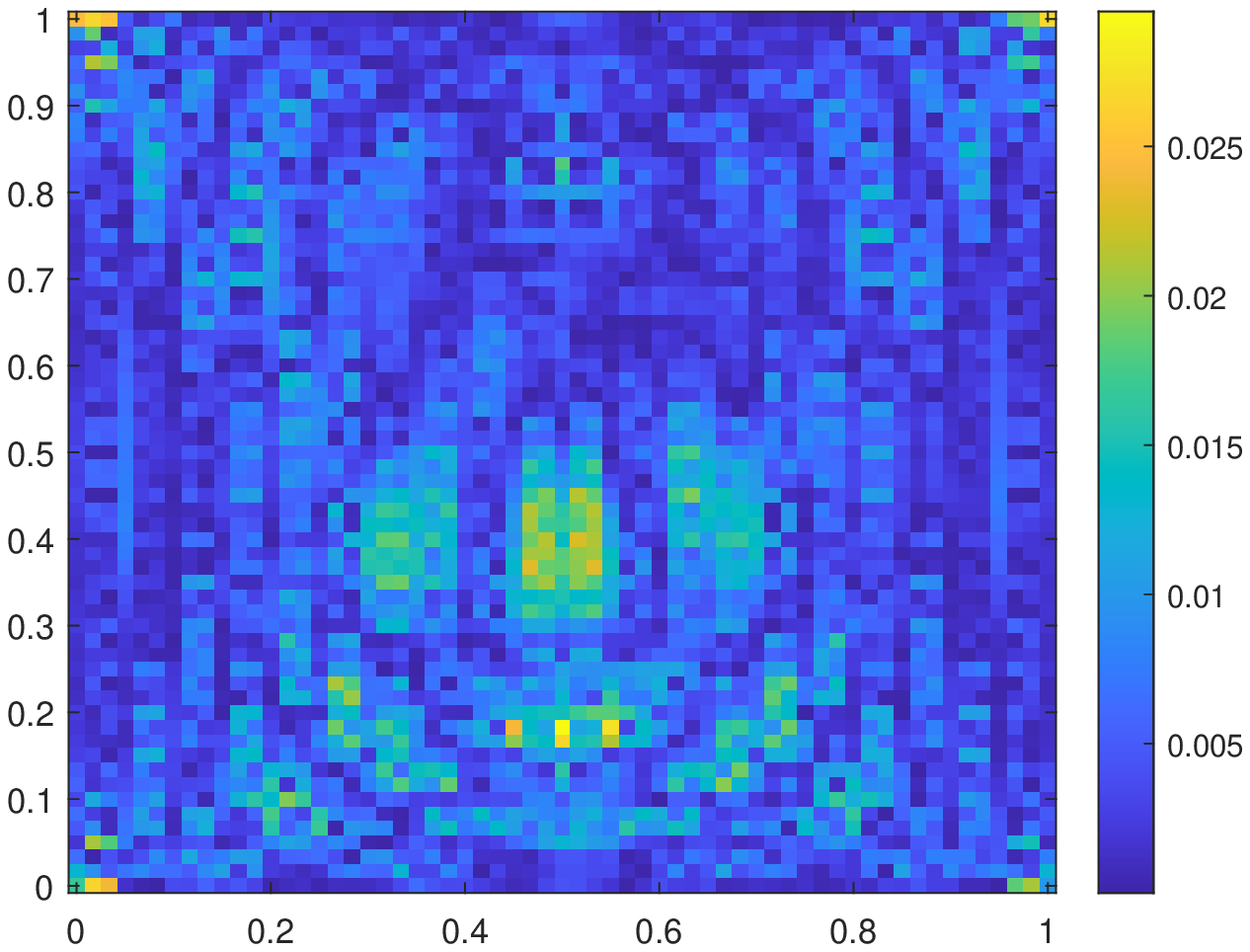}
\includegraphics[width=0.2\textwidth]{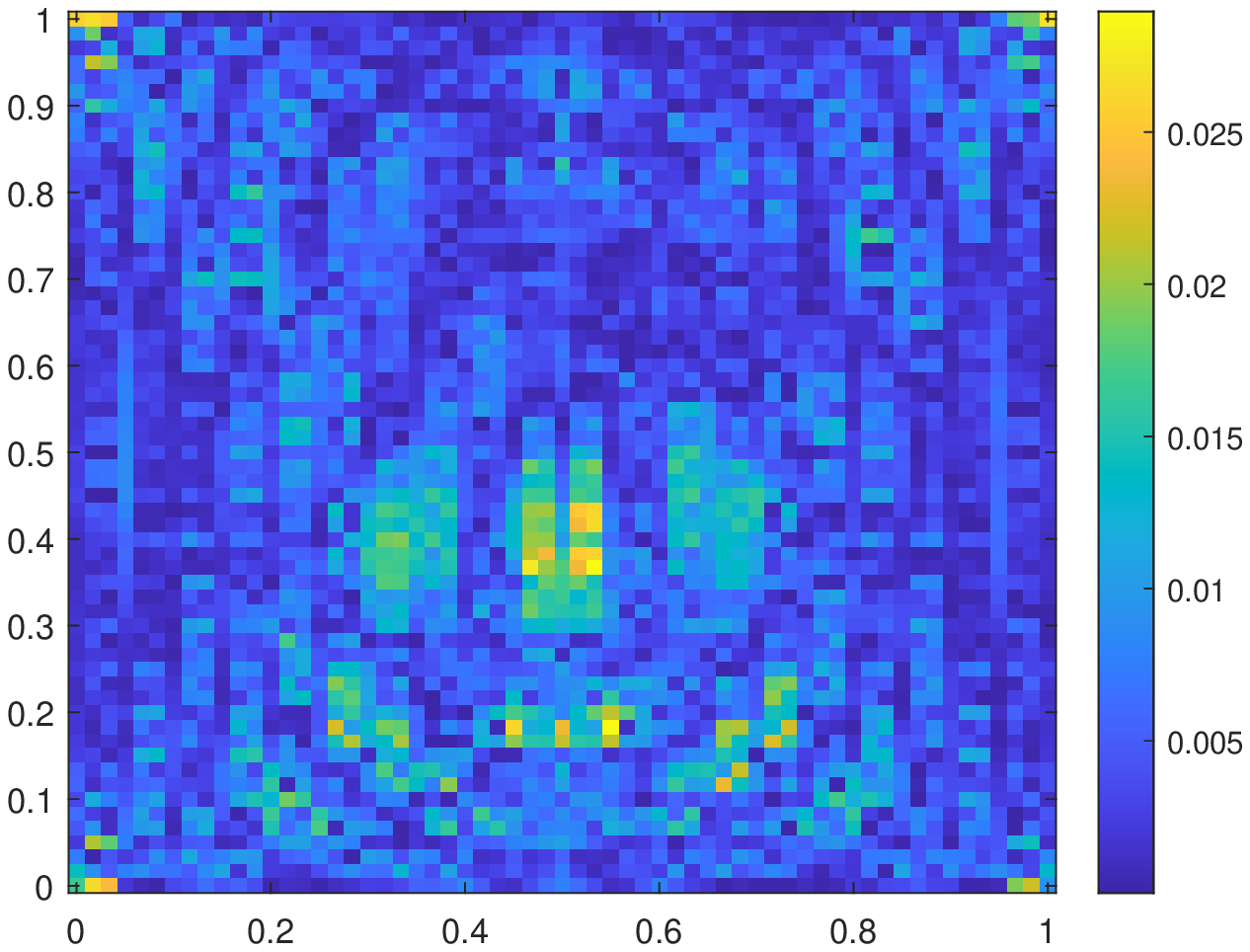}
\includegraphics[width=0.2\textwidth]{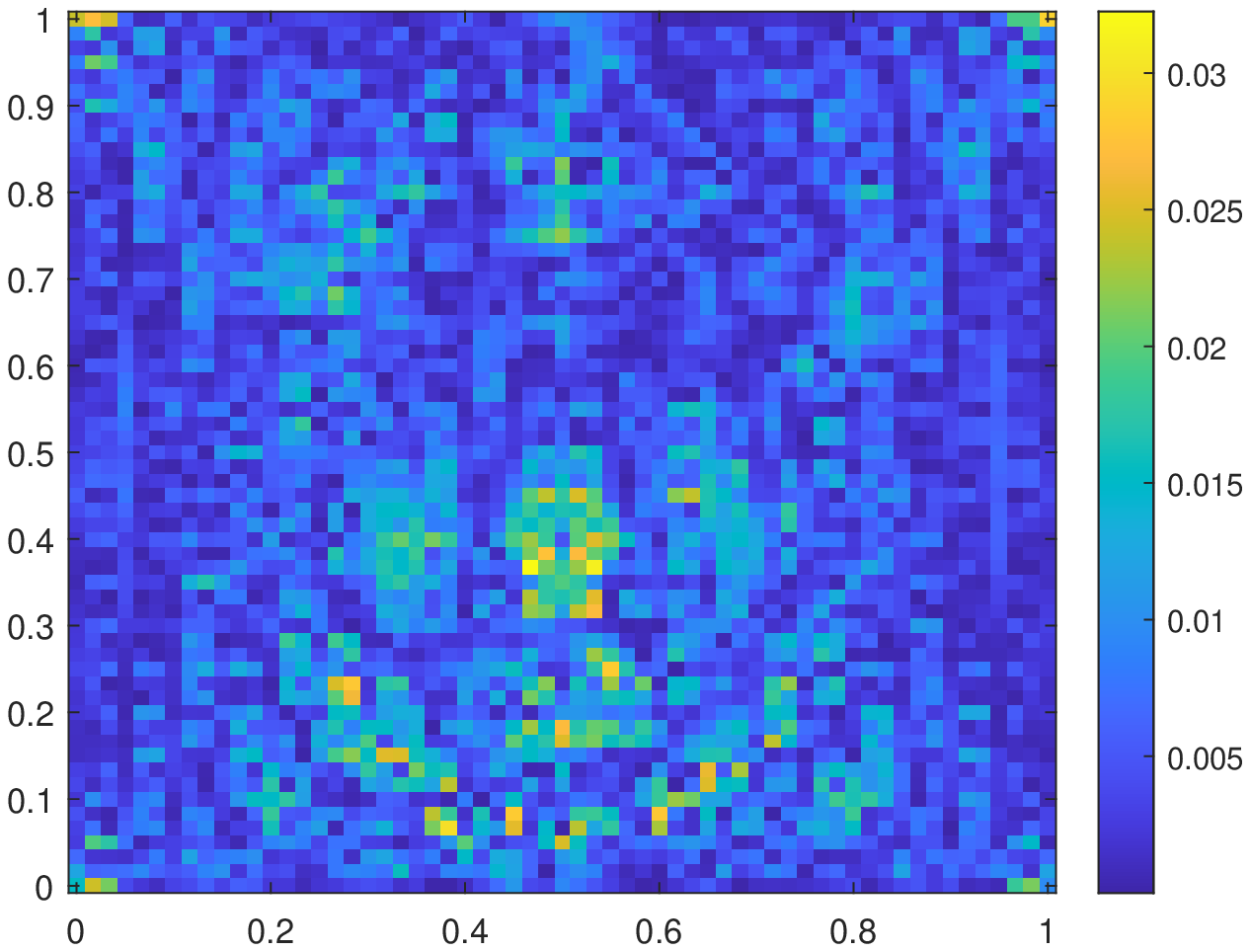}
\caption{Reconstructed $S_4$ using Fredholm inversion. For the first row, 0\%, 1\%, 2\%, 5\% random noises are added to $H_{v_0}$. The relative $L^2$ errors of the reconstructions are 4.3211\%, 4.3405\%, 4.4136\%, 5.0152\%, respectively. 
The second row displays the corresponding differences between the ground truth and the reconstructions.}
\label{fig:Fredholm3}
\end{figure}

\bigskip
\textbf{Experiment 4: Inversion with Theorem~\ref{thm:Fredholm}: Error Analysis.}
The errors in the reconstruction of the Shepp-Logan phantom is substantial. This is mostly due to our choice of the basis~\eqref{eqn:approx} in the discretization. The basis there consists only of smooth polynomials $x^ix^j$ or pyramid-shaped functions $f_{ij}$, which fail to effectively represent a discontinuous function like the Shepp-Logan phantom. In order to justify this, we filter the phantom with a 2D Gaussian kernel with standard deviation 3 to get a smoother phantom (see Figure~\ref{fig:Fcoef3}) and re-run the experiment. Then the errors in the reconstructions are greatly mitigated, as is illustrated in Figure~\ref{fig:Fredholm3}.
The experiment shows that suitable bases are critical for the success of the Fredholm inversion.

\section{Conclusion}.

In this paper, we studied the ultrasound modulated biolumnescence tomography. Assuming knowledge of the attenuation coefficient $\sigma$, the scattering kernel $k(x,\theta,\theta')$ and the domain $X$, we proved that the isotropic source $S(x)$ can be uniquely and stably reconstructed from the internal data $H_{v_0}$ in Theorem~\ref{thm:Neumann} and Theorem~\ref{thm:Fredholm}. The key step of the Fredholm method is to find a proper basis to reduce the approximation error in \eqref{eqn:approx}. The reconstructive procedures for $S$ are provided and numerically implemented in several experiments, in the presence or absence of noise, to demonstrate the efficiency of the reconstruction.

\appendix

\section{Proof of Lemma 2}

\begin{proof}[Proof of Lemma 2]
Let $T^*v = -\theta \cdot \nabla v + \sigma v$ and $Kv = \int_{S^{n-1}}k(x,\theta,\theta')v(x,\theta') \, d\theta.$
Since $\sigma$ and $k$ are uniformly positive, $K$ and the operator $T^{*-1}$ obtained by solving the transport equation $T^{*}v=w$ are uniformly positive. Moreover the solution if $J^* f_0$ is the solution to the ballistic equation $T^{*}v = 0$ with boundary condition $f_0$, then $J^*f_0$ is uniformly positive. Then it follows from the collision expansion form of the solution $v_0$
\[
v_0 = (I + T^{*-1}K + (T^{*-1}K)^2 + \ldots) J^*f_0 
\]
(see for example equation 2.28 in ~\cite{bal2016ultrasound}) that $v_0$ is uniformly positive and the result follows. 
\end{proof}

\bibliographystyle{plain}
\bibliography{ref}

\end{document}